\documentclass[12pt]{article}
\usepackage[citecolor=blue, colorlinks=true, linktoc=page, pagebackref=true, pdfauthor={Benjamin Otto}, pdftitle={Coalescence under Preimage Constraints (Otto)}]{hyperref} % load first, else wonky behavior with duplicate identifiers for different theorems
\usepackage{WideMargins2}
\usepackage{tikz,colortbl,pgfpages}
\usetikzlibrary{arrows,snakes,shapes}
\usepackage{breakcites}
\newcommand\set[1]{\lfloor #1 \rfloor} % notation for the set {1,2,...,n}
\newcommand\coeff[1]{[ #1 ]} % notation for operator [z^n]
\newcommand{\compinv}{{(-1)}} % notation for the decorator to indicate compositional inverse, what Goulden/Jackson would denote [-1]
\newcommand{\val}{{\rm val}}
\newcommand{\CYC}{{\rm CYC}}
\newcommand{\SEQ}{{\rm SEQ}}
\newcommand{\SET}{{\rm SET}}
\newcommand{\COMP}{{\rm COMP}}
\newcommand{\CONNECT}{{\rm CONNECT}}

\title{\Large\bfseries Coalescence under Preimage Constraints}
\author{Benjamin Otto}
\date{\today}
\begin {document}
\fancytitle

\begin{abstract}
  The primary goal of this document is to record the asymptotic effects that preimage constraints impose upon the sizes of the iterated images of a random function.
  Specifically, given a subset $\mathcal{P}\subseteq \mathbb{Z}_{\geq 0}$ and a finite set $S$ of size $n$, choose a function uniformly from the set of functions $f:S\rightarrow S$ that satisfy the condition that $|f^{-1}(x)|\in\mathcal{P}$ for each $x\in S$, and ask what $|f^k(S)|$ looks like as $n$ goes to infinity.
  The robust theory of singularity analysis allows one to completely answer this question if one accepts that $0\in\mathcal{P}$, that $\mathcal{P}$ contains an element bigger than 1, and that $\gcd(\mathcal{P})=1$; only the third of these conditions is a meaningful restriction.
  The secondary goal of this paper is to record much of the background necessary to achieve the primary goal.
\end{abstract}

\newpage
\tableofcontents
\newpage

% ------------------------ body of paper ------------------------------
\begin{section}{Introduction}
  \label{introduction}

  The primary goal of this document is an asymptotic description of the effects that preimage constraints impose upon the sizes of the iterated images of a random function.
  Specifically, given a subset $\mathcal{P}\subseteq \mathbb{Z}_{\geq 0}$ and a finite set $S$ of size $n$, choose a function uniformly from the set of functions $f:S\rightarrow S$ that satisfy the condition that $|f^{-1}(x)|\in\mathcal{P}$ for each $x\in S$, and ask what $|f^k(S)|$ looks like as $n$ goes to infinity.
  A complete answer to this question is found, contingent upon three conditions: that $0\in\mathcal{P}$, that $\mathcal{P}$ contains an element bigger than 1, and that $\gcd(\mathcal{P})=1$.
  Note that only the third of these conditions is a true restriction: if $0\notin \mathcal{P}$ or $\mathcal{P}\subseteq\{0,1\}$, the constraint $\mathcal{P}$ forces $f$ to be a permutation and $|f^k(S)|=n$ for all $k$.

  The very impatient reader is directed to the answer to this question in Theorem \ref{kth_image}; for their benefit, explicit reference to the notation and equations needed to implement that result are repeated immediately prior the theorem.

  The impatient, but knowledgeable, reader is invited to skim Sections \ref{singularity_analysis} and \ref{coalescence}, which follow through on applying singularity analysis to the preimage constraint problem; there are no surprises here and most of the results appear elsewhere.

  For everyone else, a more detailed outline of the paper appears following the motivation and context in the next section.

  \begin{subsection}{Background}

    This document grew out of the author's notes while tracking down the answer to the following question: given a finite set $S$ and a function $f:S\rightarrow S$, what does $|f^k(S)|$ look like? In other words, how does the set of iterated images shrink over time?

    The standard result is Direct Parameters Theorem 2 in \cite{FlajoletOdlyzko_nolink}, and it is natural to generalize the problem by only considering functions that satisfy a preimage constraint $\mathcal{P}$.
    This was done for a whole class of related problems in \cite{ArneyBender1982}, and the general machinery needed to answer the coalescence question is laid out beautifully in \cite{FlajoletSedgewick}, but this author was unable to find a mindless formula that addresses preimage constraints.
    Hence, the main result below, Theorem \ref{kth_image}, is a very modest application of well-understood techniques.
    As a neophyte working through this topic for the first time, the author stumbled a number of times along the way, and it has been a useful exercise to work through much of the theory from scratch.
    Besides the references listed above, additional insight is available from \cite{GouldenJackson}, \cite{Broder}, \cite{Kolchin}, and \cite{PemantleWilson}, and it is handy to have the general background from \cite{GrahamKnuthPatashnik}, \cite{Wilf_nolink}, \cite{SedgewickFlajolet}, and \cite{vanLintWilson} close at hand.

    This document begins by developing just enough definitions and background to put the necessary algebraic combinatorics on solid ground.
    After defining combinatorial classes in Section \ref{combinatorial_classes} and setting up combinatorial class constructions in Section \ref{constructions}, we pilfer heavily from \cite{GouldenJackson} for the formal algebraic theory necessary for generating functions.
    Section \ref{preimage_combinatorics} then provides an explicit development of many interesting functional constructions in a general environment that allows for arbitrary constraints on the sizes of preimages.
    At this point, the algebraic development is complete, and the stage is set for singularity analysis.
    Section \ref{analytic_background} gives a whirlwind tour of the necessary analytic background, developed far less carefully than the preceding algebraic development and appealing to some proofs in other documents.
    Section \ref{singularity_analysis} puts all of these pieces together, culminating in Theorem \ref{kth_image}, which definitively describes the size of the $k$th image of a function from a set to itself, averaging over all functions subject to a given preimage constraint.
    Section \ref{coalescence} contains a brief empirical exploration of the consequences of the main result.

    The knowledgeable reader is certain to find something to bristle about in this document.
    Perhaps it is the decision to pretend all generating functions are exponential, or the pedantic review of formal power series, or the spotty review of complex analytic techniques, or the decision to downplay \cite{ArneyBender1982}, or the decision to ignore \cite{Kolchin}.
    These are all legitimate complaints.
  \end{subsection}

  \begin{subsection}{Notation and Terminology}
    Sequences below are written 0-up, so $(A_n)$ denotes $A_0,A_1,A_2,\dots$.

    For $n\in\mathbb{Z}$, write $\set{n}$ for the set $\{1,2,\dots, n\}$.

    Given a multivariate indeterminate $\bold{u}=(u_1,u_2,\dots,u_r)$ and $\bold{k}=(k_1,k_2,\dots,k_r)\in(\mathbb{Z}_{n\geq0})^r$, write $\bold{u}^\bold{k}=\prod_{i=1}^r u_i^{k_i}$.
    See Section \ref{power_series}.

    If $f(z)=\sum_{n=0}^\infty f_n z^n$ is some formal power series, write
    $$\coeff{z^n}f(z)=f_n.$$
    See Section \ref{power_series}.

    A \emph{rooted tree} is a finite acyclic graph with a distinguished vertex, called the \emph{root}.
    Choosing a root imposes a direction on the edges via the convention that each edge points toward the root.
    Anticipating the function constructions of Sections \ref{treecursion} and \ref{image_combinatorics}, we say that a node $a$ is a \emph{preimage} of another node $b$ if there is a directed edge from $a$ to $b$ (in other words, if there is an edge between $a$ and $b$ and if $b$ is closer to the root than $a$).
    A node is called a \emph{leaf} if it has no preimages.
  \end{subsection}
\end{section}

\begin{section}{Combinatorial Classes}
  \label{combinatorial_classes}

  An \emph{object} or \emph{configuration} is a generic term for some mathematical object that has a finite set associated to it; the elements of this finite set are called \emph{vertices}.
  The \emph{size} of an object $a$, denoted $|a|$, is the number of vertices.
  For the purposes of this paper, a \emph{combinatorial class} or \emph{class} is a set of objects such that the number of elements of any given size is finite.

  If this vague definition is a little unsettling, consider the difficulty in defining an element of a set.
  The key feature of an element is that it is in a set.
  Analogously, the key features of an object are that it is in a class and that it has a size.
  More intuitively, one should think of objects as graphs (undirected, directed, or with whatever additional structure one might require).
  In fact, it is tempting for the specific topic of random mappings (as in this paper) to specifically define objects to be graphs, since the story begins and ends with directed graphs; however, there are points in the middle of the journey that consider tuples and equivalence classes that do not necessarily have a graph structure.

  The \emph{counting sequence} of a class $\mathfrak{A}$ is the sequence $(A_n)$ where
  $$A_n=|\{a\in\mathfrak{A}\mid |a|=n\}|.$$
  The condition that each $A_n$ be finite is exactly the criterion that a set $\mathfrak{A}$ of objects be a combinatorial class.

  For example, let $\mathfrak{A}$ be the set of all rooted trees, where the vertices are exactly the vertices in the graph sense of the word.
  Then $A_0=0$, since the root ensures there is at least one vertex.
  As in Figure \ref{rooted_trees}, the next few values $(A_n)$ are given by $A_1=1$, $A_2=1$, $A_3=2$, and $A_4=4$.

  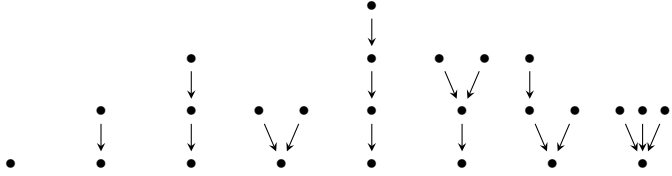
\begin{figure}
    \begin{center}
      \begin{tikzpicture}
        \foreach \h/\w/\n in {.7/1.2/.7}{ % distance between rows, distance between columns, node size
          % vertices (each row is a different graph)
          \foreach \i/\j/\l in {
            -4/0/-9,
            -3/0/-8, -3/1/-7,
            -2/0/-6, -2/1/-5, -2/2/-4,
            -1/0/-3,-1.25/1/-2,-.75/1/-1,
            0/0/1, 0/1/2, 0/2/3, 0/3/4,
            1/0/5, 1/1/6, .75/2/7, 1.25/2/8,
            2/0/9, 1.75/1/10, 2.25/1/11, 1.75/2/12,
            3/0/13, 2.75/1/14, 3/1/15, 3.25/1/16}{
            \draw(\i*\w, \j*\h) node[scale=\n] (node\l) {$\bullet$};
          }

          % edges (each row matches corresponding row above, except first row above has no edges)
          \foreach \l/\m in {
            -7/-8,
            -5/-6, -4/-5,
            -2/-3,-1/-3,
            4/3, 3/2, 2/1,
            7/6, 8/6, 6/5,
            12/10, 10/9, 11/9,
            14/13, 15/13, 16/13}{
            \draw[-stealth] (node\l) -- (node\m);
          }
        }
      \end{tikzpicture}
    \end{center}
    \caption{Rooted Trees of Size at most 4}
    \label{rooted_trees}
  \end{figure}

  An object $a$ is called \emph{weakly labeled} if its vertices are a subset of $\mathbb{Z}_{>0}$. A vertex of a weakly labeled object is also called a \emph{label}. A weakly labeled object $a$ is called \emph{well-labeled} or \emph{labeled} if its set of labels is exactly $\set{|a|}=\{1,2,\dots,|a|\}$. A class is \emph{labeled} if every object in it is labeled.

  Despite the fact that labels are defined to be vertices, one, in practice, treats a label as a feature of a vertex. Intuitively, one starts with a graph and then looks at ways to assign labels to vertices. See Figure \ref{weakly_labeled_example} for a weakly labeled rooted tree. This tree is not well-labeled, but the one in Figure \ref{labeled_example} is.

  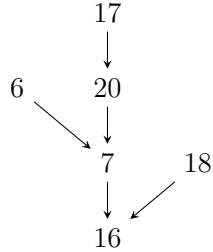
\begin{figure}
    \begin{center}
      \begin{tikzpicture}
        \foreach \h/\w/\n in {1/1.2/.9}{ % distance between rows, distance between columns, node size
          % vertices
          \foreach \i/\j/\l in {3/1/16, 3/2/7, 2/3/6, 3/3/20, 4/2/18, 3/4/17}{
            \draw (\i*\w, \j*\h) node[scale=\n] (node\l) {\l};
          }

          % edges
          \foreach \l/\m in {7/16, 6/7, 20/7, 18/16, 17/20}{
            \draw[-stealth] (node\l) -- (node\m);
          }

        }
      \end{tikzpicture}
    \end{center}
    \caption{Weakly Labeled Rooted Tree that is Not Well-Labeled}
    \label{weakly_labeled_example}
  \end{figure}

  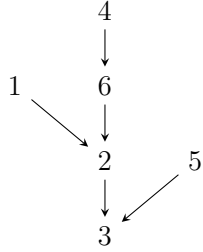
\begin{figure}
    \begin{center}
      \begin{tikzpicture}
        \foreach \h/\w/\n in {1/1.2/.9}{ % distance between rows, distance between columns, node size
          % vertices
          \foreach \i/\j/\l in {3/1/3, 3/2/2, 2/3/1, 3/3/6, 4/2/5, 3/4/4}{
            \draw (\i*\w, \j*\h) node[scale=\n] (node\l) {\l};
          }

          % edges
          \foreach \l/\m in {4/6, 1/2, 6/2, 2/3, 5/3}{
            \draw[-stealth] (node\l) -- (node\m);
          }

        }
      \end{tikzpicture}
    \end{center}
    \caption{Well-Labeled Rooted Tree}
    \label{labeled_example}
  \end{figure}

  Building on the earlier example of rooted trees $\mathfrak{A}$, let $\mathfrak{B}$ be the set of all labeled rooted trees (see, for example, Section 3.3.8 of \cite{GouldenJackson}); in other words, isomorphic rooted trees will be treated as distinct if there is no graph isomorphism that preserves the labels.
  One may check that $B_0=0$, $B_1=1$, $B_2=2$, and $B_3=9$; see Figure \ref{labeled_rooted_trees}.

  \begin{figure}
    \begin{center}
      \begin{tikzpicture}
        \foreach \h/\w/\n in {1/1.2/.9}{ % distance between rows, distance between columns, node size
          % size 1
          \draw(0*\w, 0*\h) node[scale=\n] (node1) {1};

          % size 2
          \foreach \i/\a/\b in {1/1/2, 2/2/1}{
            \draw(\i*\w, 0*\h) node[scale=\n] (node\a) {\a};
            \draw(\i*\w, 1*\h) node[scale=\n] (node\b) {\b};
            \draw[-stealth] (node\b) -- (node\a);
          }

          % size 3
          \foreach \i/\a/\b/\c in {3/1/2/3, 4/1/3/2, 5/2/1/3, 6/2/3/1, 7/3/1/2, 8/3/2/1}{
            \draw(\i*\w, 0*\h) node[scale=\n] (node\a) {\a};
            \draw(\i*\w, 1*\h) node[scale=\n] (node\b) {\b};
            \draw(\i*\w, 2*\h) node[scale=\n] (node\c) {\c};
            \draw[-stealth] (node\b) -- (node\a);
            \draw[-stealth] (node\c) -- (node\b);
          }
          \foreach \i/\a/\b/\c in {9/1/2/3, 10/2/1/3, 11/3/1/2}{
            \draw(\i*\w, 0*\h) node[scale=\n] (node\a) {\a};
            \draw(\i*\w-.25*\w, 1*\h) node[scale=\n] (node\b) {\b};
            \draw(\i*\w+.25*\w, 1*\h) node[scale=\n] (node\c) {\c};
            \draw[-stealth] (node\b) -- (node\a);
            \draw[-stealth] (node\c) -- (node\a);
          }
        }
      \end{tikzpicture}
    \end{center}
    \caption{Labeled Rooted Trees of Size up to 3}
    \label{labeled_rooted_trees}
  \end{figure}
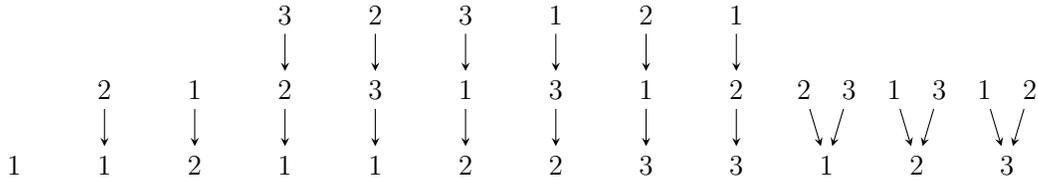

  When we wish to analyze additional structure in a class $\mathfrak{A}$ beyond the number of vertices, we employ a \emph{parameter} function $\chi:\mathfrak{A}\rightarrow (\mathbb{Z}_{\geq0})^r$.
  In explicit cases where $r=1$, we will drop the tuple notation, so that $\chi(a)\in\mathbb{Z}$.

  By allowing $\chi$ to send every object to the trivial 0-tuple, we may allow for the case $r=0$ and so we could have built the existence of a parameter function into the definition of combinatorial class.
  Alternatively, some sources view the size function as just one of the coordinates in the parameter function; here, we wish to explicitly preserve the distinction between size and the other parameters.

  An example of a parameter function is to define $\chi:\mathfrak{B}\rightarrow (\mathbb{Z}_{\geq 0})^2$ on the class of labeled rooted trees by $\chi(b)=(l_b,r_b)$, where $l_b$ is the number of leaves in $b$ and $r_b$ is the number of roots in $b$.
  Of course, this definition is a little silly, since $r_b=1$ for all $b\in\mathfrak{B}$.
\end{section}

\begin{section}{Constructions}
  \label{constructions}

  We need to understand some standard ways to construct a new labeled class from old labeled classes.

  \begin{subsection}{Combinatorial Sum}

    The \emph{disjoint union} or \emph{combinatorial sum} of two disjoint labeled classes $\mathfrak{B}$ and $\mathfrak{C}$ is the union $\mathfrak{B}+\mathfrak{C}=\mathfrak{B}\cup\mathfrak{C}$.
    It is apparent that the finiteness condition for $\mathfrak{B}+\mathfrak{C}$ to be a combinatorial class holds.
    Every object is well-labeled, since it is in one of the labeled classes $\mathfrak{B}$ or $\mathfrak{C}$.
    Thus, $\mathfrak{B}+\mathfrak{C}$ is again a labeled class.

    If $\mathfrak{B}$ and $\mathfrak{C}$ have respective parameter functions $\psi$ and $\omega$, both of the same dimension, then the sum $\mathfrak{B}+\mathfrak{C}$ inherits a parameter function $\chi$ via
    $$\chi(\alpha)=\begin{cases}
    \psi(\alpha) & \text{if } \alpha\in\mathfrak{B}\\
    \omega(\alpha) & \text{if } \alpha\in\mathfrak{C}
    \end{cases}.$$
  \end{subsection}

  \begin{subsection}{Labeled Product}
    Next, we turn to defining a product of labeled classes.
    Intuitively, we want something like a direct product, but the labelings introduce a wrinkle.
    Namely, we need to ensure that the objects in the new class cannot have repeated labels and that they are well-labeled, not just weakly labeled.
    To address this, we first need to codify some rules for relabeling an object.

    The \emph{reduction} of a weakly labeled object is the well-labeled object resulting from the process of ordering the original labels and then replacing each label with its position in the ordering.
    In other words, we are replacing the labels with the unique well-labeling that respects order.
    For example, the rooted tree in Figure \ref{labeled_example} is the reduction of the rooted tree in Figure \ref{weakly_labeled_example}.

    Going in the opposite direction, there are many ways to take a well-labeled object to a weakly labeled one.
    Given any increasing function $f:\set{n}\rightarrow\mathbb{Z}$, the \emph{expansion} with respect to $f$ of a well-labeled object of size $n$ is the weakly labeled object achieved by replacing each label $i\in\set{n}$ with $f(i)$.
    For example, the rooted tree in Figure \ref{weakly_labeled_example} is the expansion with respect to $f$ of the rooted tree in Figure \ref{labeled_example}, where $f(1)=6$, $f(2)=7$, $f(3)=16$, $f(4)=17$, $f(5)=18$, $f(6)=20$; note that, as required, this $f$ is a strictly increasing function.

    Finally, note the ordered pair $(\beta,\gamma)$ of weakly labeled objects $\beta$ and $\gamma$ that do not share labels is again a weakly labeled object with vertex set given by the union of $\beta$'s and $\gamma$'s vertex sets.
    Thus, $|(\beta,\gamma)|=|\beta|+|\gamma|$.

    Define the \emph{labeled product} of two labeled objects $\beta$ and $\gamma$ to be the set
    $$\beta\star\gamma=\{(\beta',\gamma') \mid (\beta',\gamma') \text{ is well-labeled}, \rho(\beta')=\beta, \rho(\gamma')=\gamma\},$$
    where $\rho$ denotes reduction.
    In other words, an element of $\beta\star\gamma$ is a pair $(\beta',\gamma')$ where there is no overlap between labels and, taken together, the labels are $\set{|\beta|+|\gamma|}$.
    Because an expansion function is uniquely determined by its image, and because the label set of any element in $\beta\star\gamma$ is $\set{|\beta|+|\gamma|}$, the size of $\beta\star\gamma$ is exactly $\binom{|\beta|+|\gamma|}{|\beta|}$.
    In particular, $\beta\star\gamma$ is finite, and the finiteness condition for $\beta\star\gamma$ to be a combinatorial class holds.
    Finally, every element is well-labeled by definition, so $\beta\star\gamma$ is a labeled class.

    See Figure \ref{labeled_product} for an example of a labeled product of two labeled objects that were taken from Figure \ref{labeled_rooted_trees}.

    \begin{figure}
      \begin{center}
        \begin{tikzpicture}
          \foreach \h/\w/\n in {.7/1.2/.8}{ % distance between rows, distance between columns, node size
            \draw(0*\w, 0*\h) node[scale=\n] {$\beta=$};
            \draw(.5*\w, .5*\h) node[scale=\n] (node2) {2};
            \draw(.5*\w, -.5*\h) node[scale=\n] (node1) {1};
            \draw[-stealth] (node2) -- (node1);

            \draw(4*\w, 0*\h) node[scale=\n] {$\gamma=$};
            \draw(4.25*\w, .5*\h) node[scale=\n] (node3) {1};
            \draw(4.75*\w, .5*\h) node[scale=\n] (node1) {3};
            \draw(4.5*\w, -.5*\h) node[scale=\n] (node2) {2};
            \draw[-stealth] (node1) -- (node2);
            \draw[-stealth] (node3) -- (node2);

            \draw(-2.2*\w, -3*\h) node[scale=\n] {$\beta\star\gamma=$};
            \draw(-1.5*\w, -3*\h) node[scale=2.7*\n] {$\{$};

            \foreach \i/\j/\a/\b/\c/\d/\e in {0/0/1/2/3/4/5, 1/0/1/3/2/4/5, 0/1/1/4/2/3/5, 1/1/1/5/2/3/4, 0/2/2/3/1/4/5, 1/2/2/4/1/3/5, 0/3/2/5/1/3/4, 1/3/3/4/1/2/5, 0/4/3/5/1/2/4, 1/4/4/5/1/2/3}{
              \draw(-1.25*\w+\i*4*\w, -3*\h-\j*2*\h) node[scale=2*\n]  {$($};
              \draw(-1*\w+\i*4*\w, -3.5*\h-\j*2*\h) node[scale=\n] (node\a) {\a};
              \draw(-1*\w+\i*4*\w, -2.5*\h-\j*2*\h) node[scale=\n] (node\b) {\b};
              \draw(-.8*\w+\i*4*\w, -3*\h-\j*2*\h) node[scale=2*\n]  {$,$};
              \draw(0*\w+\i*4*\w, -3.5*\h-\j*2*\h) node[scale=\n] (node\d) {\d};
              \draw(-.25*\w+\i*4*\w, -2.5*\h-\j*2*\h) node[scale=\n] (node\c) {\c};
              \draw(.25*\w+\i*4*\w, -2.5*\h-\j*2*\h) node[scale=\n] (node\e) {\e};
              \draw(.5*\w+\i*4*\w, -3*\h-\j*2*\h) node[scale=2*\n]  {$)$};
              \draw[-stealth] (node\b) -- (node\a);
              \draw[-stealth] (node\c) -- (node\d);
              \draw[-stealth] (node\e) -- (node\d);
            }

            \foreach \i/\j in {0/0, 1/0, 0/1, 1/1, 0/2, 1/2, 0/3, 1/3, 0/4}{
              \draw(1*\w+\i*4*\w, -3*\h-\j*2*\h) node[scale=2*\n] {$,$};
            }
            \draw(1.5*\w+1*4*\w, -3*\h-4*2*\h) node[scale=2.7*\n] {$\}$};
          }
        \end{tikzpicture}
      \end{center}
      \caption{An Example of $\beta\star\gamma$}
      \label{labeled_product}
    \end{figure}

    The \emph{labeled product} of two labeled classes $\mathfrak{B}$ and $\mathfrak{C}$ is the set $\mathfrak{B}\star\mathfrak{C}=\cup_{\beta\in B,\gamma\in C}\beta\star\gamma$.
    In order to ensure that $\mathfrak{B}\star\mathfrak{C}$ is a combinatorial class, we must check that, for any given size $n$, there are a finite number of objects of that size.
    Elements of $\mathfrak{B}\star\mathfrak{C}$ of size $n$ are exactly the elements of sets that look like $\beta\star\gamma$ where $|\beta|+|\gamma|=n$.
    Since there are only finitely many choices for $0\leq i\leq n$, and the fact that $\mathfrak{B}$ and $\mathfrak{C}$ are combinatorial classes ensures that there are only finitely many choices for $\beta$ of size $i$ and $\gamma$ of size $n-i$, and the size of each $\beta\star\gamma$ is finite, we are done.
    Moreover, since every element is well-labeled, the labeled product $\mathfrak{B}\star\mathfrak{C}$ is a labeled class.

    If $\mathfrak{B}$ and $\mathfrak{C}$ have respective parameter functions $\psi$ and $\omega$, both of the same dimension, then the labeled product $\mathfrak{B}\star\mathfrak{C}$ inherits a parameter function $\chi$ via
    $$\chi(\beta',\gamma')=\psi\big(\rho(\beta')\big)+\omega\big(\rho(\gamma')\big),$$
    where $\rho$ denotes reduction of a weakly labeled object.

    As presented here, $\star$ is not strictly speaking associative, since the order of execution will affect how one writes the nested tuples.
    For the applications in this paper, however, this is a silly distinction, and there is no harm in writing $\mathfrak{A}\star\mathfrak{B}\star\mathfrak{C}$.
    Similarly, $\star$ is not formally commutative, and it is tempting to say one does not care what order one writes a pair.
    In this case, though, there is a danger in being too cavalier.
    Namely, while one can choose, once and for all, the order of execution, this does not allow one to identify the results of executing $\star$ in different orders.
    For an example of why this is important, compare the constructions in Sections \ref{sequence} and \ref{cycle}.
  \end{subsection}

  \begin{subsection}{Sequence}
    \label{sequence}

    Write $\SEQ_0(\mathfrak{B})$ for the labeled class whose only object is the empty set with no vertices.
    For each $k>0$, define $\SEQ_k(\mathfrak{B})=\mathfrak{B}\star \SEQ_{k-1}(\mathfrak{B})$; recall there is no harm here in writing $\SEQ_k(\mathfrak{B})=\mathfrak{B}\star\mathfrak{B}\star\dots\star\mathfrak{B}$.
    It is immediate from the discussion of labeled products that $\SEQ_k(\mathfrak{B})$ is a labeled class.

    Let $\mathfrak{B}$ be a labeled class that does not have any objects of size $0$.
    The \emph{labeled sequence} or \emph{sequence} or \emph{labeled power} class of $\mathfrak{B}$ is $\SEQ(\mathfrak{B})=\cup_{k\geq 0}\SEQ_k(\mathfrak{B})$.
    Because there are no elements of size 0, each element of $\SEQ_k(\mathfrak{B})$ has size at least $k$.
    Thus, when checking the finiteness condition to see if $\SEQ(\mathfrak{B})$ is a combinatorial class for a given size $n$, it is enough just to consider the finite union $\cup_{0\leq k\leq n}\SEQ_k(\mathfrak{B})$ of labeled classes, an environment in which the condition clearly holds.
    Thus, $\SEQ(\mathfrak{B})$ is a labeled class.

    By way of a $\SEQ$ example, let $a$ denote the labeled rooted tree of size two whose root is 2.
    Write $\mathfrak{A}=\{a\}$.
    See Figure \ref{seq2_example} for a depiction of $\SEQ_2(\mathfrak{A})$.

    \begin{figure}
      \begin{center}
        \begin{tikzpicture}
          \foreach \h/\w/\n in {.7/1.2/.8}{ % distance between rows, distance between columns, node size
            \draw(0*\w, 0*\h) node[scale=\n] {$a=$};
            \draw(.5*\w, .5*\h) node[scale=\n] (node1) {1};
            \draw(.5*\w, -.5*\h) node[scale=\n] (node2) {2};
            \draw[-stealth] (node1) -- (node2);

            \draw(3*\w, 0*\h) node[scale=\n] {$\mathfrak{A}=\{a\}$};

            \draw(-2.7*\w, -3*\h) node[scale=\n] {$\SEQ_2(\mathfrak{A})=$};
            \draw(-1.5*\w, -3*\h) node[scale=2.7*\n] {$\{$};

            \foreach \i/\j/\a/\b/\c/\d in {0/0/2/1/3/4, 1/0/4/3/1/2, 0/1/3/1/2/4, 1/1/4/2/1/3, 0/2/4/1/2/3, 1/2/3/2/1/4}{
              \draw(-1.1*\w+\i*4*\w, -3*\h-\j*2*\h) node[scale=2*\n]  {$($};
              \draw(-.8*\w+\i*4*\w, -3.5*\h-\j*2*\h) node[scale=\n] (node\a) {\a};
              \draw(-.8*\w+\i*4*\w, -2.5*\h-\j*2*\h) node[scale=\n] (node\b) {\b};
              \draw(-.5*\w+\i*4*\w, -3*\h-\j*2*\h) node[scale=2*\n]  {$,$};
              \draw(-.2*\w+\i*4*\w, -3.5*\h-\j*2*\h) node[scale=\n] (node\d) {\d};
              \draw(-.2*\w+\i*4*\w, -2.5*\h-\j*2*\h) node[scale=\n] (node\c) {\c};
              \draw(.1*\w+\i*4*\w, -3*\h-\j*2*\h) node[scale=2*\n]  {$)$};
              \draw[-stealth] (node\b) -- (node\a);
              \draw[-stealth] (node\c) -- (node\d);
            }

            \foreach \i/\j in {0/0, 1/0, 0/1, 1/1, 0/2}{
              \draw(.4*\w+\i*4*\w, -3*\h-\j*2*\h) node[scale=2*\n] {$,$};
            }
            \draw(1.5*\w+1*4*\w, -3*\h-2*2*\h) node[scale=2.7*\n] {$\}$};
          }
        \end{tikzpicture}
      \end{center}
      \caption{An Example of $\SEQ_2(\mathfrak{A})$}
      \label{seq2_example}
    \end{figure}

    If $\mathfrak{B}$ has parameter function $\chi$, assign $\SEQ_0(\mathfrak{B})$ a parameter function of the same dimension as $\chi$ by sending the empty sequence to the 0-vector.
    For $k>0$, $\SEQ_k(\mathfrak{B})=\mathfrak{B}\star \SEQ_{k-1}(\mathfrak{B})$ inherits a parameter function from the $\star$-construction.
    Finally, since each $\SEQ_k(\mathfrak{B})$ is disjoint from every other, $\SEQ(\mathfrak{B})$ inherits a parameter function by just evaluating an element with the parameter function of the unique $\SEQ_k(\mathfrak{B})$ that it lies in; this is just an infinite version of the inherited parameter in the combinatorial sum construction, except that the dimensions of the underlying parameter function line up for free.

    Note that the condition that $\mathfrak{B}$ have no element of size 0 is necessary for the construction of $\SEQ$, but not for $\SEQ_k$.
    Specifically, if $\mathfrak{B}$ does have an element $b$ of size 0, then taking the $\star$ of $b$ with itself any number of times always yields a sequence of size 0, so $\SEQ(\mathfrak{B})$, with its infinitely many elements of size 0, is not a combinatorial class.
    There is no such issue for $\SEQ_k$, which is just an iterated product.
  \end{subsection}

  \begin{subsection}{Cycle}
    \label{cycle}

    If $\mathfrak{B}$ is a labeled class that does not have any objects of size $0$ and if $k>0$, then $\CYC_k(\mathfrak{B})$ denotes the class whose objects are equivalence classes in $\SEQ_k(\mathfrak{B})$ given by identifying those sequences that are cyclic shifts of each other.
    Since each element in such an equivalence class has the same size, and hence the same vertex set, the vertex set of the equivalence class is defined unambiguously.
    The number of elements of a given size in $\CYC_k(\mathfrak{B})$ is bounded above by the number of elements of that size in $\SEQ_k(\mathfrak{B})$.
    These two observations are enough to ensure that $\CYC_k(\mathfrak{B})$ is a labeled class.
    Finally, the \emph{cycle} class of $\mathfrak{B}$ is $\CYC(\mathfrak{B})=\cup_{k\geq 1}\CYC_k(\mathfrak{B})$; it is a labeled class by the same argument that $\CYC_k(\mathfrak{B})$ is.

    By way of a $\CYC$ example, let $a$ denote the labeled rooted tree of size two whose root is 2.
    Write $\mathfrak{A}=\{a\}$.
    One gets $\CYC_2(\mathfrak{A})$ from $\SEQ_2(\mathfrak{A})$ by identifying cyclic shifts of an object.
    Since there is no longer any first element, it is handy to represent this by placing the elements of a representative sequence on an actual cycle.
    Since the roots of the trees on the cycle are in bijective correspondence with the trees themselves, we may draw the cycle as permuting the roots of the trees.
    By doing so we may view $\CYC_2(\mathfrak{A})$ as the graphs in Figure \ref{cyc2_example}; note that the two sequence objects in each row of Figure \ref{seq2_example} are considered equivalent, yielding the single cycle object in the corresponding row of Figure \ref{cyc2_example}.
    For another example, an element of $\CYC_3(\mathfrak{A})$ is shown in Figure \ref{cyc3_example}.

    \begin{figure}
      \begin{center}
        \begin{tikzpicture}
          \foreach \h/\w/\n in {.7/1.2/.8}{ % distance between rows, distance between columns, node size
            \draw(0*\w, 0*\h) node[scale=\n] {$a=$};
            \draw(.5*\w, .5*\h) node[scale=\n] (node1) {1};
            \draw(.5*\w, -.5*\h) node[scale=\n] (node2) {2};
            \draw[-stealth] (node1) -- (node2);

            \draw(3*\w, 0*\h) node[scale=\n] {$\mathfrak{A}=\{a\}$};

            \draw(-1.4*\w, -3*\h) node[scale=\n] {$\CYC_2(\mathfrak{A})=$};
            \draw(-.5*\w, -3*\h) node[scale=2.7*\n] {$\{$};

            \foreach \j/\a/\b/\c/\d in {0/1/2/3/4, 1/1/3/2/4, 2/1/4/2/3}{
              \draw(0*\w, -3*\h-\j*2*\h) node[scale=\n] (node\a) {\a};
              \draw(1*\w, -3*\h-\j*2*\h) node[scale=\n] (node\b) {\b};
              \draw(2*\w, -3*\h-\j*2*\h) node[scale=\n] (node\d) {\d};
              \draw(3*\w, -3*\h-\j*2*\h) node[scale=\n] (node\c) {\c};
              \draw[-stealth] (node\a) -- (node\b);
              \draw[-stealth] (node\b.north east) .. controls +(0,.25) and +(0,.25) .. (node\d.north west);
              \draw[-stealth] (node\d.south west) .. controls +(0,-.25) and +(0,-.25) .. (node\b.south east);
              \draw[-stealth] (node\c) -- (node\d);
            }

            \foreach \j in {0, 1}{
              \draw(3.25*\w, -3*\h-\j*2*\h) node[scale=2*\n] {$,$};
            }
            \draw(3.5*\w, -3*\h-2*2*\h) node[scale=2.7*\n] {$\}$};
          }
        \end{tikzpicture}
      \end{center}
      \caption{An Example of $\CYC_2(\mathfrak{A})$}
      \label{cyc2_example}
    \end{figure}
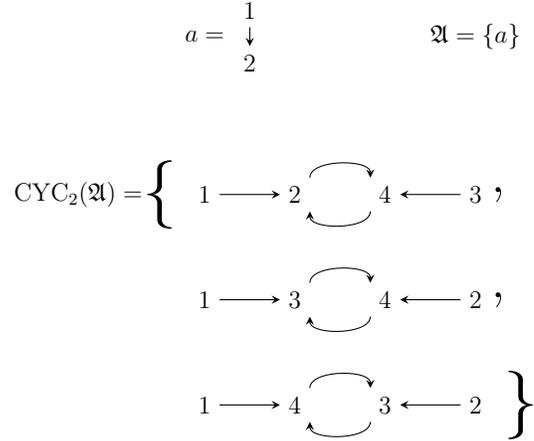

    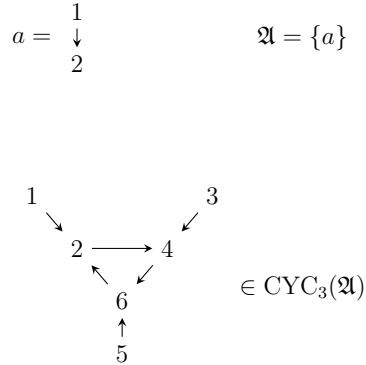
\begin{figure}
      \begin{center}
        \begin{tikzpicture}
          \foreach \h/\w/\n in {.7/1.2/.8}{ % distance between rows, distance between columns, node size
            \draw(0*\w, 0*\h) node[scale=\n] {$a=$};
            \draw(.5*\w, .5*\h) node[scale=\n] (node1) {1};
            \draw(.5*\w, -.5*\h) node[scale=\n] (node2) {2};
            \draw[-stealth] (node1) -- (node2);

            \draw(3*\w, 0*\h) node[scale=\n] {$\mathfrak{A}=\{a\}$};

            \draw(0*\w, -3*\h) node[scale=\n] (node1) {1};
            \draw(.5*\w, -4*\h) node[scale=\n] (node2) {2};
            \draw(1.5*\w, -4*\h) node[scale=\n] (node4) {4};
            \draw(2*\w, -3*\h) node[scale=\n] (node3) {3};
            \draw(1*\w, -5*\h) node[scale=\n] (node6) {6};
            \draw(1*\w, -6*\h) node[scale=\n] (node5) {5};

            \draw[-stealth] (node1) -- (node2);
            \draw[-stealth] (node2) -- (node4);
            \draw[-stealth] (node3) -- (node4);
            \draw[-stealth] (node4) -- (node6);
            \draw[-stealth] (node5) -- (node6);
            \draw[-stealth] (node6) -- (node2);
            \draw(3*\w, -4.75*\h) node[scale=\n] {$\in \CYC_3(\mathfrak{A})$};
          }
        \end{tikzpicture}
      \end{center}
      \caption{An Example of a $\CYC_3(\mathfrak{A})$ Element}
      \label{cyc3_example}
    \end{figure}

    If $\mathfrak{B}$ has parameter function $\chi$, also use $\chi$ to denote the parameter function that $\SEQ(\mathfrak{B})$ and each $\SEQ_k(\mathfrak{B})$ inherits.
    Then $\chi$ is constant on an equivalence class given by shifting a sequence (it is always the sum of $\chi$ applied to the components), and so $\chi$ also acts as a parameter function on $\CYC(\mathfrak{B})$ and $\CYC_k(\mathfrak{B})$.

    Note that the condition $k>0$ was unnecessary for any of the definitions in this section, so it is tempting to define $\CYC_0(\mathfrak{B})$ to be the set whose only element is the equivalence class $\SEQ_0(\mathfrak{B})$.
    We refrain from doing this for two reasons.
    First, the generating function for this definition does not fit the formula that all other $\CYC_k(\mathfrak{B})$ generating functions satisfy.
    Second, if we define $\CYC_0$, we are tempted to include it in the definition of $\CYC$; doing so, however, needlessly complicates the generating function of $\CYC(\mathfrak{B})$.
    Both points speak to the fact that it is not ``right'' to define a cycle with 0 terms.
    See Lemma \ref{gen_fn_constructions} for the correct formulas.

    Similarly, the condition that $\mathfrak{B}$ not have any objects of size $0$ is not necessary to define $\CYC_k(\mathfrak{B})$; however, doing so needlessly complicates the formula for the generating function described in Section \ref{construction_gen_fn}, and it is easier just to define this situation away, as the extra generality is not useful for our applications.
  \end{subsection}

  \begin{subsection}{Set}

    If $\mathfrak{B}$ is a labeled class that does not have any objects of size $0$, the class $\SET_k(\mathfrak{B})$ denotes the class whose objects are equivalence classes in $\SEQ_k(\mathfrak{B})$ given by identifying sequences that are permutations of each other.
    For $\mathcal{P}\subseteq\mathbb{Z}_{\geq0}$, define $\SET_\mathcal{P}(\mathfrak{B})=\cup_{k\in\mathcal{P}}\SET_k(\mathfrak{B})$.
    The \emph{set} class of $\mathfrak{B}$ is $\SET(\mathfrak{B})=\SET_{\mathbb{Z}_{\geq0}}(\mathfrak{B})$.

    Note that, when $|\mathcal{P}|$ is finite, the definition of $\SET_\mathcal{P}(\mathfrak{B})$ could be extended to allow $\mathfrak{B}$ to contain objects of size 0; however, like with $\CYC_k$, it needlessly complicates the formula (and intuition) for the generating function.

    See Figure \ref{set2_example} for an example $\SET_2(\mathfrak{A})$, where $\mathfrak{A}$ consists of the two specified graphs and $\SET_2(\mathfrak{A})$ can be portrayed as graphs by the simple expedient of taking the union of the input components.

    \begin{figure}
      \begin{center}
        \begin{tikzpicture}
          \foreach \h/\w/\n in {.6/1.2/.7}{ % distance between rows, distance between columns, node size

            \draw(-1.9*\w, 0*\h) node[scale=\n] {$\mathfrak{A}=$};
            \draw(-1*\w, 0*\h) node[scale=2.7*\n] {$\{$};
            \foreach \i/\j/\a in {0/0/1}{
              \draw(-.25*\w+\i*4*\w, 0*\h-\j*2*\h) node[scale=\n] (node\a) {\a};
            }
            \foreach \i/\j/\a/\b in {1/0/1/2}{
              \draw(-.25*\w+\i*4*\w, 0*\h-\j*2*\h) node[scale=\n] (node\b) {\b};
              \draw(-.25*\w+\i*4*\w, 1*\h-\j*2*\h) node[scale=\n] (node\a) {\a};
              \draw[-stealth] (node\a) -- (node\b);
            }
            \foreach \i/\j in {0/0}{
              \draw(1*\w+\i*4*\w, 0*\h-\j*2*\h) node[scale=2*\n] {$,$};
            }
            \draw(1.5*\w+1*4*\w, 0*\h-0*2*\h) node[scale=2.7*\n] {$\}$};

            \draw(-1.9*\w, -3*\h) node[scale=\n] {$\SET_2(\mathfrak{A})=$};
            \draw(-1*\w, -3*\h) node[scale=2.7*\n] {$\{$};
            \foreach \i/\j/\a/\b in {0/0/1/2}{
              \draw(-.5*\w+\i*4*\w, -3*\h-\j*2*\h) node[scale=\n] (node\a) {\a};
              \draw(0*\w+\i*4*\w, -3*\h-\j*2*\h) node[scale=\n] (node\b) {\b};
            }
            \foreach \i/\j/\a/\b/\c in {1/0/1/2/3, 0/1/2/1/3, 1/1/3/1/2}{
              \draw(-.5*\w+\i*4*\w, -3*\h-\j*2*\h) node[scale=\n] (node\a) {\a};
              \draw(0*\w+\i*4*\w, -3.5*\h-\j*2*\h) node[scale=\n] (node\c) {\c};
              \draw(0*\w+\i*4*\w, -2.5*\h-\j*2*\h) node[scale=\n] (node\b) {\b};
              \draw[-stealth] (node\b) -- (node\c);
            }
            \foreach \i/\j/\a/\b/\c/\d in {0/2/1/2/3/4, 1/2/1/3/2/4, 0/3/1/4/2/3}{
              \draw(-.5*\w+\i*4*\w, -3.5*\h-\j*2*\h) node[scale=\n] (node\a) {\a};
              \draw(-.5*\w+\i*4*\w, -2.5*\h-\j*2*\h) node[scale=\n] (node\b) {\b};
              \draw(0*\w+\i*4*\w, -3.5*\h-\j*2*\h) node[scale=\n] (node\d) {\d};
              \draw(0*\w+\i*4*\w, -2.5*\h-\j*2*\h) node[scale=\n] (node\c) {\c};
              \draw[-stealth] (node\b) -- (node\a);
              \draw[-stealth] (node\c) -- (node\d);
            }

            \foreach \i/\j in {0/0, 1/0, 0/1, 1/1, 0/2, 1/2}{
              \draw(1*\w+\i*4*\w, -3*\h-\j*2*\h) node[scale=2*\n] {$,$};
            }
            \draw(1.5*\w+1*4*\w, -3*\h-3*2*\h) node[scale=2.7*\n] {$\}$};
          }
        \end{tikzpicture}
      \end{center}
      \caption{An Example of $\SET_2(\mathfrak{A})$}
      \label{set2_example}
    \end{figure}
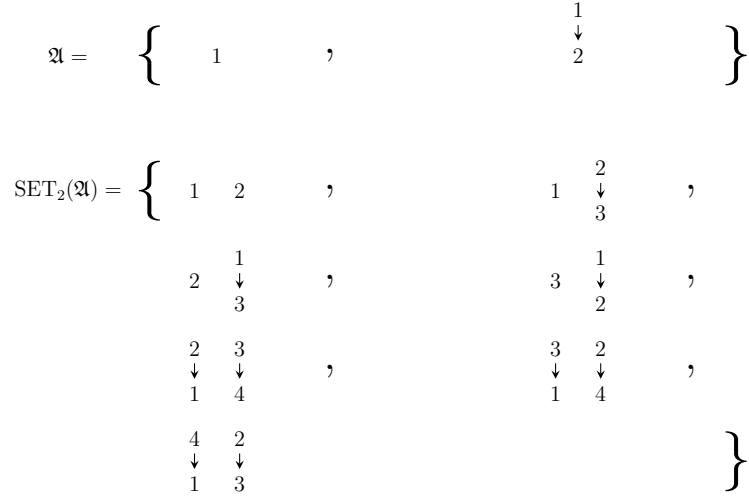
  \end{subsection}
\end{section}

\begin{section}{Generating Functions}

  We now turn to defining generating functions and translating the constructions of the previous section to this context.
  First we need to recall some facts from the theory of formal power series.

  \begin{subsection}{Formal Power Series}
    \label{power_series}

    This section, which comprises a ``just enough'' review of formal power series, was heavily influenced by the first chapter of \cite{GouldenJackson}.
    All of the results of this section will ultimately be applied to the case where the ring $R$ is $\mathbb{C}$, where standard results from complex analysis would suffice.
    Thus, if there is any section in this document that the reader should skip, it is this one.
    On the other hand, it is just plain awesome that one can get Lagrange Inversion from a purely algebraic context.

    We initially take $R$ to be a commutative ring.
    Some of the power series of particular importance will further require that $R\supseteq\mathbb{Q}$.
    In practice, the ring $R$ will be $\mathbb{C}$ or some (possibly multivariate) polynomial ring over $\mathbb{C}$ or a ring of (possibly multivariate) formal power series over $\mathbb{C}$.

    Let $\bold{u}=(u_1,u_2,\dots,u_r)$ be a multivariate indeterminate.
    The ring of formal power series $R[[\bold{u}]]$ is the set $\{\sum_{\bold{k}\in(\mathbb{Z}_{\geq 0})^r}a_\bold{k}\bold{u}^\bold{k} \mid a_\bold{k}\in R\}$ with 0 denoted $\bold{0}$.
    Addition and multiplication are the well-defined extension of the usual polynomial operations, and one can take well-defined infinite sums or products whenever any given monomial in the putative result only has a finite number of terms contributing to it; in particular, as described in Lemma \ref{comp_seq}, there are many important cases where the composition of power series is well-defined.

    For $0\leq s\leq r$, we completely ignore the distinction between $R[[(u_1,\dots, u_s)]]$ and its natural embedding into $R[[\bold{u}]]$.
    In particular, in the exposition below, every power series is over the $r$ variables $\bold{u}$; if this were not the case, one could add dummy variables as necessary.

    However, anticipating the distinction between $\bold{u}$ and $z$ in Section \ref{gen_fn_dfn} and anticipating the restriction to complex functions in Section 6, we will usually write $z$ for a univariate indeterminate.

    The multiplicative inverse of a power series does not always exist, but there is an easy characterization of when it does.

    \begin{lemma}
      \label{inv_seq}
      Let $R$ be a commutative ring.
      Then $f(\bold{u})\in R[[\bold{u}]]$ is invertible iff $f(\bold{0})\in R$ is invertible.
      In this case, $$f(\bold{u})^{-1}=f(\bold{0})^{-1}\sum_{i\geq 0}\big(1-f(\bold{0})^{-1}f(\bold{u})\big)^i.$$
    \end{lemma}
    \begin{proof}
      If $f(\bold{u})$ is invertible, there is some $g(\bold{u})$ such that
      $$f(\bold{u})g(\bold{u})=1=g(\bold{u})f(\bold{u}).$$
      But then $f(\bold{0})g(\bold{0})=1=g(\bold{0})f(\bold{0})$, and $f(\bold{0})$ is invertible.

      Conversely, suppose that $f(\bold{0})$ is invertible, and let
      $$h(\bold{u})=\sum_{i\geq 0}\big(1-f(\bold{0})^{-1}f(\bold{u})\big)^i.$$
      Then
      \begin{eqnarray*}
        \big(1-f(\bold{0})^{-1}f(\bold{u})\big)h(\bold{u}) &=& \sum_{i\geq 1}\big(1-f(\bold{0})^{-1}f(\bold{u})\big)^i\\
        &=& \sum_{i\geq 0}\big(1-f(\bold{0})^{-1}f(\bold{u})\big)^i-1\\
        &=& h(\bold{u})-1\\
        1 &=& h(\bold{u})\Big(1-\big(1-f(\bold{0})^{-1}f(\bold{u})\big)\Big)\\
        &=& h(\bold{u})f(\bold{0})^{-1}f(\bold{u}),
      \end{eqnarray*}
      so $f(\bold{u})$ has inverse $h(\bold{u})f(\bold{0})^{-1}=f(\bold{0})^{-1}\sum_{i\geq 0}\big(1-f(\bold{0})^{-1}f(\bold{u})\big)^i$.
    \end{proof}

    \begin{corollary}
      \label{geo_seq}
      Let $R$ be a commutative ring.
      For all $f(\bold{u})\in R[[\bold{u}]]$ satisfying $f(\bold{0})=0$,
      $$\big(1-f(\bold{u})\big)^{-1}=\sum_{i\geq 0}f(\bold{u})^{i}.$$
    \end{corollary}
    \begin{proof}
      Take the $f$ in Lemma \ref{inv_seq} to be $1-f(\bold{u})$.
    \end{proof}

    In the proof of Lagrange Inversion Theorem \ref{lagrange_inversion_theorem}, we will need to extend $R[[z]]$ to the ring of formal Laurent series
    $$R((z))=\{\sum_{i=j}^\infty a_i z^i | j\in\mathbb{Z}, a_i\in R\}.$$
    More generally, extend $R[[\bold{u}]]$ to the ring of formal Laurent series
    $$R((\bold{u}))=\{\sum_{\bold{i}\geq\bold{j}} a_\bold{i} \bold{u}^\bold{i} | \bold{j}\in\mathbb{Z}^r, a_\bold{i}\in R\}.$$
    Then the addition and multiplication operations extend to $R((\bold{u}))$; care must still be taken for infinite sums, infinite products, and composition.
    For $k\in \mathbb{Z}$, we define the operator $\coeff{z^k}$ on $R((z))$ via $\coeff{z^k}\sum_{i=j}^\infty a_i z^i=a_k$.

    \begin{lemma}
      \label{comp_seq}
      Let $R$ be a commutative ring.
      Let $g_1(\bold{u}), g_2(\bold{u}), \dots, g_r(\bold{u})\in R[[\bold{u}]]$.
      If $g_1(\bold{0})=g_2(\bold{0})=\dots=g_r(\bold{0})=0$, then, for all $f(\bold{u})\in R[[\bold{u}]]$, the composition $(f\circ\bold{g})(\bold{u})=f\big(g_1(\bold{u}),\cdots,g_r(\bold{u})\big)$ % \left(, \right) doesn't seem to work here
  is a power series in $R[[\bold{u}]]$.
    \end{lemma}
    \begin{proof}
      Write $f(\bold{u})=\sum_{\bold{k}\in(\mathbb{Z}_{\geq 0})^r}a_\bold{k}\bold{u}^\bold{k}$ for $a_\bold{k}\in R$.
      Then
      $$(f\circ\bold{g})(\bold{u})=\sum_{\bold{k}\in(\mathbb{Z}_{\geq 0})^r}a_\bold{k}g_1(\bold{u})^{k_1}g_2(\bold{u})^{k_2}\dots g_r(\bold{u})^{k_r},$$
      where each of the summands $a_\bold{k}g_1(\bold{u})^{k_1}g_2(\bold{u})^{k_2}\dots g_r(\bold{u})^{k_r}$ is a well-defined element in $R[[\bold{u}]]$.
      Thus, it suffices to check that any monomial $\bold{u}^\bold{l}$ only appears with nonzero coefficient in finitely many of those summands.
      But each $g_i(\bold{u})$ has no constant term, so $a_\bold{k}g_1(\bold{u})^{k_1}g_2(\bold{u})^{k_2}\dots g_r(\bold{u})^{k_r}$ has no terms of (combined) degree less than $k_1+k_2+\dots +k_r$.
      In particular, computing the coefficient of the monomial $\bold{u}^\bold{l}$, one need only sum across those $\bold{k}$ for which $k_1+k_2+\dots k_r\leq l_1+l_2+\dots l_r$.
    \end{proof}

    We say $\psi(z)\in R[[z]]$ has a compositional inverse, denoted $\psi^{\compinv}(z)\in R[[z]]$, if
    $$\psi\big(\psi^{\compinv}(z)\big)=z=\psi^{\compinv}\big(\psi(z)\big).$$

    \begin{lemma}
      \label{comp_inv_exist}
      Let $R$ be a commutative ring and $\psi(z)\in R[[z]]$, where
  $$\psi(z)=z \iota(z)$$
      for some $\iota(z)\in R[[z]]$ with a multiplicative inverse.
      Then $\psi(z)$ has a compositional inverse.
    \end{lemma}
    \begin{proof}
      The following idea is straightforward to follow, but cumbersome to read, and so we skip the details.
      Write $\psi(z)=\sum_{i\geq 1}a_i z^i$, noting $a_1$ is invertible in $R$  by Lemma \ref{inv_seq}.
      Consider a power series $\sigma(z)$ with unknown coefficients, except the constant term is 0 so that it can be plugged into a composition.
      Then write out the terms of $\psi\big(\sigma(z)\big)$ and set this equal to $z$.
      This imposes a set of constraints that, one may check, will yield a unique solution on the coefficients of $\sigma(z)$; in other words, taking these values yields that $\sigma(z)$ is a right compositional inverse of $\psi(z)$.
      But $\sigma$ also satisfies the hypotheses of the result, so it has a right compositional inverse $\upsilon$.
      Applying $\psi$ to $z=\sigma(\upsilon(z))$ shows that $\upsilon=\psi$ and $\sigma=\psi^{\compinv}(z)$.
    \end{proof}

    The partial derivative of a power series is defined as the formal object one would expect.
    Namely, for $1\leq i\leq r$,
    $$\frac{\partial}{\partial u_i}\sum_{\bold{k}\in(\mathbb{Z}_{\geq j})^r}a_\bold{k}\bold{u}^\bold{k}=\sum_{\bold{k}\in(\mathbb{Z}_{\geq j})^r, k_i>0}a_\bold{k}k_i u_1^{k_1}u_2^{k_2}\dots u_{i-1}^{k_{i-1}}u_i^{k_i-1}u_{i+1}^{k_{i+1}}\dots u_r^{k_r}.$$

    Many of the usual rules for the derivative hold over formal Laurent series.

    \begin{lemma}
      \label{basic_derivatives}
      Let $R$ be a commutative ring, $f(\bold{u}), g(\bold{u})\in R((\bold{u}))$, $h(z)\in R[[z]]$, $1\leq i \leq r$, and $n\in\mathbb{Z}$.
      Then 
      \begin{eqnarray*}
        \frac{\partial}{\partial u_i}\big(f(\bold{u})g(\bold{u})\big) &=& f(\bold{u})\frac{\partial}{\partial u_i}g(\bold{u})+g(\bold{u})\frac{\partial}{\partial u_i}f(\bold{u})\\
        \frac{\partial}{\partial u_i}h\big(f(\bold{u})\big) &=& \frac{\partial}{\partial z}h(z)|_{z=f(\bold{u})}\frac{\partial}{\partial u_i}f(\bold{u}).
      \end{eqnarray*}
      If $n\in\mathbb{Z}_{\geq 0}$ or $f$ is invertible,
      then
      \begin{eqnarray*}
        \frac{\partial}{\partial u_i} f(\bold{u})^n &=& n f(\bold{u})^{n-1} \frac{\partial}{\partial u_i} f(\bold{u}).
      \end{eqnarray*}
    \end{lemma}
    \begin{proof}
      By replacing $R$ with $R((u_1,\dots,u_{i-1},u_{i+1},\dots,u_r))$ we may, without loss, assume $r=1$ and write $u=\bold{u}$.

      For the product rule, first verify it in the simple case $f(u)=u^i$ and $g(u)=u^j$.
      Then use the linearity of the derivative to extend this case to $f(u)=u^i$ and arbitrary $g(u)$, and then use linearity again to extend the result to arbitrary $f(u)$ and $g(u)$.

      For the power rule, show it via induction for $n\geq 0$.
      If $n<0$, then $f(u)^n$ is only defined if $f(u)$ is invertible (which is why that condition is in the statement of the lemma).
      Then applying the product rule and the power rule for positive exponents to $1=f(u)^n f(u)^{-n}$ yields that $0=f(u)^n(-n)f(u)^{-n-1}+f(u)^{-n}\frac{\partial}{\partial u}f(u)^{n}$, so $\frac{\partial}{\partial u}f(u)^{n}=n f(u)^{n-1}$, as required.

      For the chain rule, use the power rule and the linearity of the derivative.
    \end{proof}

    We define the logarithm to be the univariate polynomial
    $$\ln\big((1-z)^{-1}\big)=\sum_{i\geq 1}\frac{z^i}{i}.$$
    Of course, Lemma \ref{comp_seq} lets one plug in all kinds of (possibly multivariate) power series.

    Similarly, for any $\mathcal{P}\subseteq\mathbb{Z}_{\geq0}$, write
    \begin{eqnarray}
      e_\mathcal{P}(z)&=&\label{eq:eP_def}
      \sum_{n\in\mathcal{P}}\frac{z^n}{n!};
    \end{eqnarray}
    when $\mathcal{P}=\mathbb{Z}_{\geq0}$, we get the usual exponential power series
    $$e^z=e_{\mathbb{Z}_{\geq0}}(z)=\sum_{n\geq 0}\frac{z^n}{n!}.$$

    It is convenient to write
    \begin{eqnarray}
      \mathcal{P}-i&=& \label{eq:Pminus}
      \{n-i\mid n\in\mathcal{P}\}\setminus\mathbb{Z}_{<0}
    \end{eqnarray}
    for the result of subtracting $i$ from every element of $\mathcal{P}$ and keeping only those results that are nonnegative.

    \begin{lemma}
      \label{e_identities}
      Let $R\supseteq \mathbb{Q}$ be a commutative ring, $\mathcal{P}\subseteq\mathbb{Z}_{\geq 0}$, and $f(z),g(z)\in R[[z]]$ with $f(0)=0=g(0)$.
      Then
      \begin{eqnarray*}
        \frac{\partial}{\partial z}e_\mathcal{P}(z) &=& e_{\mathcal{P}-1}(z)\\
        \frac{\partial}{\partial z}e^z &=& e^z\\
        \frac{\partial}{\partial z}\ln\big((1-z)^{-1}\big) &=& (1-z)^{-1}\\
        e^{\ln \big((1-z)^{-1}\big)} &=& (1-z)^{-1}\\
        e^{f(z)+g(z)}&=& e^{f(z)}e^{g(z)}.
      \end{eqnarray*}
    \end{lemma}
    \begin{proof}
      The first three equations are immediate from the appropriate definitions and, in the case of the third, Corollary \ref{geo_seq}.

      For the fourth, note that
      \begin{eqnarray*}
        \frac{\partial}{\partial z}(1-z)e^{\ln\big((1-z)^{-1}\big)} &=& -e^{\ln\big((1-z)^{-1}\big)}+(1-z)e^{\ln\big((1-z)^{-1}\big)}\frac{\partial}{\partial z}\ln\big((1-z)^{-1}\big)\\
        &=& -e^{\ln\big((1-z)^{-1}\big)}+(1-z)e^{\ln\big((1-z)^{-1}\big)}(1-z)^{-1}\\
        &=& 0\\
        (1-z)e^{\ln\big((1-z)^{-1}\big)}|_{z=0} &=& 1.
      \end{eqnarray*}
      It is readily verified that a power series with derivative 0 has 0 coefficients everywhere except possibly for the constant term.
      Since the constant term is 1 in this case,
      $$(1-z)e^{\ln\big((1-z)^{-1}\big)}=1.$$
      By Corollary \ref{geo_seq}, the power series $1-z$ is invertible, and the fourth result holds.

      For the fifth equation, apply the Binomial Theorem in $R[[\bold{u}]]$ to see
      \begin{eqnarray*}
        e^{u_1+u_2} &=& \sum_{i=0}^\infty\frac{(u_1+u_2)^i}{i!}\\
        &=& \sum_{i=0}^\infty\frac{\sum_{j=0}^i\binom{i}{j}u_1^j u_2^{i-j}}{i!}\\
        &=& \sum_{i=0}^\infty\sum_{j=0}^i\frac{u_1^j}{j!}\frac{u_2^{i-j}}{(i-j)!};
      \end{eqnarray*}
      note that swapping the sums in the last equality is, in fact, okay because there is exactly one term contributing to the coefficient of a given $u_1^i u_2^k$ (in either expression).
      In other words, we have $e^{u_1+u_2}=e^{u_1} e^{u_2}$; since $f(0)=0=g(0)$, one may plug in $u_1=f(z)$ and $u_2=g(z)$, giving the result.
    \end{proof}

    We close this section with Lagrange Inversion Theorem \ref{lagrange_inversion_theorem}, a tool for finding an explicit formula for the compositional inverse of a power series.
    The algebraic development below requires an exploration of residues.
    In complex analysis, residues are important for the calculation of line integrals that arise in the Taylor expansion.
    In the purely algebraic realm, their utility is that $z^{-1}$ is the only Laurent monomial that cannot be integrated.
    The following observation is Proposition 1.2.1 in \cite{GouldenJackson}.

    \begin{lemma}
      \label{zero_residue}
      Let $R$ be a commutative ring and consider $f(z),g(z)\in R((z))$.
      Then
      \begin{eqnarray*}
        \coeff{z^{-1}}\frac{\partial}{\partial z}f(z)&=&0\\
        \ \coeff{z^{-1}}\left(g(z)\frac{\partial}{\partial z}f(z)\right) &=& -\coeff{z^{-1}}\left(f(z)\frac{\partial}{\partial z}g(z)\right).
      \end{eqnarray*}
    \end{lemma}
    \begin{proof}
      For the first claim, it suffices to consider only $f(z)=z^i$ for some $i\in\mathbb{Z}$; this is clearly true for $i\neq 0$ and for $i=0$.
      For the second claim, apply the first to the product $f(z)g(z)$.
    \end{proof}

    Following Section 1.1.11 in \cite{GouldenJackson}, the \emph{valuation} of $f\in R((\bold{u}))$ is
    $$\val(f)=\begin{cases}\bold{k} & \text{if } f(\bold{u})=\bold{u}^\bold{k}g(\bold{u}) \text{ for an invertible } g\in R[[\bold{u}]] \\
    \infty & \text{otherwise}
    \end{cases}.$$

    The next result is the Residue Composition Theorem 1.2.2 in \cite{GouldenJackson}.

    \begin{theorem}
      \label{residue_composition}
            {\rm (Residue Composition Theorem)} Let $R$ be a communative ring, $f(z),\psi(z)\in R((z))$, and $\val(\psi)>0$.
            Then
            $$\val(\psi)\coeff{z^{-1}}f(z)=\coeff{z^{-1}}f\big(\psi(z)\big)\frac{\partial}{\partial z}\psi(z).$$
    \end{theorem}
    \begin{proof}
      Both sides of the claim are linear in $f$, so it suffices to prove the result for $f(z)=z^n$ for each $n\in\mathbb{Z}$.

      If $n\neq -1$, Lemma \ref{zero_residue} says that the right-hand expression $\coeff{z^{-1}}\psi(z)^n\frac{\partial}{\partial z}\psi(z)=\coeff{z^{-1}}\frac{1}{n+1}\frac{\partial}{\partial z} \psi(z)^{n+1}$ is 0, and so matches the left side.
      It remains only to consider $f(z)=z^{-1}$.

      Since $\val(\psi)>0$, one may write
      $$\psi(z)=a z^k g(z)$$
      for $a\in R$ invertible, $k=\val(\psi)$, and $g(z)\in R[[z]]$ with $g(0)=1$.
      Thus,
      \begin{eqnarray*}
        \coeff{z^{-1}}\psi(z)^{-1}\frac{\partial}{\partial z}\psi(z)&=&\coeff{z^{-1}}a^{-1} z^{-k}g(z)^{-1}\big(k a z^{k-1} g(z)+a z^k \frac{\partial}{\partial z} g(z)\big)\\
        &=&\coeff{z^{-1}}\big(k z^{-1}+ g(z)^{-1}\frac{\partial}{\partial z} g(z)\big)\\
        &=&k+ \coeff{z^{-1}}\frac{\partial}{\partial z}\ln g(z),\\
      \end{eqnarray*}
      where $\ln g(z)=\ln\Big(1-\big(1-g(z)\big)\Big)$ exists by Lemma \ref{comp_seq}.
      Then Lemma \ref{zero_residue} shows that
      $$\coeff{z^{-1}}\psi(z)^{-1}\frac{\partial}{\partial z}\psi(z)=k+0,$$
      as required.
    \end{proof}

    Finally, we get a method for explicitly computing the terms of the compositional inverse of a function $\psi(z)$ that satisfies the conditions of Lemma \ref{comp_inv_exist}.
    Write $\psi(z)=z \iota(z)$ for some multiplicatively invertible $\iota(z)$; it is handy to replace $\iota(z)$ by its multiplicative inverse $\epsilon(z)=\iota(z)^{-1}$, so we consider power series of the form $\psi(z)=\frac{z}{\epsilon(z)}$ and then look for its compositional inverse.
    Said another way, the goal is to find the explicit power series of the function $\sigma(z)$ that is implicitly defined by $\sigma(z)=z\epsilon\big(\sigma(z)\big)$.

    \begin{theorem}
      \label{lagrange_inversion_theorem}
            {\rm (Lagrange Inversion Theorem)} Let $R\supseteq \mathbb{Q}$ be a commutative ring and $\epsilon(z)\in R[[z]]$ have a multiplicative inverse.
            Then there is a unique $\sigma(z)\in R[[z]]$ with $\sigma(0)=0$ such that
            $$\sigma(z)=z\epsilon\big(\sigma(z)\big).$$
            For all $f(z)\in R[[z]]$,
            $$f\big(\sigma(z)\big)=\sum_{n\geq 1}\frac{1}{n}\coeff{\lambda^{n-1}}\epsilon(\lambda)^n \frac{\partial}{\partial \lambda}f(\lambda) z^n.$$
    \end{theorem}
    \begin{proof}
      As described above, $\sigma$ exists by considering Lemma \ref{comp_inv_exist} with $\iota=\epsilon^{-1}$.
      Let
      $$\psi(z)=\frac{z}{\epsilon(z)},$$
      so
      $$\sigma=\psi^{\compinv}.$$

      For all $n>0$, the Residue Composition Theorem \ref{residue_composition} yields
      \begin{eqnarray*}
        \coeff{z^n}f\big(\sigma(z)\big) &=& \coeff{z^{-1}}z^{-(n+1)}f(\psi^{\compinv}(z))\\
        &=& \frac{1}{\val(\psi)}\coeff{z^{-1}}\psi(z)^{-(n+1)}f\Big(\psi^{\compinv}\big(\psi(z)\big)\Big)\frac{\partial}{\partial z}\psi(z).
      \end{eqnarray*}
      Since $\psi(z)$ has the form of a multiplicatively invertible power series times $z$, it has valuation 1, and
      \begin{eqnarray*}
        \coeff{z^n}f\big(\sigma(z)\big) &=& \coeff{z^{-1}}\psi(z)^{-(n+1)}f(z)\frac{\partial}{\partial z}\psi(z)\\
        &=& \frac{-1}{n}\coeff{z^{-1}}f(z)\frac{\partial}{\partial z}\psi(z)^{-n}.
      \end{eqnarray*}
      The second result in Lemma \ref{zero_residue} now gives
      \begin{eqnarray*}
        \coeff{z^n}f\big(\sigma(z)\big) &=& -\frac{-1}{n}\coeff{z^{-1}}\psi(z)^{-n}\frac{\partial}{\partial z}f(z)\\
        &=& \frac{1}{n}\coeff{z^{-1}}\frac{\epsilon(z)^{n}}{z^n}\frac{\partial}{\partial z}f(z)\\
        &=& \frac{1}{n}\coeff{z^{n-1}}\epsilon(z)^{n}\frac{\partial}{\partial z}f(z).
      \end{eqnarray*}
      Replacing $z$ by the dummy variable $\lambda$ gives $\coeff{z^n}f\big(\sigma(z)\big)= \frac{1}{n}\coeff{\lambda^{n-1}}\epsilon(\lambda)^{n}\frac{\partial}{\partial \lambda}f(\lambda)$ and, hence, the result.
    \end{proof}
  \end{subsection}

  \begin{subsection}{Definitions}
    \label{gen_fn_dfn}

    We now have the machinery necessary to associate a combinatorial class with a power series.
    In the definitions below, we continue to work over a commutative ring $R\supseteq\mathbb{Q}$.

    The \emph{exponential generating function} or \emph{generating function} of the labeled class $\mathfrak{A}$ is the formal power series
    \begin{eqnarray}
      \label{eq:gen_fun_as_seq}
      A(z) &=& \sum_{n=0}^\infty\frac{|A_n|}{n!}z^n,
    \end{eqnarray}
    where the $A_n=|\{a\in\mathfrak{A}\mid |a|=n\}|$ are the terms in the counting sequence.
    Note that
    \begin{eqnarray}
      \label{eq:gen_fun_as_count}
      A(z) &=& \sum_{a\in\mathfrak{A}}\frac{z^{|a|}}{|a|!}.
    \end{eqnarray}
    Equation \ref{eq:gen_fun_as_seq} emphasizes that this is the usual generating function of a sequence, and Equation \ref{eq:gen_fun_as_count} emphasizes that we are associating to every object a special monomial that encodes some useful information about that object.

    For example, consider the combinatorial classes $\mathfrak{A}$ and $\CYC_2(\mathfrak{A})$ in Figure \ref{cyc2_example}.
    The generating function of $\mathfrak{A}$ is
    \begin{eqnarray}
      A(z)&=& 1\frac{z^2}{2!}\nonumber\\
      &=& \label{eq:cyc2_example_A}
      \frac{z^2}{2}
    \end{eqnarray}
    and the generating function of $\CYC_2(\mathfrak{A})$ is
    \begin{eqnarray}
      \label{eq:cyc2_exampleCYC2A}
      3\frac{z^4}{4!}&=&\frac{z^4}{8}.
    \end{eqnarray}

    Similarly, the generating function of the class $\mathfrak{A}$ shown in Figure \ref{set2_example} is
    \begin{eqnarray}
      A(z) &=& 1\frac{z}{1!}+1\frac{z^2}{2!}\nonumber\\
      &=& \label{eq:set2_example_A}
      z+\frac{z^2}{2},
    \end{eqnarray}
    while the generating function of $\SET_2(\mathfrak{A})$ in the same picture is
    \begin{eqnarray}
      \label{eq:set2_example_SET2A}
      1\frac{z^2}{2!}+3\frac{z^3}{3!}+3\frac{z^4}{4!}&=&\frac{z^2}{2}+\frac{z^3}{2}+\frac{z^4}{8}.
    \end{eqnarray}

    If we have a parameter function $\chi:\mathfrak{A}\rightarrow (\mathbb{Z}_{\geq0})^r$, we can generalize the univariate generating function $A(z)$ to incorporate the additional structure encoded in $\chi$.
    Recall that, given $\bold{k}=(k_1,k_2,\dots,k_r)\in(\mathbb{Z}_{n\geq0})^r$, we write $\bold{u}^\bold{k}$ for $\prod_{i=1}^r u_i^{k_i}$.

    The \emph{multivariate generating function} or \emph{generating function} of the labeled class $\mathfrak{A}$ and parameter $\chi:\mathfrak{A}\rightarrow (\mathbb{Z}_{n\geq0})^r$ is
    $$A(\bold{u}, z) = \sum_{(n,\bold{k})\in\mathbb{Z}_{n\geq0}\times(\mathbb{Z}_{n\geq0})^r}\frac{|A_{n,\bold{k}}|}{n!}\bold{u}^\bold{k}z^n,$$
    where $A_{n,\bold{k}}=|\{a\in\mathfrak{A}\mid |a|=n, \chi(a)=\bold{k}\}|$.
    Like before, an equivalent alternative is to express it as
    $$A(\bold{u},z) = \sum_{a\in\mathfrak{A}}\frac{\bold{u}^{\chi(a)}z^{|a|}}{|a|!}.$$

    By way of a feeble defense of the choice to call both the original and the multivariate generating functions $A$, note that evaluating some $u_i$ at 1 has the effect of ignoring that particular parameter.
    In particular, $A(z)=A(\bold{u},z)|_{\bold{u}=(1,1,\dots,1)}$.
    The applications below all consider the situation where $r=1$, in which case we will drop the vector notation for $\bold{u}$.
    That is,
    $$A(u,z)=\sum_{a\in\mathfrak{A}}\frac{{u}^{\chi(a)}z^{|a|}}{|a|!}.$$

    For an example of how multivariate generating functions encode additional structure, consider the combinatorial class of labeled rooted trees of size at most three (these are all shown explicitly in Figure \ref{labeled_rooted_trees}).
    If the parameter function $\chi$ sends a tree to the number of leaves in that tree, then the corresponding multivariate generating function is $1\frac{u z}{1!}+2\frac{u z^2}{2!}+6\frac{u z^3}{3!}+3\frac{u^2 z^3}{3!}$.
  \end{subsection}

  \begin{subsection}{Constructions}
    \label{construction_gen_fn}

    Of course, it is not feasible to compute generating functions from an explicit list of the elements of a combinatorial class, as was done in the examples of the previous section.
    So, our next task is to address how the constructions of Section \ref{constructions} translate to generating functions.

    \begin{lemma}
      Let $\mathfrak{B}$ and $\mathfrak{C}$ be disjoint labeled classes with parameter functions of the same dimension.
      Write $A(\bold{u},z), B(\bold{u},z)$, and $C(\bold{u},z)$ for respective generating functions of $\mathfrak{B}+\mathfrak{C}, \mathfrak{B}$, and $\mathfrak{C}$ over a commutative ring $R\supseteq\mathbb{Q}$.
      Then
      $$A(\bold{u},z)=B(\bold{u},z)+C(\bold{u},z).$$
    \end{lemma}
    \begin{proof}
      Write $\psi$ and $\omega$ for the respective parameter functions of $\mathfrak{B}$ and $\mathfrak{C}$, and write $\chi$ for the parameter function they induce on $\mathfrak{B}+\mathfrak{C}$.
      Then
      \begin{eqnarray*}
        A(\bold{u},z) &=& \sum_{\alpha\in\mathfrak{B}\cup\mathfrak{C}}\frac{{\bold{u}}^{\chi(a)}z^{|a|}}{|a|!}\\
        &=& \sum_{\alpha\in\mathfrak{B}}\frac{{\bold{u}}^{\chi(\alpha)}z^{|\alpha|}}{|\alpha|!}+\sum_{\alpha\in\mathfrak{C}}\frac{{\bold{u}}^{\chi(\alpha)}z^{|\alpha|}}{|\alpha|!}\\
        &=& \sum_{\alpha\in\mathfrak{B}}\frac{{\bold{u}}^{\psi(\alpha)}z^{|\alpha|}}{|\alpha|!}+\sum_{\alpha\in\mathfrak{C}}\frac{{\bold{u}}^{\omega(\alpha)}z^{|\alpha|}}{|\alpha|!}\\
        &=& B(\bold{u},z)+C(\bold{u},z).
      \end{eqnarray*}
    \end{proof}

    The crux of the section is probably the following result.
    The proof also reveals why terms in generating functions must be scaled by $n!$; it is so there are no leftover terms when computing or unwrapping the convolutional product.

    \begin{lemma}
      \label{seq_gen_fn}
      Let $\mathfrak{B}$ and $\mathfrak{C}$ be labeled classes with parameter functions of the same dimension.
      Write $A(\bold{u},z)$, $B(\bold{u},z)$, and $C(\bold{u},z)$ for respective generating functions of $\mathfrak{B}\star\mathfrak{C}$, $\mathfrak{B}$, and $\mathfrak{C}$ over a commutative ring $R\supseteq\mathbb{Q}$.
      Then
      $$A(\bold{u},z)=B(\bold{u},z)C(\bold{u},z).$$
    \end{lemma}
    \begin{proof}
      Write $\psi$ and $\omega$ for the respective parameter functions of $\mathfrak{B}$ and $\mathfrak{C}$; recall the induced parameter function $\chi$ on $\mathfrak{B}\star\mathfrak{C}$ sends every element of $\beta\star\gamma$ to $\psi(\beta)+\omega(\gamma)$.
      Then
      \begin{eqnarray*}
        A(\bold{u},z) &=& \sum_{\alpha\in\mathfrak{B}\star\mathfrak{C}}\frac{{\bold{u}}^{\chi(\alpha)}z^{|\alpha|}}{|\alpha|!}\\
        &=& \sum_{\beta\in\mathfrak{B}}\sum_{\gamma\in\mathfrak{C}}\sum_{\alpha\in\beta\star\gamma}\frac{{\bold{u}}^{\chi(\alpha)}z^{|\alpha|}}{|\alpha|!}\\
        &=& \sum_{\beta\in\mathfrak{B}}\sum_{\gamma\in\mathfrak{C}}\sum_{\alpha\in\beta\star\gamma}\frac{{\bold{u}}^{\psi(\beta)+\omega(\gamma)}z^{|\beta|+|\gamma|}}{(|\beta|+|\gamma|)!}\\
        &=& \sum_{\beta\in\mathfrak{B}}\sum_{\gamma\in\mathfrak{C}}\binom{|\beta|+|\gamma|}{|\beta|}\frac{{\bold{u}}^{\psi(\beta)+\omega(\gamma)}z^{|\beta|+|\gamma|}}{(|\beta|+|\gamma|)!}\\
        &=& \sum_{\beta\in\mathfrak{B}}\sum_{\gamma\in\mathfrak{C}}\frac{{\bold{u}}^{\psi(\gamma)}{\bold{u}}^{\omega(\gamma)}z^{|\beta|}z^{|\gamma|}}{|\beta|!|\gamma|!}\\
        &=&B(\bold{u},z)C(\bold{u},z),
      \end{eqnarray*}
      as claimed.
    \end{proof}

    The remaining constructions now quickly fall into place.

    \begin{lemma}
      \label{gen_fn_constructions}
      Let $\mathfrak{B}$ be a labeled class with multivariate generating function $B(\bold{u},z)$ over a commutative ring $R\supseteq\mathbb{Q}$.
      Then the generating function of $\SEQ_k(\mathfrak{B})$ is $B(\bold{u},z)^k$.
      If $\mathfrak{B}$ has no objects of size 0 and if $k>0$, then the generating function of $\CYC_k(\mathfrak{B})$ is $\frac{B(\bold{u},z)^k}{k}$.
      If $\mathfrak{B}$ has no objects of size 0, then the generating function of $\SEQ(\mathfrak{B})$ is $\frac{1}{1-B(\bold{u},z)}$ and the generating function of $\CYC(\mathfrak{B})$ is $\ln\Big(\big(1-B(\bold{u},z)\big)^{-1}\Big)$.
      If $\mathfrak{B}$ has no objects of size 0 and if $\mathcal{P}\subseteq\mathbb{Z}_{\geq0}$, then the generating function of $\SET_\mathcal{P}(\mathfrak{B})$ is $e_{\mathcal{P}}\big(B(\bold{u},z)\big)$; in particular, the generating function of $\SET(\mathfrak{B})$ is $e^{B(\bold{u},z)}$.
    \end{lemma}
    \begin{proof}
      Recall that $\SEQ_k(\mathfrak{B})$ is $\mathfrak{B}\star\mathfrak{B}\star\dots\star\mathfrak{B}$, so, by Lemma \ref{seq_gen_fn}, the generating function is $B(\bold{u},z)B(\bold{u},z)\dots B(\bold{u},z)=B(\bold{u},z)^k$.

      Hereafter, assume $\mathfrak{B}$ has no objects of size 0.
      Then its generating function evaluated at $z=0$ is $B(\bold{u},0)=0$; in particular, $B(\bold{0},0)=0$, and Lemma \ref{comp_seq} ensures it can be plugged into any other power series.
      In particular, this ensures the existence of $\ln\Big(\big(1-B(\bold{u},z)\big)^{-1}\Big)$.

      The class $\CYC_k(\mathfrak{B})$ consists of equivalence classes of elements in $\SEQ_k(\mathfrak{B})$ where cyclic shifts of a sequence are identified.
      Note that there are always $k$ elements in an equivalence class, since a shift of exactly $k$ always gives the original sequence and since the fact that sequences are labeled and have no entries of size 0 ensures there is no way for a shift of less than $k$ to yield the original sequence.
      In other words, the generating function for $\CYC_k(\mathfrak{B})$ is the generating function of $\SEQ_k(\mathfrak{B})$ with every term scaled down by a factor of $k$.

      $\SEQ(\mathfrak{B})$, when it is defined, is the union of the $\SEQ_k(\mathfrak{B})$ for $k\geq 0$; these unions are disjoint, so the overall generating function is just the sum.
      In other words, the generating function of $\SEQ(\mathfrak{B})$ is given by $1+B(\bold{u},z)+B(\bold{u},z)^2+B(\bold{u},z)^3+\dots$, which can be rewritten as $\big(1-B(\bold{u},z)\big)^{-1}$ by Corollary \ref{geo_seq}.

      For $\CYC(\mathfrak{B})$, note that $\CYC(\mathfrak{B})$ is the union of the $\CYC_k(\mathfrak{B})$ for $k\geq 1$; these unions are disjoint, so the overall generating function is just the sum.
      In other words, the generating function of $\SEQ(\mathfrak{B})$ is
      $$B(\bold{u},z)+\frac{B(\bold{u},z)^2}{2}+\frac{B(\bold{u},z)^3}{3}+\dots=\ln\Big(\big(1-B(\bold{u},z)\big)^{-1}\Big).$$

      The class $\SET_{k}(\mathfrak{B})$ consists of equivalence classes of elements in $\SEQ_k(\mathfrak{B})$ where permutations on the coordinates of the sequence are identified.
      Note that there are always $k!$ elements in an equivalence class, since there are $k!$ permutations on the coordinates and since the fact that sequences are labeled with no entries of size 0 ensures there is no way for the same entry to appear twice in a single sequence.
      In other words, the generating function for $\SET_{k}(\mathfrak{B})$ is the generating function of $\SEQ_k(\mathfrak{B})$ with every term scaled down by a factor of $k!$.
      But $\SET_\mathcal{P}(\mathfrak{B})$ is the union of the $\SET_{k}(\mathfrak{B})$ for $k\in\mathcal{P}$; these unions are disjoint, so the overall generating function is just the sum $\sum_{k\in\mathcal{P}}\frac{B(\bold{u},z)^k}{k!}=e_\mathcal{P}\big(B(\bold{u},z)\big)$.
    \end{proof}

    Recall that Equation \ref{eq:cyc2_example_A} found the generating function of the $\mathfrak{A}$ in Figure \ref{cyc2_example} to be $A(z)=\frac{z^2}{2}$.
    By Lemma \ref{gen_fn_constructions}, the generating function of $\CYC_2(\mathfrak{A})$ is $\frac{A(z)^2}{2}=\frac{z^4}{8}$, which matches the computation, done via in inspection, in Equation \ref{eq:cyc2_exampleCYC2A}.
    The generating function of $\CYC_3(\mathfrak{A})$ is $\frac{A(z)^3}{3}=30\frac{z^6}{6!}$; one of the 30 elements in this class was shown in Figure \ref{cyc3_example}.

    Similarly, Equation \ref{eq:set2_example_A} gives the generating function of the $\mathfrak{A}$ shown in Figure \ref{set2_example} as $A(z)=z+\frac{z^2}{2}$.
    Thus, Lemma \ref{gen_fn_constructions} gives the generating function of $\SET_2(\mathfrak{A})$ as $\frac{A(z)^2}{2!}=\frac{z^2}{2}+\frac{z^3}{2}+\frac{z^4}{8}$, as shown in Equation \ref{eq:set2_example_SET2A}.
  \end{subsection}

  \begin{subsection}{Treecursion and Functions as Graphs}
    \label{treecursion}

    The machinery for the main combinatorial results are now in place.
    Before tackling them, we work through a simpler derivation of the univariate generating function for the class of rooted labeled trees and use this to get the generating function for functions from some $\set{n}$ to itself.

    This goal is a little contrived (there are plainly $n^n$ functions from $\set{n}$ to itself, so the generating function is $\sum_{n\geq 0}n^n\frac{z^n}{n!}$), but the approach taken below foreshadows the machinery for counting images in the next section.

    Write $\mathfrak{F}$ for the combinatorial class of functions from $\set{n}$ to itself for all $n\geq 0$.
    By identifying such a function $f$ with the digraph consisting of edges $i\rightarrow j$ iff $f(i)=j$, we see that we are counting all labeled graphs where each vertex has outdegree 1.
    See Figure \ref{fn_graph}.
    Figure \ref{fn_as_setcycletree} shows the same function decomposed into a set of cycles of rooted trees; this decomposition plainly works for any function and any such (well-labeled) decomposition is such a function.

    \begin{figure}
      \begin{center}
        \begin{tikzpicture}
          \foreach \h/\w/\n in {.5/1/.6}{ % distance between rows, distance between columns, node size

            % vertices
            \foreach \i/\j/\l in {0/0/19, 1/0/10, 2/1/3, 3/0/16, 2/-1/5, 2/2/4, 1/2/15, 2/3/12, 3/-2/7, 2/-3/6, 3/-3/20, 4/-3/18, 3/-4/17, 4/1/14, 5/1/9, 5/0/8, 6/0/2, 7/2/11, 7/1/1, 7/0/13}{
              \draw (\i*\w, \j*\h) node[scale=\n] (node\l) {\l};
            }

            % edges
            \foreach \l/\m in {19/10, 10/3, 3/16, 16/5, 5/10, 4/3, 15/4, 12/4, 7/16, 6/7, 20/7, 18/16, 17/20, 14/8, 9/8, 11/1, 1/13%, 2/8, 8/2, 13/13
            }{
              \draw[-stealth] (node\l) -- (node\m);
            }
            \draw[-stealth] (node2.north west) .. controls +(0,.25) and +(0,.25) .. (node8.north east);
            \draw[-stealth] (node8.south east) .. controls +(0,-.25) and +(0,-.25) .. (node2.south west);
            \draw[-stealth] (node13.south east) .. controls +(0,-.25) and +(0,-.25) .. (node13.south west);
          }
        \end{tikzpicture}
      \end{center}
      \caption{A Function as a Graph}
      \label{fn_graph}
    \end{figure}
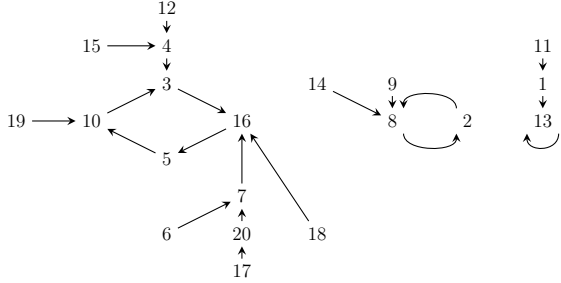

    \begin{figure}
      \begin{center}
        \begin{tikzpicture}
          \foreach \h/\w/\n in {.5/1/.6}{ % distance between rows, distance between columns, node size

            % vertices
            \foreach \i/\j/\l in {0/0/19, 1/0/10, 2/1/3, 3/0/16, 2/-1/5, 2/2/4, 1/2/15, 2/3/12, 3/-2/7, 2/-3/6, 3/-3/20, 4/-3/18, 3/-4/17, 4/1/14, 5/1/9, 5/0/8, 6/0/2, 7/2/11, 7/1/1, 7/0/13}{
              \draw (\i*\w, \j*\h) node[scale=\n] (node\l) {\l};
            }

            % edges
            \foreach \l/\m in {19/10, 10/3, 3/16, 16/5, 5/10, 4/3, 15/4, 12/4, 7/16, 6/7, 20/7, 18/16, 17/20, 14/8, 9/8, 11/1, 1/13%, 2/8, 8/2, 13/13
            }{
              \draw[-stealth] (node\l) -- (node\m);
            }
            \draw[-stealth] (node2.north west) .. controls +(0,.25) and +(0,.25) .. (node8.north east);
            \draw[-stealth] (node8.south east) .. controls +(0,-.25) and +(0,-.25) .. (node2.south west);
            \draw[-stealth] (node13.south east) .. controls +(0,-.25) and +(0,-.25) .. (node13.south west);

            % label the trees
            \draw(-1*\w,3*\h) node {{\color{green} tree}};
            % root 10
            \draw[green] (-.5*\w,0*\h) .. controls +(0,.25) and +(0,.25) .. (1.5*\w,0*\h); % 19 -> 10
            \draw[green] (-.5*\w,0*\h) .. controls +(0,-.25) and +(0,-.25) .. (1.5*\w,0*\h); % 19 -> 10
            % root 3
            \draw[green] (.5*\w,2*\h) .. controls +(0,.25) and +(0,.25) .. (2*\w,3.5*\h); % 15 -> 12
            \draw[green] (.5*\w,2*\h) .. controls +(0,-.25) and +(0,-.25) .. (2*\w,.5*\h); % 15 -> 3
            \draw[green] (2*\w,3.5*\h) .. controls +(.25,0) and +(.25,0) .. (2*\w,.5*\h); % 12 -> 3
            % root 16
            \draw[green] (3*\w,.5*\h) .. controls +(0,.25) and +(0,.25) .. (1.5*\w,-3*\h); % 16 -> 6
            \draw[green] (3*\w,.5*\h) .. controls +(0,.25) and +(0,.25) .. (4.5*\w,-3*\h); % 16 -> 18
            \draw[green] (3*\w,-4.5*\h) .. controls +(0,-.25) and +(0,-.25) .. (1.5*\w,-3*\h); % 17 -> 6 
            \draw[green] (3*\w,-4.5*\h) .. controls +(0,-.25) and +(0,-.25) .. (4.5*\w,-3*\h); % 17 -> 18
            % root 5
            \draw[green] (1.5*\w,-1*\h) .. controls +(0,.25) and +(0,.25) .. (2.5*\w,-1*\h); % 5 -> 5
            \draw[green] (1.5*\w,-1*\h) .. controls +(0,-.25) and +(0,-.25) .. (2.5*\w,-1*\h); % 5 -> 5
            % root 8
            \draw[green] (3.5*\w,1*\h) .. controls +(0,.25) and +(0,.25) .. (5*\w,1.5*\h); % 14 -> 9
            \draw[green] (3.5*\w,1*\h) .. controls +(0,-.25) and +(0,-.25) .. (5*\w,-.5*\h); % 14 -> 8
            \draw[green] (5*\w,1.5*\h) .. controls +(.25,0) and +(.25,0) .. (5*\w,-.5*\h); % 9 -> 8
            % root 2
            \draw[green] (5.5*\w,0*\h) .. controls +(0,.25) and +(0,.25) .. (6.5*\w,0*\h); % 2 -> 2
            \draw[green] (5.5*\w,0*\h) .. controls +(0,-.25) and +(0,-.25) .. (6.5*\w,0*\h); % 2 -> 2
            % root 13
            \draw[green] (7*\w,2.5*\h) .. controls +(.25,0) and +(.25,0) .. (7*\w,-.5*\h); % 11 -> 13
            \draw[green] (7*\w,2.5*\h) .. controls +(-.25,0) and +(-.25,0) .. (7*\w,-.5*\h); % 11 -> 13

            % label the cycles
            \draw(-1*\w,2*\h) node {{\color{blue} cycle}};
            % cycle 3
            \draw[blue] (-.5*\w,0*\h) .. controls +(0,.25) and +(0,.25) .. (2*\w,3.5*\h); % 19 -> 12
            \draw[blue] (-.5*\w,0*\h) .. controls +(-.25,0) and +(-.25,0) .. (3*\w,-5*\h); % 19 -> 17
            \draw[blue] (2*\w,3.5*\h) .. controls +(0,.25) and +(0,.25) .. (4.5*\w,-3*\h); % 12 -> 18
            \draw[blue] (3*\w,-5*\h) .. controls +(0,-.25) and +(0,-.25) .. (4.5*\w,-3*\h); % 17 -> 18
            % cycle 2
            \draw[blue] (3.5*\w,1*\h) .. controls +(0,.25) and +(0,.25) .. (5*\w,2*\h); % 14 -> 9
            \draw[blue] (3.5*\w,1*\h) .. controls +(0,-.25) and +(0,-.25) .. (5*\w,-1*\h); % 14 -> 8
            \draw[blue] (5*\w,-1*\h) .. controls +(0,-.25) and +(0,-.25) .. (6.5*\w,0*\h); % 8 -> 2
            \draw[blue] (5*\w,2*\h) .. controls +(0,.25) and +(0,.25) .. (6.5*\w,0*\h); % 9 -> 2
            % cycle 13
            \draw[blue] (7*\w,2.5*\h) .. controls +(.25,.25) and +(.25,.25) .. (7*\w,-.5*\h); % 11 -> 13
            \draw[blue] (7*\w,2.5*\h) .. controls +(-.25,-.25) and +(-.25,-.25) .. (7*\w,-.5*\h); % 11 -> 13

            % label the set
            \draw(-1*\w,1*\h) node {{\color{red} set}};
            % set 2
            \draw[red] (-.5*\w,0*\h) .. controls +(.25,.25) and +(.25,.25) .. (2*\w,3.5*\h); % 19 -> 12
            \draw[red] (-.5*\w,0*\h) .. controls +(-.25,-.25) and +(-.25,-.25) .. (3*\w,-5*\h); % 19 -> 17
            \draw[red] (7*\w,2.5*\h) .. controls +(.25,.25) and +(.25,.25) .. (2*\w,3.5*\h); % 11 -> 12
            \draw[red] (7*\w,2.5*\h) .. controls +(.25,.25) and +(.25,.25) .. (7*\w,-.5*\h); % 11 -> 13
            \draw[red] (3*\w,-5*\h) .. controls +(-.25,-.25) and +(-.25,-.25) .. (7*\w,-.5*\h); % 17 -> 13
          }
        \end{tikzpicture}
      \end{center}
      \caption{A Function as a Set of Cycles of Rooted Trees}
      \label{fn_as_setcycletree}
    \end{figure}

    In other words, write $\mathfrak{T}$ for the combinatorial class of rooted labeled trees and $T(z)$ for its generating function.
    We then have $\mathfrak{F}=\SET\big(\CYC(\mathfrak{T})\big)$, noting that this construction is well-defined since every object in $\mathfrak{T}$ has a root, and so is nonempty, and since the $\CYC$ construction does not produce any objects of size 0.
    Thus, Lemma \ref{gen_fn_constructions} says that
    $$F(z)=e^{\ln\Big(\big(1-T(z)\big)^{-1}\Big)}.$$
    By Lemma \ref{e_identities}, this simplifies to
    \begin{eqnarray}
      F(z)&=& \label{eq:unconstrained_function}
      \big(1-T(z)\big)^{-1}.
    \end{eqnarray}

    As a quick aside, note that by saying a function is a set of cycles of trees (as opposed to a non-empty set of cycles of trees), we are implicitly adopting the convention that there is a unique function from the empty set to itself.
    This has the advantage of fitting with the convention that $0^0=1$, fitting with the usual function-as-ordered-pairs definition, and, most importantly, making the resulting generating functions look cleaner.

    We have reduced understanding $F(z)$ to understanding $T(z)$.
    Write $\mathfrak{V}$ (for vertex) for the labeled class whose only object has size 1.
    The generating function for $\mathfrak{V}$ is
    $$V(z)=z.$$
    Note every rooted labeled tree can be written as an element in the labeled product of a vertex (the root) and a set of rooted labeled trees, and conversely; for example, see Figure \ref{tree_decomp} for a decomposition of one of the (weakly-labeled) trees in Figure \ref{fn_as_setcycletree}.
    Thus, one observes that $\mathfrak{T}=\mathfrak{V}\star \SET(\mathfrak{T})$ and
    \begin{eqnarray}
      T(z)&=& \label{eq:unconstrained_tree}
      z e^{T(z)}.
    \end{eqnarray}
    This is a special case of the second result in Theorem \ref{when_u_equals_1_tree}.

    \begin{figure}
      \begin{center}
        \begin{tikzpicture}
          \foreach \h/\w/\n in {.5/1/.6}{ % distance between rows, distance between columns, node size
            % vertices
            \foreach \i/\j/\l in {3/0/16, 3/-2/7, 2/-3/6, 3/-3/20, 4/-3/18, 3/-4/17}{
              \draw (\i*\w, \j*\h) node[scale=\n] (node\l) {\l};
            }

            % edges
            \foreach \l/\m in {7/16, 6/7, 20/7, 18/16, 17/20}{
              \draw[-stealth] (node\l) -- (node\m);
            }

            % label the trees
            \draw(1*\w,0*\h) node {{\color{green} tree}};
            % root 7
            \draw[green] (3*\w,-4.5*\h) .. controls +(0,-.25) and +(0,-.25) .. (1.5*\w,-3*\h); % 17 -> 6
            \draw[green] (3*\w,-4.5*\h) .. controls +(.25,0) and +(.25,0) .. (3*\w,-1.5*\h); % 17 -> 7
            \draw[green] (1.5*\w,-3*\h) .. controls +(0,.25) and +(0,.25) .. (3*\w,-1.5*\h); % 6 -> 7
            % root 18
            \draw[green] (3.5*\w,-3*\h) .. controls +(0,.25) and +(0,.25) .. (4.5*\w,-3*\h); % 18 -> 18
            \draw[green] (3.5*\w,-3*\h) .. controls +(0,-.25) and +(0,-.25) .. (4.5*\w,-3*\h); % 18 -> 18

            % label the set
            \draw(1*\w,-1*\h) node {{\color{blue} set}};
            % set 7
            \draw[blue] (3*\w,-4.5*\h) .. controls +(-.25,-.25) and +(-.25,-.25) .. (1.5*\w,-3*\h); % 17 -> 6
            \draw[blue] (1.5*\w,-3*\h) .. controls +(.25,.25) and +(.25,.25) .. (3*\w,-1.5*\h); % 6 -> 7
            \draw[blue] (3*\w,-1.5*\h) .. controls +(0,.25) and +(0,.25) .. (4.5*\w,-3*\h); % 7 -> 18
            \draw[blue] (3*\w,-4.5*\h) .. controls +(0,-.25) and +(0,-.25) .. (4.5*\w,-3*\h); % 17 -> 18

            % label *the* root
            \draw(1*\w,-2*\h) node {{\color{red} *the* root}};
            % root 16
            \draw[red] (2.5*\w,-0*\h) .. controls +(.25,.25) and +(.25,.25) .. (3.5*\w,0*\h); % 16 -> 16
            \draw[red] (2.5*\w,-0*\h) .. controls +(-.25,-.25) and +(-.25,-.25) .. (3.5*\w,0*\h); % 16 -> 16

            % label the big tree
            % root 16
            \draw[green, thick] (3*\w,-4.5*\h) .. controls +(-.25,0) and +(-.25,0) .. (1.5*\w,-3*\h); % 17 -> 6
            \draw[green, thick] (1.5*\w,-3*\h) .. controls +(-.25,0) and +(-.25,0) .. (3*\w,.5*\h); % 6 -> 16
            \draw[green, thick] (3*\w,.5*\h) .. controls +(.25,0) and +(.25,0) .. (4.5*\w,-3*\h); % 16 -> 18
            \draw[green, thick] (3*\w,-4.5*\h) .. controls +(-.25,-.25) and +(-.25,-.25) .. (4.5*\w,-3*\h); % 17 -> 18
          }
        \end{tikzpicture}
      \end{center}
      \caption{A Tree is a Root and a Set of Subtrees}
      \label{tree_decomp}
    \end{figure}

    Applying Lagrange Inversion Theorem \ref{lagrange_inversion_theorem} with $\epsilon=e$ (so that $\sigma=T$) and $f(z)=z$ yields that
    \begin{eqnarray}
      \coeff{z^n}T(z) &=& \label{eq:tree_coefficient}
      \frac{1}{n}\coeff{\lambda^{n-1}}(e^{\lambda})^n.
    \end{eqnarray}
    The fifth result in Lemma \ref{e_identities} implies $(e^{\lambda})^n=e^{n\lambda}$, so
    \begin{eqnarray*}
      \coeff{z^n}T(z) &=&  \frac{1}{n}\coeff{\lambda^{n-1}}e^{n\lambda}\\ % by t
      &=& \frac{1}{n}\frac{n^{n-1}}{(n-1)!}\\
      &=& \frac{n^{n-1}}{n!}.
    \end{eqnarray*}
    That is,
    $$T(z) = \sum_{n\geq 1} n^{n-1}\frac{z^n}{n!},$$
    and so there are $n^{n-1}$ rooted labeled trees with $n$ vertices.
    This observation is equivalent to Cayley's formula that there are $n^{n-2}$ labeled trees of size $n$, a standard example to demonstrate the utility of Lagrange Inversion; see Section 3.3.10 of \cite{GouldenJackson} or Section 2.1 of \cite{FlajoletOdlyzko_nolink} or Equation 47 in Section II.5.1 of \cite{FlajoletSedgewick}.
    See Corollary \ref{tree_coefficient} below for a generalization of Equation \ref{eq:tree_coefficient}.

    Differentiating Equation \ref{eq:unconstrained_tree} and solving for $\frac{\partial}{\partial z} T(z)$ yields
    \begin{eqnarray*}
      \frac{\partial}{\partial z} T(z) &=&  e^{T(z)}\left(1-z e^{T(z)}\right)^{-1}.
    \end{eqnarray*}
    By Equations \ref{eq:unconstrained_function} and \ref{eq:unconstrained_tree}, this can be written as
    \begin{eqnarray*}
      z \frac{\partial}{\partial z} T(z) &=& T(z)F(z),
    \end{eqnarray*}
    so
    \begin{eqnarray*}
      \coeff{z^n}F(z) &=& \coeff{z^n} z \frac{\frac{\partial}{\partial z} T(z)}{T(z)}\\
      &=& \coeff{z^{-1}} z^{-n} \frac{\frac{\partial}{\partial z} T(z)}{T(z)}.
    \end{eqnarray*}
    Since Equation \ref{eq:unconstrained_tree} says that $T(z)$ is a compositional inverse of $\psi(z)=\frac{z}{e^z}$, this can be written
    \begin{eqnarray*}
      \coeff{z^n}F(z) &=& \coeff{z^{-1}} \left(\frac{T(z)}{e^{T(z)}}\right)^{-n} \frac{\frac{\partial}{\partial z} T(z)}{T(z)}.
    \end{eqnarray*}
    Taking $\psi(z)=T(z)$ and $f(z)=\frac{e^{n z}}{z^{n+1}}$ in the Residue Composition Theorem \ref{residue_composition} simplifies this to
    \begin{eqnarray*}
      \coeff{z^n}F(z) &=& \val(T)\coeff{z^{-1}} \frac{e^{n z}}{z^{n+1}}\\
      &=& \coeff{z^n} e^{n z}\\
      &=& \frac{n^n}{n!}.
    \end{eqnarray*}
    That is,
    $$F(z) = \sum_{n\geq 0} n^{n}\frac{z^n}{n!},$$
    and, as expected, there are $n^{n}$ functions from $\set{n}$ to itself.
    See Corollary \ref{function_coefficients} for a generalization of these calculations.
\end{subsection}
\end{section}

\begin{section}{Combinatorics with Preimage Constraints}
  \label{preimage_combinatorics}

  This section is devoted to deriving generating functions that encode parameters related to functions from a finite set to itself subject to the preimage constraint $\mathcal{P}$, by which we mean $\mathcal{P}\subseteq\mathbb{Z}_{\geq 0}$ and we are restricting our attention to functions $f:\set{n}\rightarrow\set{n}$ such that for any $1\leq x\leq n$,
  $$|f^{-1}(x)|\in\mathcal{P}.$$

  For $\mathcal{P}\subseteq\mathbb{Z}_{\geq0}$, recall from Equations \ref{eq:eP_def} and \ref{eq:Pminus} that $$e_\mathcal{P}(z)=\sum_{n\in\mathcal{P}}\frac{z^n}{n!}$$ and, if $i\in\mathbb{Z}_{\geq0}$,
  $$\mathcal{P}-i=\{n-i\mid n\in\mathcal{P}\}\setminus\mathbb{Z}_{<0}.$$

  For $\mathcal{P}\subseteq\mathbb{Z}_{\geq0}$ with $0\in\mathcal{P}$, define
  \begin{eqnarray}
    \label{eq:psi}
    \psi_\mathcal{P}(z)&=&\frac{z}{e_\mathcal{P}(z)}.
  \end{eqnarray}
  Note that, by Lemma \ref{comp_inv_exist}, the generating function $\psi_\mathcal{P}^{\compinv}$ exists.

  \begin{subsection}{Trees}
    \label{trees}

    The basic building blocks for all of the constructions that follow are rooted trees subject to a preimage constraint $\mathcal{P}$, by which we mean a rooted tree such that for every node, the number of neighbors that are further from the root than the node itself is in $\mathcal{P}$.

    The additional label $u$ in the next definition is not needed until Section \ref{image_combinatorics}, but we include it here, a little early, in order to avoid repeating all of the definitions and proofs verbatim.
    (See Theorem \ref{when_u_equals_1_tree} to formally justify that one can get an analysis analogous to that of Section \ref{treecursion} by simply plugging in $u=1$.)

    For $h,k\geq 0$ and $\mathcal{P}\subseteq\mathbb{Z}_{\geq0}$, write $\mathfrak{T}_h^{k,\mathcal{P}}$ for the combinatorial class of rooted trees of height exactly $h$ where the number of preimages of any given node is an element of $\mathcal{P}$, where nodes are marked by $z$ and nodes whose maximal distance from a leaf is less than $k$ are marked by $u$.
    Write $\mathfrak{T}_{\leq h}^{k,\mathcal{P}}=\cup_{0\leq i\leq h}\mathfrak{T}_i^{k,\mathcal{P}}$ for the class of these trees of height at most $h$ and $\mathfrak{T}^{k,\mathcal{P}}=\cup_{0\leq i}\mathfrak{T}_i^{k,\mathcal{P}}$ for the class of these trees of any height.
    Write ${T}_h^{k,\mathcal{P}}(u,z)$, ${T}_{\leq h}^{k,\mathcal{P}}(u,z)=\sum_{i=0}^h T_i^{k,\mathcal{P}}(u,z)$, ${T}^{k,\mathcal{P}}(u,z)=\sum_{i=0}^\infty T_i^{k,\mathcal{P}}(u,z)$, respectively, for the generating functions of $\mathfrak{T}_{h}^{k,\mathcal{P}}$, $\mathfrak{T}_{\leq h}^{k,\mathcal{P}}$, and $\mathfrak{T}^{k,\mathcal{P}}$.

    For example, the class $\mathfrak{T}_{1}^{1,\{0,1,2\}}$ consists of all labeled rooted trees where the roots have 1 or 2 neighbors which are leaves marked with $u$; the height constraint $h=1$ precludes the root from having 0 preimages and the nonroots from having any preimages.
    See Figure \ref{T_ht1_k1_max2preimages}, where vertices marked with $u$ are circled.
    Thus,
    \begin{eqnarray}
      {T}_{1}^{1,\{0,1,2\}}(u,z) &=& 2\frac{u z^2}{2!}+3\frac{u^2z^3}{3!}\nonumber\\
      &=& \label{eq:T_1_ex}
      u z^2 + \frac{u^2z^3}{2}.
    \end{eqnarray}

    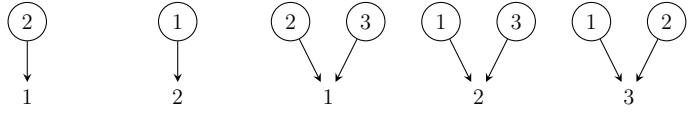
\begin{figure}
      \begin{center}
        \begin{tikzpicture}
          \foreach \h/\w/\n in {1/2/.7}{ % distance between rows, distance between columns, node size
            % size 2
            \foreach \i/\a/\b in {1/1/2, 2/2/1}{
              \draw(\i*\w, 0*\h) node[scale=\n] (node\a) {\a};
              \draw(\i*\w, 1*\h) node[draw, circle, scale=\n] (node\b) {\b};
              \draw[-stealth] (node\b) -- (node\a);
            }

            % size 3
            \foreach \i/\a/\b/\c in {3/1/2/3, 4/2/1/3, 5/3/1/2}{
              \draw(\i*\w, 0*\h) node[scale=\n] (node\a) {\a};
              \draw(\i*\w-.25*\w, 1*\h) node[draw, circle, scale=\n] (node\b) {\b};
              \draw(\i*\w+.25*\w, 1*\h) node[draw, circle, scale=\n] (node\c) {\c};
              \draw[-stealth] (node\b) -- (node\a);
              \draw[-stealth] (node\c) -- (node\a);
            }
          }
        \end{tikzpicture}
      \end{center}
      \caption{$\mathfrak{T}_{1}^{1,\{0,1,2\}}$}
      \label{T_ht1_k1_max2preimages}
    \end{figure}

    The following result gives the recursion for trees of bounded and unbounded height.

    \begin{lemma}
      \label{tree_recursion_formula}
      Let $k\geq 0$ and $\mathcal{P}\subseteq\mathbb{Z}_{\geq0}$, and work over a commutative ring $R\supseteq\mathbb{Q}$.
      The recursion for trees of bounded height is
      $$T_{\leq h}^{k,\mathcal{P}}(u,z)=\begin{cases}
      0 & h<0\\
      u z e_\mathcal{P}\big(T_{\leq h-1}^{k,\mathcal{P}}(u,z)\big) & 0\leq h<k\\
      z e_\mathcal{P}\big(T_{\leq h-1}^{k,\mathcal{P}}(u,z)\big) & 0=k\leq h\\
      z e_\mathcal{P}\big(T_{\leq h-1}^{k,\mathcal{P}}(u,z)\big) + (u-1)z e_\mathcal{P}\big(T_{\leq k-2}^{k,\mathcal{P}}(u,z)\big) & 0<k\leq h
      \end{cases}.$$
      The recursion for trees of all heights is
      $$T^{k,\mathcal{P}}(u,z)=\begin{cases}
      z e_\mathcal{P}\big(T^{k,\mathcal{P}}(u,z)\big) & 0=k\\
      z e_\mathcal{P}\big(T^{k,\mathcal{P}}(u,z)\big) + (u-1)z e_\mathcal{P}\big(T_{\leq k-2}^{k,\mathcal{P}}(u,z)\big)& 0<k
      \end{cases}.$$
    \end{lemma}
    \begin{proof}
      No rooted tree can have negative height, giving the first case of the first equality.

      If $0\leq h<k$, then the root has maximum distance from a leaf which is less than $k$, and hence is marked with both $u$ and $z$.
      Then the root has a set of subtrees subject to the preimage constraint $\mathcal{P}$; each subtree has height at most $h-1$.

      Similarly, if $k\leq h$, mark the root by $z$ and take a set of subtrees subject to the preimage constraint, where each has height at most $h-1$.
      But in the case where each has height at most $k-2$, the root will have height at most $k-1$; if $0<k$, this means the root has been incorrectly marked by not including the $u$, so subtract off those terms and add the corrected version back in.

      This last argument also works for unbounded trees.
      Namely, mark the root by $z$ and take a set of subtrees subject to the preimage constraint.
      But in the case where each has height at most $k-2$, the root will have height at most $k-1$; if $0<k$, this means the root has been incorrectly marked by not including the $u$, and these terms need to be corrected.
    \end{proof}

    The proof of Lemma \ref{tree_recursion_formula} handles the edge case where $0\notin\mathcal{P}$ just fine, but the results all turn out to be trivial.
    See Section \ref{examples}.

    In the case $k=0$, where there are no $u$-labels, we suppress the $k$ superscript to write
    \begin{eqnarray*}
      T_{\leq h}^{\mathcal{P}}(z) &=& T_{\leq h}^{0,\mathcal{P}}(u,z)\\
      \mathfrak{T}^\mathcal{P}&=&\mathfrak{T}^{0,\mathcal{P}}\\
      T^{\mathcal{P}}(z) &=& T^{0,\mathcal{P}}(u,z).
    \end{eqnarray*}
    In particular, as in the discussion prior to Lagrange Inversion Theorem \ref{lagrange_inversion_theorem}, $T^{\mathcal{P}}(z)$ is exactly $\psi_\mathcal{P}^{\compinv}$.

    The next result records the observation that taking $k=0$ gives the same functions one gets by plugging $u=1$ into the respective generating function defined for an arbitrary $k$.
    Notice that the second result is a version of Equation \ref{eq:unconstrained_tree} that allows for preimage constraints.

    \begin{theorem}
      \label{when_u_equals_1_tree}
      Let $h\geq 0$ and $\mathcal{P}\subseteq\mathbb{Z}_{\geq0}$, and work over a commutative ring $R\supseteq\mathbb{Q}$.
      Then
      \begin{eqnarray*}
        T_{\leq h}^{\mathcal{P}}(z) &=& z e_\mathcal{P}(T_{\leq h-1}^{\mathcal{P}}(z))\\
        T^{\mathcal{P}}(z)&=&z e_\mathcal{P}\big(T^{\mathcal{P}}(z)\big)\\
        T_{\leq h}^{\mathcal{P}}(z) &=& T_{\leq h}^{k,\mathcal{P}}(1,z)\\
        T^{\mathcal{P}}(z) &=& T^{k,\mathcal{P}}(1,z)\\
      \end{eqnarray*}
      for all $k\geq 0$.
    \end{theorem}
    \begin{proof}
      All four statements are immediate from Lemma \ref{tree_recursion_formula}.
      For the first two, take $k=0$.
      For the third, induct on $h$ and plug $u=1$ into the lemma; regardless of whether $h<k$ or $h\geq k$, the outcome is the same.
      For the fourth, plug in $u=1$ to see $T^{k,\mathcal{P}}(1,z)$ satisfies exactly the recursion that defines $T^{\mathcal{P}}(z)$; when $0\in\mathcal{P}$, use the uniqueness of compositional inverses in, say, Lagrange Inversion Theorem \ref{lagrange_inversion_theorem}, to see that the two functions must match; when $0\notin\mathcal{P}$, everything is 0.
    \end{proof}

    The following fact is a generalization of Equation \ref{eq:tree_coefficient} to the situation allowing preimage constraints.
    See Corollary \ref{F_Xi_coefficients} for an asymptotic estimate of these coefficients.

    \begin{corollary}
      \label{tree_coefficient}
      Let $\mathcal{P}\subseteq\mathbb{Z}_{\geq0}$ with $0\in\mathcal{P}$ and work over a commutative ring $R\supseteq\mathbb{Q}$.
      Then
      \begin{eqnarray*}
        \ \coeff{z^n}T^\mathcal{P}(z) &=& \frac{1}{n}\coeff{z^{n-1}}e_\mathcal{P}(z)^n
      \end{eqnarray*}
      for $n\in\mathbb{Z}_{>0}$.
    \end{corollary}
    \begin{proof}
      This is immediate from Lagrange Inversion Theorem \ref{lagrange_inversion_theorem} with $\epsilon=e_\mathcal{P}$ and $f(z)=z$.
    \end{proof}

    Consider a few sanity checks of Corollary \ref{tree_coefficient}.
    Taking $\mathcal{P}=\{0\}$ yields $e_\mathcal{P}(z)=1$ and $\coeff{z^n}T^\mathcal{P}(z) = \begin{cases} 1 & n=1\\
0 & n>1
    \end{cases}$, so
    $$T^{\{0\}}(z)=z,$$
    confirming the fact that the only rooted tree where every vertex has no preimages is the trivial one with a single node.
    Taking $\mathcal{P}=\{0,1\}$ yields $e_\mathcal{P}(z)=1+z$ and $\coeff{z^n}T^\mathcal{P}(z) = \frac{1}{n}n= \frac{n!}{n!}$, so
$$T^{\{0,1\}}(z)=\frac{z}{1-z},$$
    confirming the fact that there is only one rooted tree of size $n$ where each vertex has at most 1 preimage, and there are $n!$ ways to label such a graph (note that the root distinguishes one end of the graph from the other).
    Taking $\mathcal{P}=\mathbb{Z}_{\geq 0}$ gives the aforementioned specialization to Equation \ref{eq:tree_coefficient}.

    Corollary \ref{tree_coefficient} says that the number of rooted labeled trees subject to the preimage constraint $\mathcal{P}$ and of size $n$ is $n!\frac{1}{n}\coeff{z^{n-1}}e_\mathcal{P}(z)^n$.
    This is easy to compute for relatively small $n$, so Figure \ref{exact_tree_counts} shows the $\log_2$ of these counts for every subset $\mathcal{P}$ of $\{0,1,2,3,4\}$ which satisfies the hypotheses of Corollary \ref{tree_coefficient}.
    Anticipating the concept of periodicity that will arise repeatedly in the sections on singularity analysis, note that none of the curves for which $\gcd(\mathcal{P})>1$ show up in the plot.
    This is because such a tree must have a number of nodes that is congruent to 1 modulo $\gcd(\mathcal{P})$, and so most of the points on those curves are $-\infty$.
    See Figure \ref{exact_tree_counts_scatter} for a discrete version of the missing curves.

    \begin{figure}
      \begin{center}
        \includegraphics[width=7in]{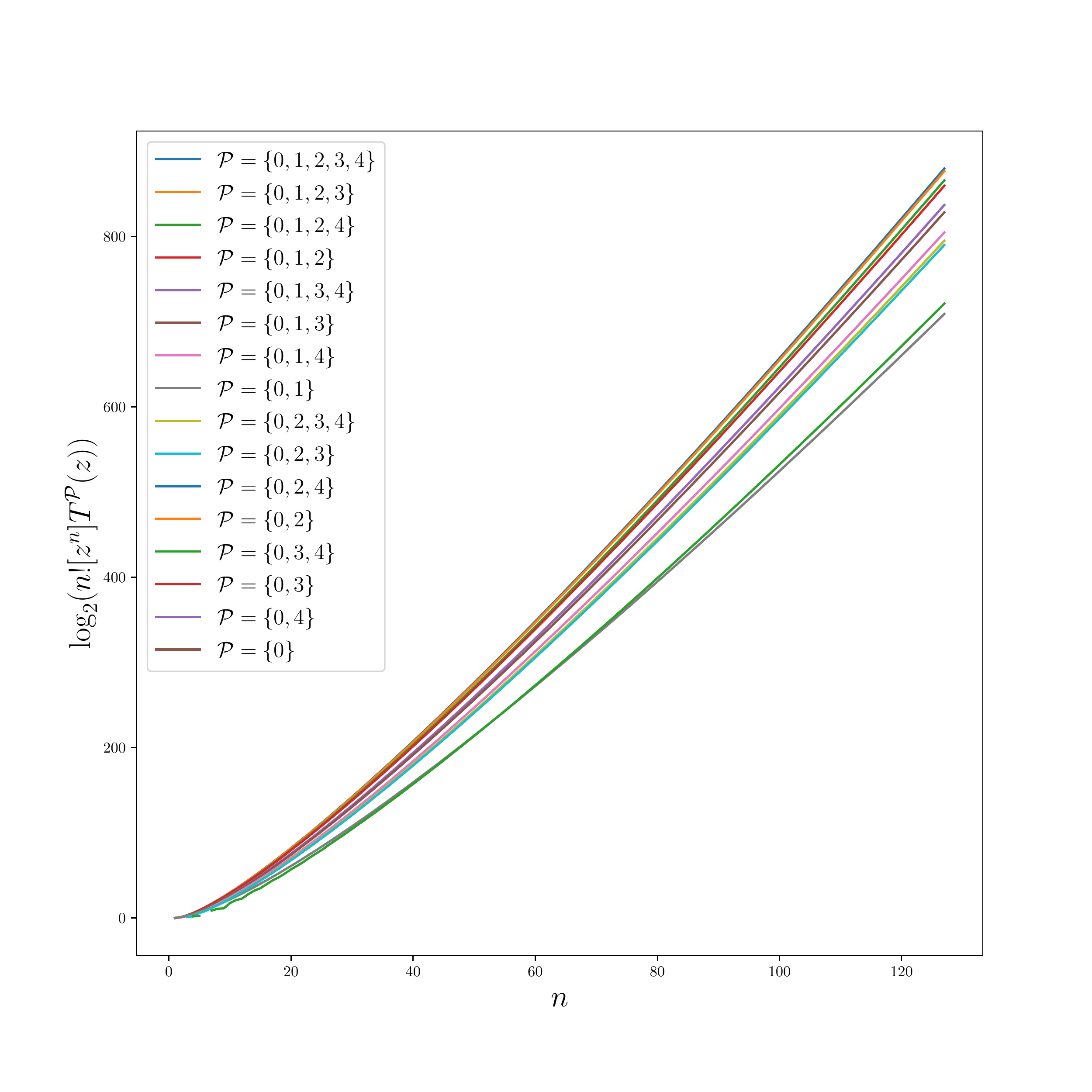}
      \end{center}
      \caption{$\log_2$ of Number of Trees of Size $n$}
      \label{exact_tree_counts}
    \end{figure}

    \begin{figure}
      \begin{center}
        \includegraphics[width=7in]{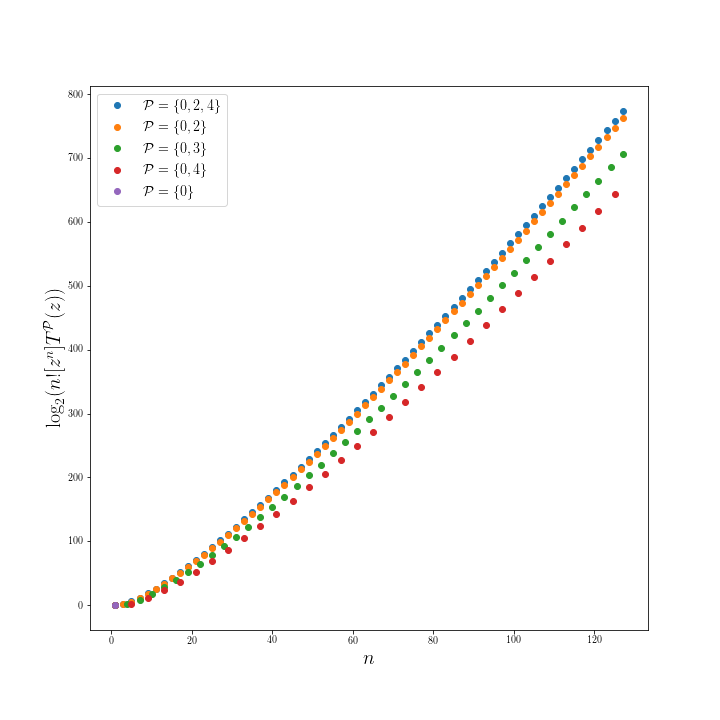}
      \end{center}
      \caption{$\log_2$ of Number of Trees of Size $n$ when $\gcd(\mathcal{P})>1$}
      \label{exact_tree_counts_scatter}
    \end{figure}
  \end{subsection}

  \begin{subsection}{Functions}
    \label{functions}

    While not quite as basic as rooted trees, the combinatorial class of functions subject to preimage constraints will appear in all of the constructions that follow.
    We could view it as a special case of any one of the function classes defined in the next three sections, but, following the choice for rooted trees, we approach it from the view of what is needed for Section \ref{image_combinatorics}.
    As with trees, we may then plug in $u=1$ to get the foundational, univariate version.

    For $k\geq 0$ and $\mathcal{P}\subseteq\mathbb{Z}_{\geq0}$, write $\mathfrak{F}^{k,\mathcal{P}}$ for the combinatorial class of functions $f:\set{n}\rightarrow \set{n}$, where $|f^{-1}(x)|\in\mathcal{P}$ for each $1\leq x\leq n$, where nodes are marked by $z$ and nodes not in $f^k(\set{n})$ are marked by $u$.
    Write $F^{k,\mathcal{P}}(u,z)$ for the generating function of $\mathfrak{F}^{k,\mathcal{P}}$.

    It is now easy to apply the generating function recipe to functions.

    \begin{theorem}
      \label{function_recursion_formula}
      Let $k\geq 0$ and $\mathcal{P}\subseteq\mathbb{Z}_{\geq0}$, and work over a commutative ring $R\supseteq\mathbb{Q}$.
      Then
      $$F^{k,\mathcal{P}}(u,z)=\Big(1-z e_{\mathcal{P}-1}\big(T^{k,\mathcal{P}}(u,z)\big)\Big)^{-1}.$$
    \end{theorem}
    \begin{proof}
      Every function is a set of cycles of rooted trees.
      In $\mathfrak{F}^{k,\mathcal{P}}$, the roots of those trees are subject to the preimage constraint $\mathcal{P}-1$, since one of their preimages will be accounted for on the cycle.
      Also, the roots of the trees are never marked by $u$, since they are in every iterated image.
      Each of the subtrees of the root are subject to the preimage constraint $\mathcal{P}$, and hence lie in $\mathfrak{T}^{k,\mathcal{P}}$.
      In other words,
      $$\mathfrak{F}^{k,\mathcal{P}}=\SET\Big(\CYC\big(\mathfrak{V}\star \SET_{\mathcal{P}-1}(\mathfrak{T}^{k,\mathcal{P}})\big)\Big),$$
      where $\mathfrak{V}$ denotes the class whose only object is the graph with a single vertex.
      Thus,
      $$F^{k,\mathcal{P}}(u,z)=e^{\ln\left(\left(1-z e_{\mathcal{P}-1}\left(T^{k,\mathcal{P}}(u,z)\right)\right)^{-1}\right)},$$
      and the result follows from Lemma \ref{e_identities}.
    \end{proof}

    In the case $k=0$, where there are no $u$-labels, we again suppress the $k$ superscript, to write
    \begin{eqnarray*}
      \mathfrak{F}^\mathcal{P}&=&\mathfrak{F}^{0,\mathcal{P}}\\
      F^{\mathcal{P}}(z) &=& F^{0,\mathcal{P}}(1,z).
    \end{eqnarray*}
    As before, taking $k=0$ gives the same functions you get by plugging $u=1$ into the respective generating function defined for an arbitrary $k$; this is pedantically noted in the next result.

    \begin{corollary}
      \label{when_u_equals_1_function}
      Let $\mathcal{P}\subseteq\mathbb{Z}_{\geq0}$, and work over a commutative ring $R\supseteq\mathbb{Q}$.
      Then
      \begin{eqnarray*}
        F^{\mathcal{P}}(z)&=&\Big(1-z e_{\mathcal{P}-1}\big(T^{\mathcal{P}}(z)\big)\Big)^{-1}\\
        F^{\mathcal{P}}(z) &=& F^{k,\mathcal{P}}(1,z)
      \end{eqnarray*}
      for all $k\geq 0$.
    \end{corollary}
    \begin{proof}
      This is immediate from Theorem \ref{function_recursion_formula} and Theorem \ref{when_u_equals_1_tree}.
    \end{proof}

    See Corollary \ref{F_Xi_coefficients} for an asymptotic estimate of the coefficients described in the next result.

    \begin{corollary}
      \label{function_coefficients}
      Let $\mathcal{P}\subseteq\mathbb{Z}_{\geq0}$ with $0\in\mathcal{P}$, and work over a commutative ring $R\supseteq\mathbb{Q}$.
      Then
      \begin{eqnarray*}
        F^{\mathcal{P}}(z)&=&z\frac{\partial}{\partial z}\ln T^\mathcal{P}(z)\\
        \ \coeff{z^n}F^\mathcal{P}(z) &=& \coeff{z^{n}}e_\mathcal{P}(z)^n.
      \end{eqnarray*}
    \end{corollary}
    \begin{proof}
      By Theorem \ref{when_u_equals_1_tree}, we have
      \begin{eqnarray*}
        T^\mathcal{P}(z) &=& z e_\mathcal{P}\big(T^{\mathcal{P}}(z)\big),
      \end{eqnarray*}
      so
      \begin{eqnarray*}
        \frac{\partial}{\partial z} T^\mathcal{P}(z) &=&  e_\mathcal{P}\big(T^{\mathcal{P}}(z)\big)+ z e_{\mathcal{P}-1}\big(T^{\mathcal{P}}(z)\big) \frac{\partial}{\partial z} T^\mathcal{P}(z)\\
        \frac{\partial}{\partial z} T^\mathcal{P}(z) \Big(1-z e_{\mathcal{P}-1}\big(T^{\mathcal{P}}(z)\big)\Big)&=&  e_\mathcal{P}\big(T^{\mathcal{P}}(z)\big)\\
        \frac{\partial}{\partial z} T^\mathcal{P}(z) &=&  e_\mathcal{P}\big(T^{\mathcal{P}}(z)\big)\Big(1-z e_{\mathcal{P}-1}\big(T^{\mathcal{P}}(z)\big)\Big)^{-1}.
      \end{eqnarray*}
      By Corollary \ref{when_u_equals_1_function}, this gives
      \begin{eqnarray*}
        \frac{\partial}{\partial z} T^\mathcal{P}(z)&=&  e_\mathcal{P}\big(T^{\mathcal{P}}(z)\big)F^{\mathcal{P}}(z)\\
        z \frac{\partial}{\partial z} T^\mathcal{P}(z)&=& z e_\mathcal{P}\big(T^{\mathcal{P}}(z)\big)F^{\mathcal{P}}(z)\\
        &=& T^{\mathcal{P}}(z)F^{\mathcal{P}}(z).
      \end{eqnarray*}
      Solving for $F^\mathcal{P}(z)$ gives the first result.

      For the second result, recall $\psi_\mathcal{P}=\frac{z}{e_\mathcal{P}(z)}$, so that $\psi_\mathcal{P}^{\compinv}=T^\mathcal{P}$.
      Then consider
      \begin{eqnarray}
        \coeff{z^n}F^\mathcal{P}(z) &=& \nonumber \coeff{z^n} z \frac{\frac{\partial}{\partial z} T^\mathcal{P}(z)}{T^\mathcal{P}(z)}\\
        &=& \nonumber \coeff{z^{-1}} z^{-n-1} z \frac{\frac{\partial}{\partial z} T^\mathcal{P}(z)}{T^\mathcal{P}(z)}\\
        &=& \label{eq:der_T}
        \coeff{z^{-1}} \Big(\psi_\mathcal{P}\big(T^\mathcal{P}(z)\big)\Big)^{-n} \frac{\frac{\partial}{\partial z} T^\mathcal{P}(z)}{T^\mathcal{P}(z)}.
      \end{eqnarray}
      Taking $\psi(z)=T^\mathcal{P}(z)$ and $f(z)=\psi_\mathcal{P}(z)^{-n}\frac{1}{z}$ in the Residue Composition Theorem \ref{residue_composition} simplifies Equation \ref{eq:der_T} to
      \begin{eqnarray*}
        \coeff{z^n}F^\mathcal{P}(z) &=& \val(T^\mathcal{P})\coeff{z^{-1}} \psi_\mathcal{P}(z)^{-n} \frac{1}{z}\\
        &=& \coeff{z^{-1}} \frac{e_\mathcal{P}(z)^n}{z^{n+1}}\\
        &=& \coeff{z^n} e_\mathcal{P}(z)^n.
      \end{eqnarray*}
    \end{proof}

    In the case $\mathcal{P}=\mathbb{Z}_{\geq 0}$, where there are no preimage constraints, the second result in Corollary \ref{function_coefficients} yields that
    \begin{eqnarray*}
      \coeff{z^n} F^{\mathbb{Z}_{\geq 0}}(z) &=& \coeff{z^n} e_{\mathbb{Z}_{\geq 0}}(z)^n\\
      &=& \coeff{z^n} e^{z n}\\
      &=& \frac{n^n}{n!},
    \end{eqnarray*}
    confirming the well-known fact that there are $n^n$ functions from $\set{n}$ to itself.
    Or taking $\mathcal{P}=\{0,1\}$ gives $\coeff{z^n}F^{\{0,1\}}(z)=1=\frac{n!}{n!}$ and $F^{\{0,1\}}(z)=\frac{1}{1-z}$, confirming that there are $n!$ permutations on $\set{n}$.
    Finally, taking $\mathcal{P}=\{0\}$ gives $F^{\{0\}}(z)=1$ and the rather boring observation that the only function for which every point has no preimages is the trivial function from the empty set to itself.

    Since Corollary \ref{function_coefficients} says that the number of functions from $\set{n}$ to $\set{n}$ subject to the preimage constraint $\mathcal{P}$ is $n!\coeff{z^{n}}e_\mathcal{P}(z)^n$, and since this is easy to compute for relatively small $n$, Figure \ref{exact_function_counts} shows the $\log_2$ of these counts for every subset $\mathcal{P}$ of $\{0,1,2,3,4\}$ which satisfies the hypotheses of Corollary \ref{function_coefficients}.
    In another anticipation of the concept of periodicity, note such a function must have number of nodes divisible by $\gcd(\mathcal{P})$, and so most of the points on those curves are $-\infty$.
    See Figure \ref{exact_function_counts_scatter} for a discrete version of the missing curves.

    \begin{figure}
      \begin{center}
        \includegraphics[width=7in]{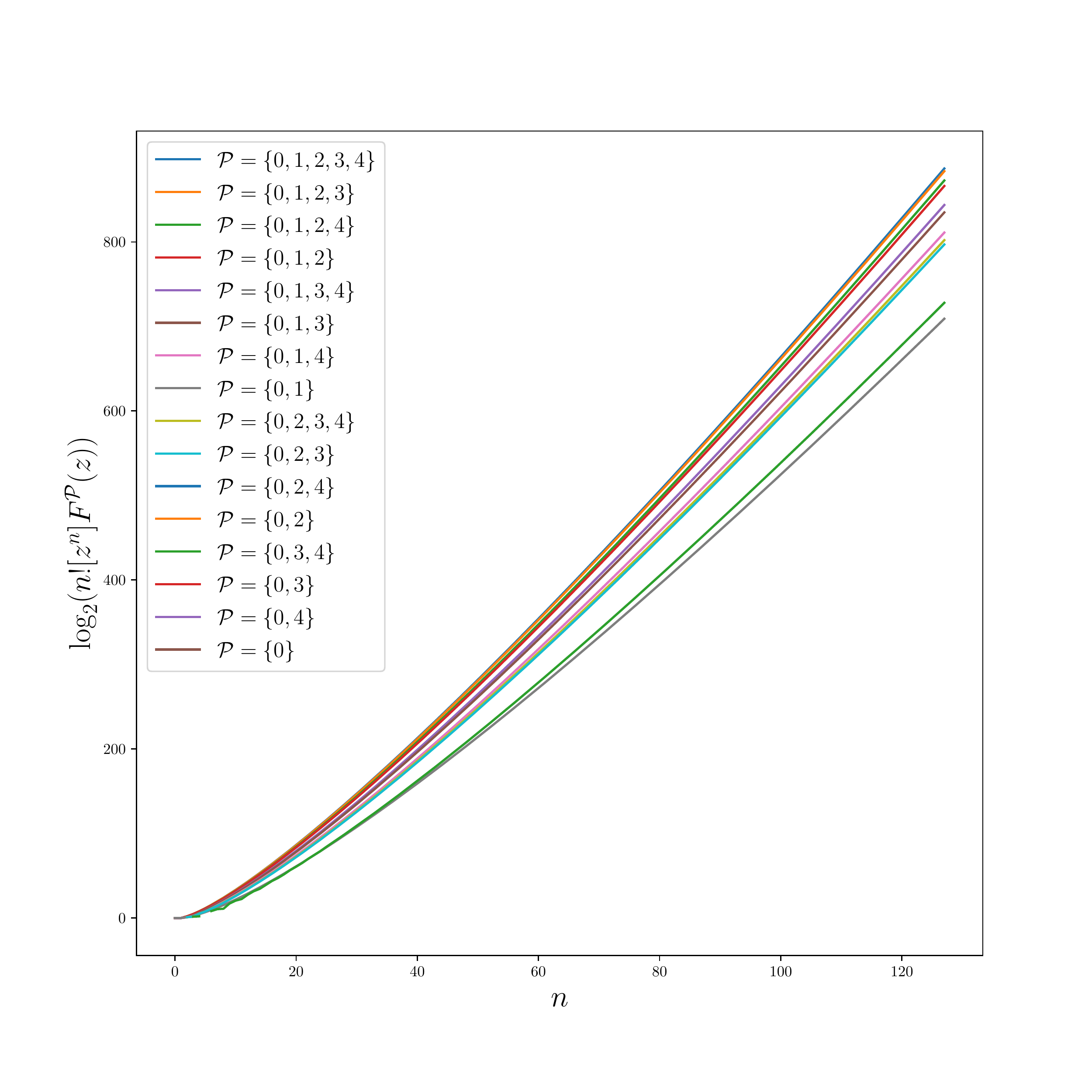}
      \end{center}
      \caption{$\log_2$ of Number of Functions on $\set{n}$}
      \label{exact_function_counts}
    \end{figure}

    \begin{figure}
      \begin{center}
        \includegraphics[width=7in]{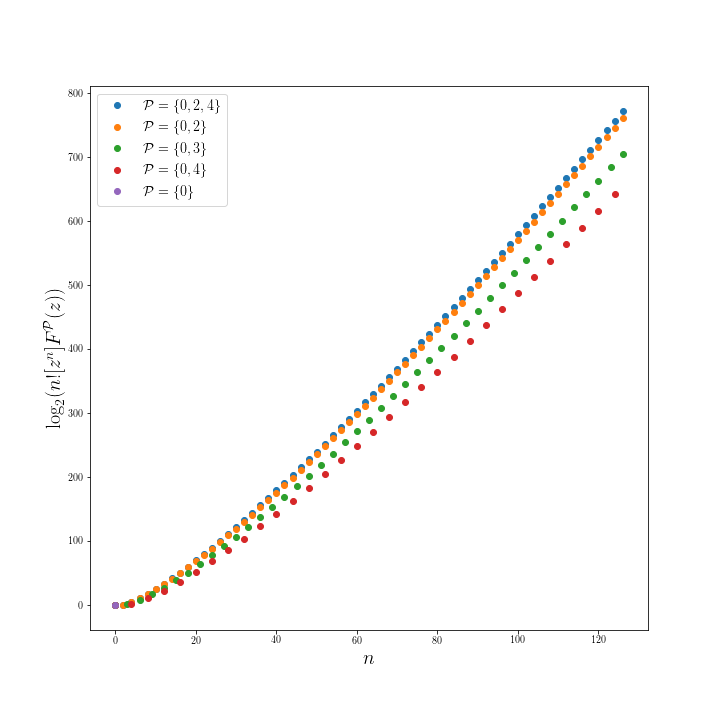}
      \end{center}
      \caption{$\log_2$ of Number of Functions on $\set{n}$ when $\gcd(\mathcal{P})>1$}
      \label{exact_function_counts_scatter}
    \end{figure}
  \end{subsection}

  \begin{subsection}{Components}

    For $\mathcal{P}\subseteq\mathbb{Z}_{\geq0}$, write $\mathfrak{F}^{\COMP,\mathcal{P}}$ for the combinatorial class of functions $f:\set{n}\rightarrow \set{n}$, where $|f^{-1}(x)|\in\mathcal{P}$ for each $1\leq x\leq n$, where nodes are marked by $z$ and connected components are marked by $u$.
    Write $F^{\COMP,\mathcal{P}}(u,z)$ for the generating function of $\mathfrak{F}^{\COMP,\mathcal{P}}$.

    Write $\mathfrak{F}^{\CONNECT,\mathcal{P}}$ for the result of taking the combinatorial subclass of $\mathfrak{F}^{\COMP,\mathcal{P}}$ consisting only of those functions whose graph is connected, and since they all have exactly one component, removing the parameter function that counts components.
    Write $F^{\CONNECT,\mathcal{P}}(z)$ for the corresponding generating function.

    \begin{lemma}
      \label{component_function_formula}
      Let $\mathcal{P}\subseteq\mathbb{Z}_{\geq0}$, and work over a commutative ring $R\supseteq\mathbb{Q}$.
      Then
      \begin{eqnarray*}
        F^{\CONNECT,\mathcal{P}}(z)&=&\ln F^\mathcal{P}(z)\\
        F^{\COMP,\mathcal{P}}(u,z)&=&e^{u F^{\CONNECT, \mathcal{P}}(z)}.
      \end{eqnarray*}
    \end{lemma}
    \begin{proof}
      Every function is a set of cycles of rooted trees.
      In $\mathfrak{F}^{\COMP,\mathcal{P}}$, the roots of those trees are subject to the preimage constraint $\mathcal{P}-1$, since one of their preimages is already accounted for on the cycle.
      Each of the subtrees of the root are subject to the preimage constraint $\mathcal{P}$, and hence lie in $\mathfrak{T}^{\mathcal{P}}$.
      The generating functions for the cycles of trees needs a single divisor $u$ to mark the cycle; from the viewpoint of a combinatorial class, this can be accomplished via a dummy node $U$ of size 0 which the parameter function marks as having cycle-size 1.
      In other words, $\mathfrak{F}^{\COMP,\mathcal{P}}=\SET\Big(\mathfrak{U}\star\CYC\big(\mathfrak{V}\star \SET_{\mathcal{P}-1}(\mathfrak{T}^{\mathcal{P}})\big)\Big)$, where $\mathfrak{V}$ denotes the class whose only object is the graph with a single vertex and $\mathfrak{U}=\{U\}$ is the class whose only object is the dummy node.
      That is,
      \begin{eqnarray}
        F^{\COMP,\mathcal{P}}(u,z)&=&\nonumber
        e^{u\ln\left(\left(1-z e_{\mathcal{P}-1}\left(T^{\mathcal{P}}(z)\right)\right)^{-1}\right)}.
      \end{eqnarray}
      By Corollary \ref{when_u_equals_1_function}, this is
      \begin{eqnarray}
        F^{\COMP,\mathcal{P}}(u,z)&=& \label{eq:component_function} e^{u \ln F^\mathcal{P}(z)}.
      \end{eqnarray}

      The first statement is now immediate from the observation that
      $$F^{\CONNECT,\mathcal{P}}(z)=\coeff{u^1}F^{\COMP,\mathcal{P}}(u,z).$$
      The second statement follows by substituting the first into Equation \ref{eq:component_function}.
    \end{proof}

    For $\mathcal{P}\subseteq\mathbb{Z}_{\geq0}$, define
    $$\Xi^{\COMP,\mathcal{P}}(z)=\Big(\frac{\partial}{\partial u}F^{\COMP,\mathcal{P}}(u,z)\Big)\mid_{u=1}.$$
    Note that $\frac{\partial}{\partial u}u^i z^n|_{u=1}=i z^n$, so $\coeff{z^n}\Xi^{\COMP,\mathcal{P}}(z)$ is the number of components of $f$, divided by $n!$, as one sums across all $f\in\mathfrak{F}^{\COMP, \mathcal{P}}$.

    \begin{theorem}
      \label{component_Xi_formula}
      Let $\mathcal{P}\subseteq\mathbb{Z}_{\geq0}$, and work over a commutative ring $R\supseteq\mathbb{Q}$.
      Then
      $$\Xi^{\COMP,\mathcal{P}}(z)=F^{\mathcal{P}}(z)\ln F^{\mathcal{P}}(z).$$
    \end{theorem}
    \begin{proof}
      Plug Lemma \ref{component_function_formula} into the definition of $\Xi^{\COMP,\mathcal{P}}$ and simplify, using Lemmas \ref{basic_derivatives} and \ref{e_identities} as appropriate, to see that
      \begin{eqnarray*}
        \Xi^{\COMP,\mathcal{P}}(z)&=&e^{u\ln F^\mathcal{P}(z)}\ln F^\mathcal{P}(z)\mid_{u=1}\\
        &=&F^\mathcal{P}(z)\ln F^\mathcal{P}(z).
      \end{eqnarray*}
    \end{proof}
  \end{subsection}

  \begin{subsection}{Cycle Points}
    \label{cycle_points}

    For $\mathcal{P}\subseteq\mathbb{Z}_{\geq0}$, write $\mathfrak{F}^{\CYC,\mathcal{P}}$ for the combinatorial class of functions $f:\set{n}\rightarrow \set{n}$, where $|f^{-1}(x)|\in\mathcal{P}$ for each $1\leq x \leq n$, where nodes are marked by $z$ and nodes that lie on a cycle are marked by $u$.
    Write $F^{\CYC,\mathcal{P}}(u,z)$ for the generating function of $\mathfrak{F}^{\CYC,\mathcal{P}}$.

    \begin{lemma}
      \label{cyclic_function_formula}
      Let $\mathcal{P}\subseteq\mathbb{Z}_{\geq0}$, and work over a commutative ring $R\supseteq\mathbb{Q}$.
      Then
      $$F^{\CYC,\mathcal{P}}(u,z)=\Big(1-u z e_{\mathcal{P}-1}\big(T^{\mathcal{P}}(z)\big)\Big)^{-1}.$$
    \end{lemma}
    \begin{proof}
      Every function is a set of cycles of rooted trees.
      In $\mathfrak{F}^{\CYC,\mathcal{P}}$, the roots of those trees are subject to the preimage constraint $\mathcal{P}-1$, since one of their preimages is already accounted for on the cycle.
      Also, the roots of the trees are exactly the points marked by $u$, since they are on a cycle; from the viewpoint of a combinatorial class, this can be accomplished via a dummy node $U$ of size 0 which the parameter function marks with 1 to denote being on a cycle.
      In other words, if $\mathfrak{V}$ denotes the class whose only object is the graph with a single vertex and $\mathfrak{U}=\{U\}$ is the class whose only object is the dummy node, then the combinatorial class for a node not on a cycle is again $\mathfrak{V}$ and the class for a node on a cycle is $\mathfrak{U}\star\mathfrak{V}$.
      Each of the subtrees of the root are subject to the preimage constraint $\mathcal{P}$, and hence lie in $\mathfrak{T}^{\mathcal{P}}$.

      Putting this all together, $\mathfrak{F}^{\CYC,\mathcal{P}}=\SET\left(\CYC\left(\left(\mathfrak{U}\star\mathfrak{V}\right)\star \SET_{\mathcal{P}-1}(\mathfrak{T}^{\mathcal{P}})\right)\right)$ and
      $$F^{\CYC,\mathcal{P}}(u,z)=e^{\ln\left(\left(1-u z e_{\mathcal{P}-1}\left(T^{\mathcal{P}}(z)\right)\right)^{-1}\right)},$$
      and the result follows from Lemma \ref{e_identities}.
    \end{proof}

    For $\mathcal{P}\subseteq\mathbb{Z}_{\geq0}$, define
    $$\Xi^{\CYC,\mathcal{P}}(z)=\Big(\frac{\partial}{\partial u}F^{\CYC,\mathcal{P}}(u,z)\Big)\mid_{u=1}.$$
    Note that $\frac{\partial}{\partial u}u^i z^n|_{u=1}=i z^n$, so the coefficient of $z^n$ in $\Xi^{\CYC,\mathcal{P}}(z)$ is the number of cyclic nodes in $f$, divided by $n!$, as one sums across all $f\in\mathfrak{F}^{\CYC, \mathcal{P}}$.

    \begin{theorem}
      \label{cyclic_Xi_formula}
      Let $\mathcal{P}\subseteq\mathbb{Z}_{\geq0}$, and work over a commutative ring $R\supseteq\mathbb{Q}$.
      Then
      $$\Xi^{\CYC,\mathcal{P}}(z)=z e_{\mathcal{P}-1}\big(T^{\mathcal{P}}(z)\big)F^{\mathcal{P}}(z)^2.$$
    \end{theorem}
    \begin{proof}
      Plug Lemma \ref{cyclic_function_formula} into the definition of $\Xi^{\CYC,\mathcal{P}}$ and simplify, using Lemmas \ref{basic_derivatives} and \ref{e_identities} as appropriate, to see that
      \begin{eqnarray*}
        \Xi^{\CYC,\mathcal{P}}(z)&=&\Big(\frac{\partial}{\partial u}F^{\CYC,\mathcal{P}}(u,z)\Big)\mid_{u=1}\\
        &=&-\Big(1-u z e_{\mathcal{P}-1}\big(T^{\mathcal{P}}(z)\big)\Big)^{-2}\Big(-z e_{\mathcal{P}-1}\big(T^{\mathcal{P}}(z)\big)\Big)\mid_{u=1}\\
        &=&\frac{z e_{\mathcal{P}-1}\big(T^{\mathcal{P}}(z)\big)}{\Big(1-z e_{\mathcal{P}-1}\big(T^{\mathcal{P}}(z)\big)\Big)^{2}}.
      \end{eqnarray*}
      The result now follows from Corollary \ref{when_u_equals_1_function}.
    \end{proof}

    An asymptotic estimate of the coefficients of $\Xi^{\CYC,\mathcal{P}}(z)$ appears in Corollary \ref{F_Xi_coefficients}.
  \end{subsection}

  \begin{subsection}{Partial Functions}
    Since the background is already in place, it is a straightforward task to apply this machinery to partial functions subject to preimage constraints.

    Recall that a partial function on $\set{n}$ is a function from some subset of $S\subseteq\set{n}$ to $\set{n}$.
    The elements of $\set{n}\setminus S$ are not mapped anywhere; thus, in the graph view of functions, partial functions are the result of relaxing the condition that every node have an out-edge to the condition that every node have at most one out-edge.

    For $k\geq 0$ and $\mathcal{P}\subseteq\mathbb{Z}_{\geq0}$, write $\mathfrak{P}^{k,\mathcal{P}}$ for the combinatorial class of partial functions $f$ from $\set{n}\rightarrow \set{n}$, where $|f^{-1}(x)|\in\mathcal{P}$ for each $1\leq x\leq n$, where nodes are marked by $z$ and nodes not in $f^k(\set{n})$ are marked by $u$.
    Write $P^{k,\mathcal{P}}(u,z)$ for the generating function of $\mathfrak{P}^{k,\mathcal{P}}$.

    \begin{theorem}
      \label{partial_function_recursion_formula}
      Let $k\geq 0$ and $\mathcal{P}\subseteq\mathbb{Z}_{\geq0}$, and work over a commutative ring $R\supseteq\mathbb{Q}$.
      Then
      $$P^{k,\mathcal{P}}(u,z)=F^{k,\mathcal{P}}(u,z)e^{T^{k,\mathcal{P}}(u,z)}.$$
    \end{theorem}
    \begin{proof}
      The components of a partial function are either rooted trees or cycles of rooted trees.
      Like in $\mathfrak{F}^{k,\mathcal{P}}$, the roots of those trees that lie on cycles are subject to the preimage constraint $\mathcal{P}-1$, since one of their preimages will be accounted for on the cycle, and the roots of the trees are never marked by $u$, since they are in every iterated image.
      The trees that are not on cycles (including the proper subtrees of the roots on the cycle) are subject to the preimage constraint $\mathcal{P}$, and hence lie in $\mathfrak{T}^{k,\mathcal{P}}$.
      In other words,
      $$\mathfrak{P}^{k,\mathcal{P}}=\SET\Big(\CYC\big(\mathfrak{V}\star \SET_{\mathcal{P}-1}(\mathfrak{T}^{k,\mathcal{P}})\big) \cup \mathfrak{T}^{k,\mathcal{P}}\Big),$$
      where $\mathfrak{V}$ denotes the class whose only object is the graph with a single vertex.
      Thus,
      $$P^{k,\mathcal{P}}(u,z)=e^{\ln\left(\left(1-z e_{\mathcal{P}-1}\left(T^{k,\mathcal{P}}(u,z)\right)\right)^{-1}\right) + T^{k,\mathcal{P}}(u,z)},$$
      and the result follows from Lemma \ref{e_identities} and Theorem \ref{function_recursion_formula}.
    \end{proof}

    In the case $k=0$, where there are no $u$-labels, we again suppress the $k$ superscript, to write
    \begin{eqnarray*}
      \mathfrak{P}^\mathcal{P}&=&\mathfrak{P}^{0,\mathcal{P}}\\
      P^{\mathcal{P}}(z) &=& P^{0,\mathcal{P}}(1,z).
    \end{eqnarray*}
    As before, taking $k=0$ gives the same functions you get by plugging $u=1$ into the respective generating function defined for an arbitrary $k$; this is pedantically noted in the next result.

    \begin{corollary}
      \label{when_u_equals_1_partial_function}
      Let $\mathcal{P}\subseteq\mathbb{Z}_{\geq0}$, and work over a commutative ring $R\supseteq\mathbb{Q}$.
      Then
      \begin{eqnarray*}
        P^{\mathcal{P}}(z)&=& F^{\mathcal{P}}(z) e^{T^{\mathcal{P}}(z)}\\
        P^{\mathcal{P}}(z) &=& P^{k,\mathcal{P}}(1,z)
      \end{eqnarray*}
      for all $k\geq 0$.
    \end{corollary}
    \begin{proof}
      This is immediate from Theorem \ref{partial_function_recursion_formula}, Theorem \ref{when_u_equals_1_tree}, and Theorem \ref{when_u_equals_1_tree}.
    \end{proof}

    See Corollary \ref{F_Xi_coefficients} for an asymptotic estimate of the coefficients described in the next result.

    \begin{corollary}
      \label{partial_function_coefficients}
      Let $\mathcal{P}\subseteq\mathbb{Z}_{\geq0}$ with $0\in\mathcal{P}$, and work over a commutative ring $R\supseteq\mathbb{Q}$.
      Then
      \begin{eqnarray*}
        \ \coeff{z^n}P^\mathcal{P}(z) &=& \coeff{z^{n}}e_\mathcal{P}(z)^n e^z.
      \end{eqnarray*}
    \end{corollary}
    \begin{proof}
      Recall $\psi_\mathcal{P}=\frac{z}{e_\mathcal{P}(z)}$, so that $\psi_\mathcal{P}^{\compinv}=T^\mathcal{P}$.
      First Corollary \ref{when_u_equals_1_partial_function} and then Corollary \ref{function_coefficients} together give that
      \begin{eqnarray}
        \coeff{z^n}P^\mathcal{P}(z) &=& \nonumber \coeff{z^n} z \frac{\frac{\partial}{\partial z} T^\mathcal{P}(z)}{T^\mathcal{P}(z)} e^{T^\mathcal{P}(z)}\\
        &=& \nonumber \coeff{z^{-1}} z^{-n-1} z \frac{\frac{\partial}{\partial z} T^\mathcal{P}(z)}{T^\mathcal{P}(z)}e^{T^\mathcal{P}(z)}\\
        &=& \label{eq:der_P}
        \coeff{z^{-1}} \Big(\psi_\mathcal{P}\big(T^\mathcal{P}(z)\big)\Big)^{-n} \frac{\frac{\partial}{\partial z} T^\mathcal{P}(z)}{T^\mathcal{P}(z)}e^{T^\mathcal{P}(z)}.
      \end{eqnarray}
      Taking $\psi(z)=T^\mathcal{P}(z)$ and $f(z)=\psi_\mathcal{P}(z)^{-n}\frac{e^z}{z}$ in the Residue Composition Theorem \ref{residue_composition} simplifies Equation \ref{eq:der_P} to
      \begin{eqnarray*}
        \coeff{z^n}F^\mathcal{P}(z) &=& \val(T^\mathcal{P})\coeff{z^{-1}} \psi_\mathcal{P}(z)^{-n} \frac{e^z}{z}\\
        &=& \coeff{z^{-1}} \frac{e_\mathcal{P}(z)^n e^z}{z^{n+1}}\\
        &=& \coeff{z^n} e_\mathcal{P}(z)^n e^z.
      \end{eqnarray*}
    \end{proof}

    In the case $\mathcal{P}=\mathbb{Z}_{\geq 0}$, where there are no preimage constraints, Corollary \ref{partial_function_coefficients} yields that
    \begin{eqnarray*}
      \coeff{z^n} P^{\mathbb{Z}_{\geq 0}}(z) &=& \coeff{z^n} e_{\mathbb{Z}_{\geq 0}}(z)^n e^z\\
      &=& \coeff{z^n} e^{(n+1)z}\\
      &=& \frac{(n+1)^n}{n!},
    \end{eqnarray*}
    yielding that there are $(n+1)^n$ partial functions from $\set{n}$ to itself; this is consistent with the observation that, for a given partial function, each node has $n+1$ possible outcomes: it is either mapped to one of the $n$ nodes, or it is not mapped anywhere.
    Or taking $\mathcal{P}=\{0,1\}$ gives
    $$\coeff{z^n}P^{\{0,1\}}(z)=\frac{\sum_{k=o}^n\binom{n}{k}\frac{n!}{(n-k)!}}{n!};$$
    this is consistent with letting $k$ be the number of nodes that have a single preimage, so there are $\frac{n!}{(n-k)!}$ ways to assign preimages.
    Finally, taking $\mathcal{P}=\{0\}$ gives
    $$[z^n] P^{\{0\}}(z)=\frac{1}{n!},$$
    consistent with the observation that there is a unique partial function on $n$ nodes that does not map any element to any other.

    Corollary \ref{partial_function_coefficients} says that the number of partial functions from $\set{n}$ to $\set{n}$ subject to the preimage constraint $\mathcal{P}$ is $n!\coeff{z^{n}}e_\mathcal{P}(z)^n$.
    This is easy to compute for relatively small $n$, so Figure \ref{exact_partial_function_counts} shows the $\log_2$ of these counts for every subset $\mathcal{P}$ of $\{0,1,2,3,4\}$ which satisfies the hypotheses of Corollary \ref{partial_function_coefficients}.

    \begin{figure}
      \begin{center}
        \includegraphics[width=7in]{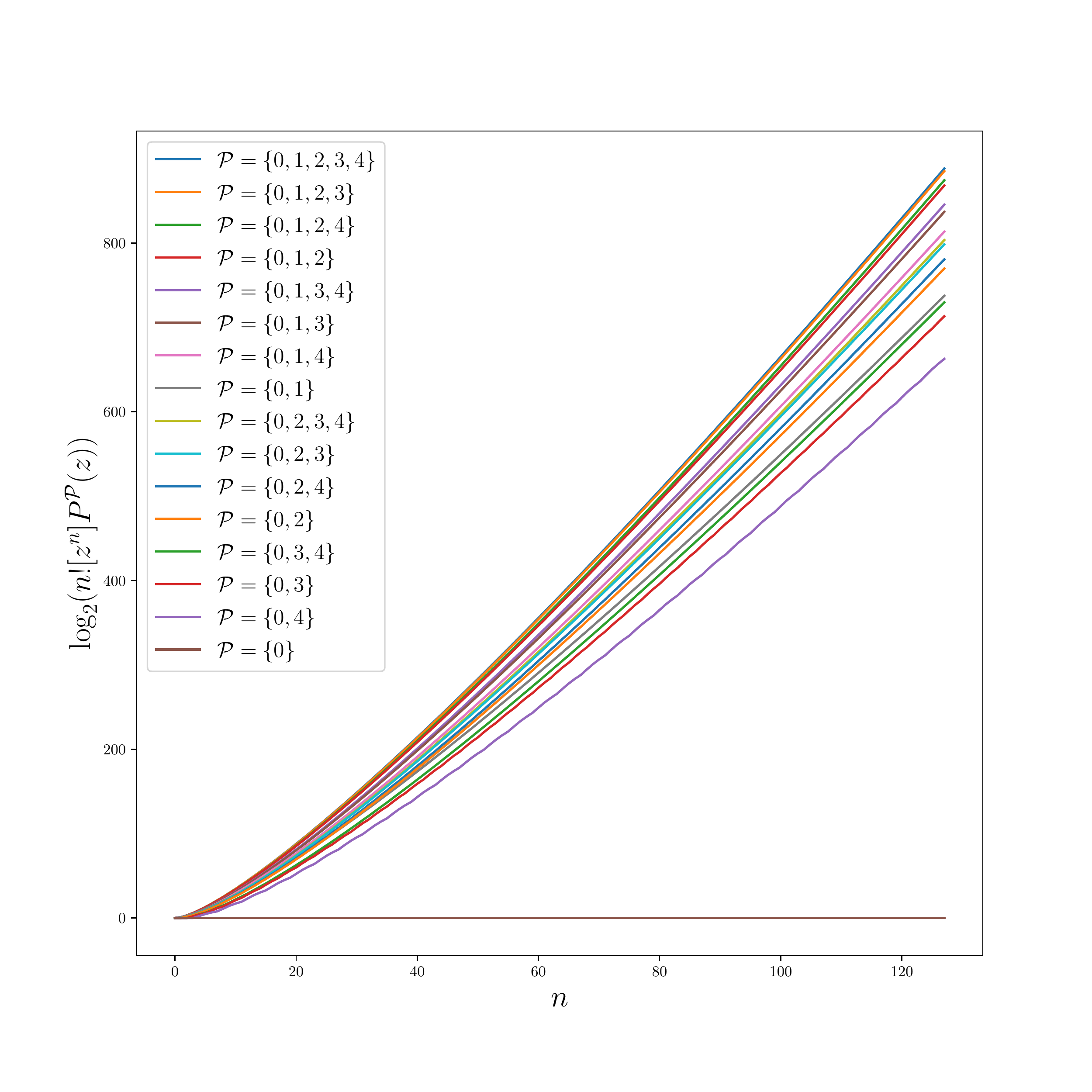}
      \end{center}
      \caption{$\log_2$ of Number of Partial Functions on $\set{n}$}
      \label{exact_partial_function_counts}
    \end{figure}
  \end{subsection}
  
  \begin{subsection}{Iterated Images}
    \label{image_combinatorics}

    In order to examine the $k$th image of a function, we need a class that marks all vertices with $z$ and that marks vertices that are too close to leaves (and hence not in the $k$th image) by a variable $u$.
    This was defined in Section \ref{trees}, iteratively building things up in such a way that we could keep track of the height of trees as we do the recursion.

    Before getting to the main combinatorial result, we record a technical observation.

    \begin{lemma}
      \label{tree_tech2}
      Let $k\geq 0$, $h\geq 0$, and $\mathcal{P}\subseteq\mathbb{Z}_{\geq0}$, and work over a commutative ring $R\supseteq\mathbb{Q}$.
      Then
      \begin{eqnarray*}
        \Big(\frac{\partial}{\partial u}T^{k,\mathcal{P}}(u,z)\Big)\mid_{u=1} &=& T_{\leq k-1}^{\mathcal{P}}(z)F^\mathcal{P}(z).
      \end{eqnarray*}
    \end{lemma}
    \begin{proof}
      When $k=0$, the result simplifies to the claim $0=0$, so suppose $k>0$.

      Take the derivative with respect to $u$ of the result in Lemma \ref{tree_recursion_formula}, using Lemmas \ref{basic_derivatives} and \ref{e_identities} as necessary, to see that
      \begin{eqnarray*}
        \frac{\partial}{\partial u}T^{k,\mathcal{P}}(u,z) &=& z e_{\mathcal{P}-1}\big(T^{k,\mathcal{P}}(u,z)\big)\frac{\partial}{\partial u}T^{k,\mathcal{P}}(u,z)+\Big(z e_\mathcal{P}\big(T_{\leq k-2}^{k,\mathcal{P}}(u,z)\big)\\
        && +(u-1)z e_{\mathcal{P}-1}\big(T_{\leq k-2}^{k,\mathcal{P}}(u,z)\big)\frac{\partial}{\partial u}T_{\leq k-2}^{k,\mathcal{P}}(u,z)\Big).
      \end{eqnarray*}
      Plug in $u=1$ and apply Theorem \ref{when_u_equals_1_tree} to get
      \begin{eqnarray*}
        \Big(\frac{\partial}{\partial u}T^{k,\mathcal{P}}(u,z)\Big)\mid_{u=1} &=& z e_{\mathcal{P}-1}\big(T^{\mathcal{P}}(z)\big)\Big(\frac{\partial}{\partial u}T^{k,\mathcal{P}}(u,z)\Big)|_{u=1}+z e_\mathcal{P}\big(T_{\leq k-2}^{\mathcal{P}}(z)\big),
      \end{eqnarray*}
      so
      \begin{eqnarray*}
        \Big(\frac{\partial}{\partial u}T^{k,\mathcal{P}}(u,z)\Big)\mid_{u=1}\Big(1-z e_{\mathcal{P}-1}\big(T^{\mathcal{P}}(z)\big)\Big) &=& z e_\mathcal{P}\big(T_{\leq k-2}^{\mathcal{P}}(z)\big)\\
        \Big(\frac{\partial}{\partial u}T^{k,\mathcal{P}}(u,z)\Big)\mid_{u=1} &=& \frac{z e_\mathcal{P}\big(T_{\leq k-2}^{\mathcal{P}}(z)\big)}{\Big(1-z e_{\mathcal{P}-1}\big(T^{\mathcal{P}}(z)\big)\Big)}.
      \end{eqnarray*}
      Taking $h=k-1$ in Theorem \ref{when_u_equals_1_tree} simplifies this to
      \begin{eqnarray*}
        \Big(\frac{\partial}{\partial u}T^{k,\mathcal{P}}(u,z)\Big)\mid_{u=1} &=& \frac{T_{\leq k-1}^{\mathcal{P}}(z)}{\Big(1-z e_{\mathcal{P}-1}\big(T^{\mathcal{P}}(z)\big)\Big)},
      \end{eqnarray*}
      and the result now follows from Corollary \ref{when_u_equals_1_function}.
    \end{proof}

    For $k\geq 0$ and $\mathcal{P}\subseteq\mathbb{Z}_{\geq0}$, define
    $$\Xi^{k,\mathcal{P}}(z)=\Big(\frac{\partial}{\partial u}F^{k,\mathcal{P}}(u,z)\Big)\mid_{u=1}.$$
    Note that $\frac{\partial}{\partial u}u^i z^n|_{u=1}=i z^n$, so the coefficient of $z^n$ in $\Xi^{k,\mathcal{P}}(z)$ is the number of nodes not in $f^k(\set{n})$, divided by $n!$, as one sums across all $f\in\mathfrak{F}^{\mathcal{P}}$.
    Similarly, define
    $$\Xi^{{\rm partial}, k,\mathcal{P}}(z)=\Big(\frac{\partial}{\partial u}P^{k,\mathcal{P}}(u,z)\Big)\mid_{u=1},$$
    noting it is the analog of $\Xi^{k,\mathcal{P}}$ for partial functions.

    \begin{theorem}
      \label{Xi_formula}
      Let $k\geq 0$ and $\mathcal{P}\subseteq\mathbb{Z}_{\geq0}$, and work over a commutative ring $R\supseteq\mathbb{Q}$. Then
      \begin{eqnarray*}
        \Xi^{k,\mathcal{P}}(z) &=& z T_{\leq k-1}^{\mathcal{P}}(z) e_{\mathcal{P}-2}\big(T^{\mathcal{P}}(z)\big)F^{\mathcal{P}}(z)^3\\
        \Xi^{{\rm partial}, k,\mathcal{P}}\left(z\right) &=& T_{\leq k-1}^{\mathcal{P}}\left(z\right) F^{\mathcal{P}}\left(z\right)^2 e^{T^{\mathcal{P}}\left(z\right)} \Big(z e_{\mathcal{P}-2}\left(T^{\mathcal{P}}\left(z\right)\right)F^{\mathcal{P}}\left(z\right)+1\Big).
      \end{eqnarray*}
    \end{theorem}
    \begin{proof}
      Plug Theorem \ref{function_recursion_formula} into the definition of $\Xi^{k,\mathcal{P}}$ and simplify, using Lemmas \ref{basic_derivatives}, \ref{e_identities}, and \ref{tree_tech2} and Theorem \ref{when_u_equals_1_tree}, as appropriate, to see that
      \begin{eqnarray*}
        \Xi^{k,\mathcal{P}}(z)&=&\Big(\frac{\partial}{\partial u}F^{k,\mathcal{P}}(u,z)\Big)\mid_{u=1}\\
        &=&\left(\frac{\partial}{\partial u}\left(1-z e_{\mathcal{P}-1}\left(T^{k,\mathcal{P}}(u,z)\right)\right)^{-1}\right)\mid_{u=1}\\
        &=&\left(-\left(1-z e_{\mathcal{P}-1}\left(T^{k,\mathcal{P}}(u,z)\right)\right)^{-2}(-z) \frac{\partial}{\partial u}e_{\mathcal{P}-1}\left(T^{k,\mathcal{P}}(u,z)\right)\right)\mid_{u=1}\\
        &=&\Big(1-z e_{\mathcal{P}-1}\big(T^{\mathcal{P}}(z)\big)\Big)^{-2}z\Big(e_{\mathcal{P}-2}\big(T^{k,\mathcal{P}}(u,z)\big)\frac{\partial}{\partial u}T^{k,\mathcal{P}}(u,z)\big)\Big)\mid_{u=1}\\
        &=&\frac{z e_{\mathcal{P}-2}\big(T^{\mathcal{P}}(z)\big)}{\Big(1-z e_{\mathcal{P}-1}\big(T^{\mathcal{P}}(z)\big)\Big)^{2}}T_{\leq k-1}^{\mathcal{P}}(z)F^\mathcal{P}(z).
      \end{eqnarray*}
      An application of Corollary \ref{when_u_equals_1_function} now gives the first result.

      An easier calculation, starting with Theorem \ref{partial_function_recursion_formula} and using the first result, gives the second result.
    \end{proof}

    An asymptotic estimate of the coefficients of $\Xi^{k,\mathcal{P}}(z)$ appears in Corollary \ref{F_Xi_coefficients}.
  \end{subsection}

  \begin{subsection}{Examples}
    \label{examples}
    We close the section with some sanity checks and explicit calculations.

    First, if $0\notin \mathcal{P}$, any rooted tree with the preimage constraint $\mathcal{P}$ must have infinite height, since the root has at least one preimage, one of those in turn has at least one preimage, and so on; this, of course, violates the requirement that a tree be finite, and so $T_{\leq h}^{k,\mathcal{P}}(u,z)=0$ for all $h,k$, just as Lemma \ref{tree_recursion_formula} implies.
    Thus, $T^{k,\mathcal{P}}(u,z)=0$ and $F^{k,\mathcal{P}}(u,z)=\big(1-z e_{\mathcal{P}-1}(0)\big)^{-1}.$
    Now we have two subcases.
    If $1\notin \mathcal{P}$, then $e_{\mathcal{P}-1}\big(0\big)=0$ and $F^{k,\mathcal{P}}(u,z)=1$; in other words, the only way a function from a finite set to itself can satisfy the condition that every element have at least 2 preimages is if the set is the empty set and the preimage condition is vacuous.
    On the other hand, if $1\in \mathcal{P}$, then $e_{\mathcal{P}-1}\big(0\big)=1$ and $F^{k,\mathcal{P}}(u,z)=(1-z)^{-1}=\sum_{n=0}^\infty z^n=\sum_{n=0}^\infty n!\frac{z^n}{n!}$; in other words, there are $n!$ functions where every element has 1 or more preimages.
    These are clearly the permutations.
    Regardless of whether or not $1\in\mathcal{P}$, note that $\Xi^{k,\mathcal{P}}(z)=0$, implying that $f^k(\set{n})=\set{n}$ for all $f\in\mathfrak{F}^{k,\mathcal{P}}$.
    This is vacuously true in the empty case and obvious in the permutation case.

    As a second sanity check, take $k=1$ and $\mathcal{P}=\mathbb{Z}_{\geq 0}$.
    In this case, Lemma \ref{tree_recursion_formula} simplifies to
    $$T^{1,\mathbb{Z}_{\geq 0}}(u,z)=z e^{T^{1,\mathbb{Z}_{\geq 0}}(u,z)} + (u-1)z;$$
    this matches Equation 24 in Section 3.1 of \cite{FlajoletOdlyzko_nolink}, where $T^{1,\mathbb{Z}_{\geq 0}}(u,z)$ is denoted by $t(u,z)$.
    Theorem \ref{function_recursion_formula}, with the observation $e_{\mathbb{Z}_{\geq 0}-1}(z)=e^z$ and with Lemma \ref{tree_recursion_formula}, simplifies to
    \begin{eqnarray}
      \label{eq:F_gen_seq}
      F^{1,\mathbb{Z}_{\geq 0}}(u,z)&=&\frac{1}{1-T^{1,\mathbb{Z}_{\geq 0}}(u,z) + (u-1)z}
    \end{eqnarray}
    and Theorem \ref{Xi_formula} simplifies to
    $$\Xi^{1,\mathbb{Z}_{\geq 0}}(z)=\frac{z T^{\mathbb{Z}_{\geq 0}}(z)}{(1-T^{\mathbb{Z}_{\geq 0}}(z))^3}.$$
    Equation \ref{eq:F_gen_seq} does not match Equation 24 in Section 3.1 of \cite{FlajoletOdlyzko_nolink}, where $F^{1,\mathbb{Z}_{\geq 0}}(u,z)$ is denoted by $\xi_3(u,z)$.
    The error seems to be neglecting to account for the fact that the roots of the trees, as cyclic points, should never be marked with $u$.
    This discrepancy propagates forward to Equation 25 of \cite{FlajoletOdlyzko_nolink}, where $\Xi^{1,\mathbb{Z}_{\geq 0}}(z)$ is denoted by $\Xi_3(z)$, but is washed away when one passes to the asymptotics; see the discussion containing Equation \ref{eq:FO_iterates} below.
  \end{subsection}
\end{section}

\begin{section}{Analytic Background}
  \label{analytic_background}

  This section summarizes some definitions and results necessary to apply singularity analysis to the combinatorial functions derived above.

  \begin{subsection}{Asymptotic Notation}
    Following Section A.2 of \cite{FlajoletSedgewick}, define $f\sim g$ to mean
    \begin{eqnarray}
      \label{eq:sim_def}
      \lim_{s\rightarrow s_0} \frac{f(s)}{g(s)}&=&1,
    \end{eqnarray}
    where both $s_0$ and the approach to $s_0$ on which $f$ and $g$ are defined are implicit from the context.
    In this document, $s_0$ will frequently be the $\rho_\mathcal{P}$ defined in Equation \ref{eq:rho_def} of Section \ref{singularity_analysis}, and the limit will be taken for real $s$ approaching $s_0$ from below.
    On occasion, $s_0$ will be infinity with $s$ positive integers.

    \begin{lemma}
      \label{sim_calculus}
      Suppose $f\sim g$ and $h\sim i$. Then
      \begin{eqnarray*}
        g &\sim& f\\
        f h &\sim& g i.
      \end{eqnarray*}
      If $g\sim h$, then
      $$f\sim h.$$
      If $\lim_{s\rightarrow s_0}\frac{h(s)}{g(s)}$ exists and is not $-1$, then
      $$f+h\sim g+h.$$
    \end{lemma}
    \begin{proof}
      Just check.
    \end{proof}
  \end{subsection}

  \begin{subsection}{Complex Analysis}

    Following Section IV.2 of \cite{FlajoletSedgewick}, call a subset $\Omega$ of $\mathbb{C}$ a \emph{region} if it is open and connected.

    A function $f:\Omega\rightarrow\mathbb{C}$, where $\Omega$ is a region, is analytic at $z_0\in\Omega$ if there is an open disc in $\Omega$ around $z_0$ in which $f(z)$ is equal to a convergent power series of the form $\sum_{n=0}^\infty c_n(z-z_0)^n$.

    Following Definition VI.1 on page 389 of \cite{FlajoletSedgewick}, for real numbers $1<R$ and $0<\phi<\frac{\pi}{2}$, define
    $$\Delta(\phi, R)=\{z\in\mathbb{C}\mid |z|<R, z\neq 1, \phi<|\arg(z-1)|\};$$
    see Figure \ref{Delta_domain}.
    A function is $\Delta$-analytic if it is analytic in a set of the form $\Delta=\Delta(\phi,R)$.
    The book \cite{PemantleWilson} uses the colorful phrase \emph{Camembert-shaped region} to describe a very similar construction.

    \begin{figure}
      \begin{center}
        \begin{tikzpicture}[scale=2];
          \draw (-1.5,0) -- (1.5,0);
          \draw (0,-1.5) -- (0,1.5);
          \draw[fill=gray!20] (.4,0) -- (30:1.2cm) arc (30:330:1.2cm) -- cycle;
          \draw (.4:.75cm) arc (0:22.:.75cm);
          \draw(.35,.1) node {1};
          \draw(1.15,.1) node {$R$};
          \draw(-1.35,.1) node {$-R$};
          \draw(.8,.2) node {$\phi$};
        \end{tikzpicture}
      \end{center}
      \caption{$\Delta(\phi,R)$}
      \label{Delta_domain}
    \end{figure}
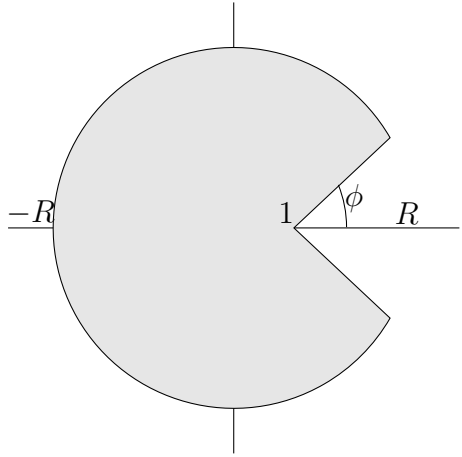

    Following conditions $H_i$ on pages 402-403 and on page 453 of \cite{FlajoletSedgewick}, a function $\epsilon(z)$ analytic at 0 with radius of convergence $R$ is said to satisfy the \emph{H-schema} if it satisfies the following three conditions.
    \begin{itemize}
    \item $\epsilon(0)\neq 0$ and $\epsilon(z)$ cannot be written as $\epsilon_0+\epsilon_1 z$ for any $\epsilon_i\in\mathbb{C}$
    \item $\coeff{z^n}\epsilon(z)\geq 0$ for all $n\geq 0$
    \item there is a unique $\tau\in\mathbb{R}_{>0}$ such that $\tau<R$ and $\epsilon(\tau)-\tau\epsilon'(\tau)=0$
    \end{itemize}

    Following Definition IV.5 on page 266 of \cite{FlajoletSedgewick}, the \emph{support} of a power series $f(z)=\sum_{n=0}^\infty f_n z^n$ is the set
    $${\rm supp}(f)=\{n\geq 0\mid f_n\neq 0\}.$$
    The \emph{period} of $f$ is the largest $d$ for which there is some $r$ such that ${\rm supp}(f)\subseteq r+d\mathbb{Z}_{\geq 0}$.
    If the period of $f$ is 1, then $f$ is said to be \emph{aperiodic}.

    The point of the previous definition is that if $f(z)$ has period $p>1$, then there is an integer $r$ and a function $g$ such that $f(z)=z^r g(z^p)$, and the $z^p$ term wraps $z$ around the origin $p$ times, invalidating the integral evaluation that underpins singularity analysis.

    All of the singularity analysis in this paper ultimately relies on the Singular Inversion Theorem \ref{SIT}, which we reprove below.
    Before doing so, we record needed results, some without proof.

    The first is the Daffodil Lemma IV.1 on page 266 of \cite{FlajoletSedgewick}.

    \begin{lemma}
      \label{daffodil}
            {\rm (Daffodil Lemma)} Let $\sigma(z)=\sum_{n\geq 0} \sigma_n z^n$ have nonnegative coefficients, have at least two elements in its support, and converge for all $|z|\leq \rho$, where $\rho\in\mathbb{R}_{>0}$.
            If $\xi\notin\mathbb{R}_{\geq 0}$ satisfies
            \begin{eqnarray*}
              |\xi| &\leq&\rho\\
              \ |\sigma(\xi)| &=& \sigma(|\xi|),
            \end{eqnarray*}
            there is a $d\in\mathbb{Z}_{\geq 1}$ and some $0<r<d$ in $\mathbb{Z}$ such that $\arg(\xi)=2\pi\frac{r}{d}$.
            In particular, if $\sigma$ is aperiodic, then
            $$|\sigma(\xi)|< \sigma(\rho)$$
            for all $\xi\notin\mathbb{R}_{\geq 0}$ satisfying $|\xi|\leq\rho$.
    \end{lemma}
    \begin{proof}
      Write $\theta=\arg(\xi)$.

      Since $\sigma$ has more than one element in its support, the triangle inequality gives $|\sum_{n\geq 0} \sigma_n \xi^n|\leq  \sum_{n\geq 0} \sigma_n |\xi|^n$ with equality iff $\sigma_n \xi^n=\sigma_n |\xi|^n$ for all $n$.
      Thus, the hypothesis $|\sigma(\xi)| = \sigma(|\xi|)$ ensures $\xi^n=|\xi|^n$ for all $n\in{\rm supp}(\sigma)$.
      In other words,
      \begin{eqnarray}
        e^{i\theta n}&=&\label{eq:real_arg}
        1
      \end{eqnarray}
      for all $n\in{\rm supp}(\sigma)$.
      Fixing one such nonzero $n$, we have $\theta n\in 2\pi\mathbb{Z}$ and $\frac{\theta}{2\pi}\in\mathbb{Q}$.
      Write $\theta=2\pi\frac{r}{d}$ with $0\leq r <d$; since $\xi\notin\mathbb{R}_{\geq 0}$, $r\neq 0$.

      For the second claim, note that Equation \ref{eq:real_arg} generalizes to all $n$ that are a $\mathbb{Z}$-linear combination of elements in ${\rm supp}(\sigma)$.
      Thus, if $\sigma$ is aperiodic, we may take $n=1$, concluding $\frac{\theta}{2\pi}$ is an integer and $\xi\in\mathbb{R}_{\geq 0}$.
      This contradiction shows that the condition $|\sigma(\xi)| = \sigma(|\xi|)$ has been violated.
      Since all of the coefficients $\sigma_n$ are nonnegative, we further conclude $|\sigma(\xi)|< \sigma(|\xi|)\leq\sigma(\rho)$.
    \end{proof}

    \begin{lemma}
      \label{analytic_inversion}
            {\rm (Analytic Inversion Lemma)} Let $\psi(z)$ be analytic at $y_0$ with $z_0=\psi(y_0)$.
            If $\psi'(y_0)\neq 0$, then there is function $\sigma(z)$ analytic in a neighborhood of $z_0$ that satisfies
            \begin{eqnarray*}
              \psi\big(\sigma(z)\big)&=&z\\
              \sigma(z_0)&=&y_0.
            \end{eqnarray*}
    \end{lemma}
    \begin{proof}
      This is Analytic Inversion Lemma IV.2 on page 275 of \cite{FlajoletSedgewick}.
    \end{proof}

    \begin{lemma}
      \label{singular_inversion}
            {\rm (Singular Inversion Lemma)} Let $\psi(z)$ be analytic at $y_0$ with $z_0=\psi(y_0)$.
            If $\psi'(y_0)= 0$ and $\psi''(y_0)\neq 0$, then there is neighborhood of $z_0$ such that, for any ray emanating from $z_0$, there are two functions $y_i(z)$ analytic in the neighborhood slit along the ray such that $y_i$ has a singularity at $z_0$.
    \end{lemma}
    \begin{proof}
      This is a slightly weaker statement of Singular Inversion Lemma IV.3 on page 277 of \cite{FlajoletSedgewick}.
    \end{proof}

    The following result, Singular Inversion Theorem VI.6 on page 404 of \cite{FlajoletSedgewick}, has been the primary goal of this section.

    \begin{theorem}
      \label{SIT}
            {\rm (Singular Inversion Theorem)} Let $\epsilon(z)$ satisfy the H-schema; in particular, there is a unique $\tau\in \mathbb{R}_{>0}$ such that
            \begin{eqnarray*}
              \epsilon(\tau)-\tau\epsilon'(\tau)&=&0\\
              \tau&<&R,
            \end{eqnarray*}
            where $R$ is the radius of convergence of $\epsilon(z)$ at 0.
            Then there is a solution $\sigma(z)$ to
            \begin{eqnarray*}
              \sigma(z)&=&z\epsilon\big(\sigma(z)\big)\\
              \sigma(0)&=&0.
            \end{eqnarray*}
            The radius of convergence of $\sigma(z)$ is
            $$\rho=\frac{\tau}{\epsilon(\tau)},$$
            and
            \begin{eqnarray*}
              \sigma(z) &\sim& \tau-\sqrt{\frac{2\epsilon(\tau)}{\epsilon''(\tau)}}\Big(1-\frac{z}{\rho}\Big)^{\frac{1}{2}}\\
              \tau-\sigma(z) &\sim& \sqrt{\frac{2\epsilon(\tau)}{\epsilon''(\tau)}}\Big(1-\frac{z}{\rho}\Big)^{\frac{1}{2}}.
            \end{eqnarray*}
            If $\epsilon$ is aperiodic, then there is some $r>\rho$ such that $\sigma(z)$ is analytic on the open disc of radius $r$ around the origin, except slit along the ray $\mathbb{R}_{\geq\rho}$.
            In particular, in this case, $\sigma(\rho z)$ is $\Delta$-analytic.
    \end{theorem}
    \begin{proof}
      First, recall that the H-schema ensures the coefficients of $\epsilon$ are nonnegative reals and $\epsilon$ is not constant or linear.
      Since $\tau>0$, this implies
      \begin{eqnarray}
        \tau\epsilon''(\tau)&>& \label{eq:eps_sq}
        0.
      \end{eqnarray}

      Now, Lagrange Inversion Theorem \ref{lagrange_inversion_theorem} ensures the formal power series
      \begin{eqnarray}
        \sigma(z)&=& \label{lit_sigma} \sum_{n= 1}^\infty\frac{1}{n}\coeff{\lambda^{n-1}}\epsilon(\lambda)^n z^n
      \end{eqnarray}
      is the unique solution to
      \begin{eqnarray}
        \sigma(z)&=& \label{sigma_def} z\epsilon\big(\sigma(z)\big)
      \end{eqnarray}
      that sends 0 to 0, and it follows from this formula and the H-schema assumption that the coefficients of $\sigma(z)$ are nonnegative.
      By Analytic Inversion Lemma \ref{analytic_inversion}, this $\sigma(z)$ is analytic at $0$.

      Let $r$ denote the radius of convergence of $\sigma(z)$, and we work to show that $L=\lim_{x\rightarrow r^-}\sigma(x)$ is equal to $\tau$.
      Because $\sigma(z)$ has nonnegative coefficients, $L$ either exists and is a nonnegative real or is infinity.
      If $L<\tau$, then the fact that $\tau<R$ implies that $\frac{\partial}{\partial z}\frac{z}{\epsilon(z)}=\frac{\epsilon(z)-z\epsilon'(z)}{\epsilon(z)^2}$ is nonzero at $L$ allows Analytic Inversion Lemma \ref{analytic_inversion} to show that $\sigma(z)$ is analytic at $r$; since $\sigma(z)$ has nonnegative coefficients, it is then analytic on the entire circle of radius $r$, a contradiction.
      If $L>\tau$, then there is a $0<s<r$ such that $\sigma(s)=\tau$, so $(\frac{\partial}{\partial \lambda}\frac{\lambda}{\epsilon(\lambda)})|_{\sigma(s)}=\frac{\epsilon(\tau)-\tau\epsilon'(\tau)}{\epsilon(\tau)^2}$ and
      \begin{eqnarray}
        \left(\frac{\partial^2}{\partial \lambda^2}\frac{\lambda}{\epsilon(\lambda)}\right)|_{\sigma(s)}&=&\label{taylor2} -\frac{\tau\epsilon''(\tau)}{\epsilon(\tau)^2};
      \end{eqnarray}
      the definition of $\tau$ then ensures that the first of these expressions is 0 and Equation \ref{eq:eps_sq} says the second is nonzero, so Singular Inversion Lemma \ref{singular_inversion} implies that $\sigma(z)$ has a singularity at $s$, which is again a contradiction.
      That is, $\lim_{x\rightarrow r^-}\sigma(x)=\tau$, as claimed.

      Now, since $\sigma(z)$ and $\frac{z}{\epsilon(z)}$ are compositional inverses of each other,
$$\rho=\frac{\tau}{\epsilon(\tau)}=\lim_{x\rightarrow r^-}\frac{\sigma(x)}{\epsilon\big(\sigma(x)\big)}=\lim_{x\rightarrow r^-}x=r,$$
      and we have confirmed that $\rho$ is the radius of convergence of $\sigma(z)$.
      In particular, another application of Singular Inversion Lemma \ref{singular_inversion}, with $\psi(z)=\frac{z}{\epsilon(z)}$ and $y_0=\tau$, says that $\sigma(z)$ is analytic in a neighborhood of $\rho$ slit along $\mathbb{R}_{\geq \rho}$.

      Next, for the asymptotic behavior of $\sigma(z)$ for $z$ near $\rho$, consider the Taylor expansion of $\frac{\tau}{\epsilon(\tau)}-\frac{z}{\epsilon(z)}$ around $z=\tau$.
      Since $\epsilon(z)$ is analytic at $\tau$, we have
      \begin{eqnarray*}
        \frac{\tau}{\epsilon(\tau)}-\frac{z}{\epsilon(z)} &=& \sum_{n=0}^\infty\frac{1}{n!}\Big(\frac{\partial^n}{\partial\lambda^n}\big(\frac{\tau}{\epsilon(\tau)}-\frac{\lambda}{\epsilon(\lambda)}\big)\Big)\Big|_{\lambda=\tau}(z-\tau)^n.
      \end{eqnarray*}
      The $n=0$ term in the Taylor expansion is trivially 0, and the definition of $\tau$ kills the $n=1$ term.
      The $n=2$ term was, in effect, calculated in Equation \ref{taylor2}.
      As $z\rightarrow\tau$, we now have
      \begin{eqnarray*}
        \frac{\tau}{\epsilon(\tau)}-\frac{z}{\epsilon(z)} &=& \frac{1}{2}\frac{\tau\epsilon''(\tau)}{\epsilon(\tau)^2}(z-\tau)^2+\sum_{n=3}^\infty\frac{1}{n!}\Big(\frac{\partial^n}{\partial\lambda^n}\big(\frac{\tau}{\epsilon(\tau)}-\frac{\lambda}{\epsilon(\lambda)}\big)\Big)\Big|_{\lambda=\tau}(z-\tau)^n\\
        &\sim& \frac{\tau\epsilon''(\tau)}{2\epsilon(\tau)^2}(z-\tau)^2\\
        &=& \frac{\rho\epsilon''(\tau)}{2\epsilon(\tau)}(z-\tau)^2.
      \end{eqnarray*}
      Since $\sigma(z)$ is the inverse of $\frac{z}{\epsilon(z)}$ in a slit neighborhood of $z=\rho$, we get the local approximation
      \begin{eqnarray*}
        \rho-z &=& \frac{\tau}{\epsilon(\tau)}-\frac{\sigma(z)}{\epsilon(\sigma(z))}\\
        &\sim& \frac{\rho\epsilon''(\tau)}{2\epsilon(\tau)}\big(\sigma(z)-\tau\big)^2\\
        \frac{2\epsilon(\tau)}{\epsilon''(\tau)}\left(1-\frac{z}{\rho}\right) &\sim& \big(\sigma(z)-\tau\big)^2.
      \end{eqnarray*}
      Note that for real $x$, $x$ approaches $\tau$ from below as $\sigma(x)$ approaches $\rho$ from below, which tells us which branch to take when solving the asymptotic behavior equation for $\sigma(z)$.
      That is,
$$\sigma(z)-\tau\sim-\sqrt{\frac{2\epsilon(\tau)}{\epsilon''(\tau)}}\Big(1-\frac{z}{\rho}\Big)^{\frac{1}{2}},$$
      which gives the second claim about the asymptotic behavior.
      The first asymptotic claim then follows immediately by applying Lemma \ref{sim_calculus}.

      Now assume $\epsilon(z)$ is aperiodic.

      We first argue that $\sigma(z)$ is aperiodic.
      Suppose there are $d,r\in\mathbb{Z}_{\geq 0}$ such that ${\rm supp}(\sigma)\subseteq r+d\mathbb{Z}_{\geq 0}$; the goal is to show that $d=1$.
      Write
      \begin{eqnarray*}
        \sigma(z)&=&z^r\sum_{k=0}^\infty\sigma_k z^{d k}\\
        \epsilon(z)&=&\sum_{n=0}^\infty\epsilon_n\frac{z^n}{n!},
      \end{eqnarray*}
      noting, as in the paragraph containing Equation \ref{lit_sigma}, that each $\sigma_k$ is nonnegative.
      Since these coefficients are nonnegative, there can be no cancellation of coefficients when expanding out the fact, immediate from Equation \ref{sigma_def}, that $\sigma(z)= z\sum_{n=0}^\infty\epsilon_n\frac{\left(z^r\sum_{k=0}^\infty\sigma_k z^{d k}\right)^n}{n!}$.
      In particular,
      \begin{eqnarray}
        1+r n +d k_1+\dots + d k_n&\in& \label{mod_supp}{\rm supp}(\sigma)\subseteq r+d\mathbb{Z}_{\geq 0}
      \end{eqnarray}
      for every $n\in{\rm supp}(\epsilon)$ and $k_1,\dots,k_n\in{\rm supp}(\sigma)$.
      Taking $n=0$ in Equation \ref{mod_supp} gives that $1\in r+d\mathbb{Z}_{\geq 0}$.
      If $d\neq 1$, we have $r=1$ and wish to derive a contradiction.
      But then Equation \ref{mod_supp} becomes $1+n +d k_1+\dots + d k_n\in 1+d\mathbb{Z}_{\geq 0}$, so $d|n$ for all $n$ in the support of $\epsilon$, and this contradicts that $\epsilon$ is aperiodic.

      We next show that $\sigma(z)$ is analytic on a region containing the punctured circle $\{z\neq\rho\mid |z|=\rho\}$.
      Let $\xi$ be a point on that punctured circle, and let $\eta=\sigma(\xi)$.
      Since $\sigma(z)$ is aperiodic, Daffodil Lemma \ref{daffodil} yields that $|\eta|=|\sigma(\xi)|<\sigma(\rho)=\tau$.
      Write $\epsilon(z)=\sum_{n=0}^\infty\epsilon_n\frac{z^n}{n!}$ and note that
      \begin{eqnarray*}
        |\epsilon(\eta)-\eta\epsilon'(\eta)| &=& |\sum_{n=0}^\infty\epsilon_n\frac{\eta^n}{n!}-\eta\sum_{n=1}^\infty\epsilon_n\frac{\eta^{n-1}}{(n-1)!} |\\
        &=& |\epsilon_0-\sum_{n=1}^\infty\epsilon_n(n-1)\frac{\eta^n}{n!}|\\
        &\geq& |\epsilon_0|-|\sum_{n=1}^\infty\epsilon_n(n-1)\frac{\eta^n}{n!}|\\
        &\geq& \epsilon_0-\sum_{n=1}^\infty\epsilon_n(n-1)\frac{|\eta|^n}{n!}\\
        &>& \epsilon_0-\sum_{n=1}^\infty\epsilon_n(n-1)\frac{\tau^n}{n!},
      \end{eqnarray*}
      where we have used that $\epsilon$ is not linear from the H-schema condition in order to ensure that the last inequality is strict.
      Continuing,
      \begin{eqnarray*}
        |\epsilon(\eta)-\eta\epsilon'(\eta)| &>& \epsilon_0-\sum_{n=1}^\infty\epsilon_n(n-1)\frac{\tau^n}{n!}\\
        &=& \epsilon(\tau)-\tau\epsilon'(\tau)\\
        &=&0.
      \end{eqnarray*}
      Thus, $(\frac{\partial}{\partial \lambda}\frac{\lambda}{\epsilon(\lambda)})|_{\eta}=\frac{\epsilon(\eta)-\eta\epsilon'(\eta)}{\epsilon(\eta)^2}\neq 0$, and Analytic Inversion Lemma \ref{analytic_inversion} again applies to ensure that $\sigma(z)$ is analytic at $z=\eta$.

      Finally, we need to extend this to show $\sigma(z)$ is analytic on a $\Delta$-domain; we do so by finding some radius larger than $\rho$ for which $\sigma(z)$ is analytic for all $z\notin \mathbb{R}_{\geq\rho}$.
      The only issue is to ensure that the radii of the open balls around each point on the boundary do not shrink as the point gets close to $\rho$.
      But recall that $\sigma(z)$ is analytic on the slit disc near $z=\rho$, so one does not need to check the radii of points on the circle arbitrarily close to $\rho$, but rather some closed subset of the punctured circle.
      Compactness yields a finite cover, and the rest is standard.
    \end{proof}

    The standard way one uses the previous result is to plug its implication into the next result.

    \begin{corollary}
      \label{sim_transfer}
            {\rm (sim-transfer Corollary)} Let $\alpha\notin \mathbb{Z}_{\leq 0}$.
            If $f(z)$ is $\Delta$-analytic and
            $$f(z)\sim(1-z)^{-\alpha}$$
            as $z$ in the $\Delta$-domain approach 1, then
            $$\coeff{z^n} f(z) \sim \frac{n^{\alpha-1}}{\Gamma(\alpha)}.$$
    \end{corollary}
    \begin{proof}
      This is sim-transfer Corollary VI.1 on page 392 of \cite{FlajoletSedgewick}.
    \end{proof}
  \end{subsection}
\end{section}

\begin{section}{Singularity Analysis}
  \label{singularity_analysis}

  The goal of this section is to apply singularity analysis to the combinatorial functions derived in Section \ref{preimage_combinatorics}.
  In particular, doing this for $\Xi^{k,\mathcal{P}}(z)$ accomplishes the primary goal of this paper, namely, to obtain asymptotic results about iterated images of a random mapping that is subject to the constraint that all preimage sizes are in $\mathcal{P}$.
  These mappings are explicitly explored in the seminal paper \cite{ArneyBender1982}, and singularity analysis is explicitly applied to them in Example VII.10 of \cite{FlajoletSedgewick}, but neither of these references address the size of iterated images.
  This was done for $\mathcal{P}=\mathbb{Z}_{\geq 0}$ in the Direct Parameters Theorem 2 of \cite{FlajoletOdlyzko_nolink}, and the treatment below follows that approach.

  \begin{subsection}{Preliminaries}

    Throughout this section, we work over the ring $R=\mathbb{C}$.
    This allows us to define the formal generating functions as in the previous sections.
    Given such a function, say $\sigma(z)$, we can examine its radius of convergence $\rho$ and identify the power series with the function defined on the open disc $\{z\in\mathbb{C}\mid |z|<\rho\}$ given by evaluating the power series.
    However, it is crucial in singularity analysis to be able to extend the function beyond this disc, else there would be no \emph{analysis} of the \emph{singularity} that lies on the boundary of the disc.
    In particular, the name $\sigma(z)$ will now refer to some analytic continuation of the power series, supplanting the original reference to the power series.
    Of course, the point is, and will remain, to understand the original generating function.

    Recall that for $\mathcal{P}\subseteq\mathbb{Z}_{\geq 0}$, the function $e_\mathcal{P}(z)$ is the truncated exponential which only includes the monomials designated by elements of $\mathcal{P}$; see Equation \ref{eq:eP_def}.
    Recall that $T^\mathcal{P}(z)$ denotes the generating function for rooted, labeled trees subject to the preimage constraint $\mathcal{P}$; see Section \ref{trees}, especially the second result in Theorem \ref{when_u_equals_1_tree}.

    We start with observations about $e_\mathcal{P}$ and $T^\mathcal{P}$ needed to apply the singularity analysis results.
    For $\mathcal{P}\subseteq\mathbb{Z}_{\geq 0}$ with $0\in\mathcal{P}$, recall that Equation \ref{eq:psi} defined
    \begin{eqnarray*}
      \psi_\mathcal{P}(z)&=&\frac{z}{e_\mathcal{P}(z)}.
    \end{eqnarray*}
    The original motivation for this additional notation was described prior to Lagrange Inversion Theorem \ref{lagrange_inversion_theorem}.
    Now that the discussion involves complex analysis, this motivation is complemented by the fact that analytic functions are locally invertible whenever their derivative is nonzero; see Analytic Inversion Lemma \ref{analytic_inversion}.

    \begin{lemma}
      \label{e_analytic}
      Let $\mathcal{P}\subseteq\mathbb{Z}_{\geq 0}$ and work over the ring $R=\mathbb{C}$.
      Then
      \begin{itemize}
      \item $e_\mathcal{P}(z)$ is analytic everywhere
      \item $\psi_\mathcal{P}(z)$ is analytic everywhere that $e_\mathcal{P}(z)\neq 0$
      \item $T^\mathcal{P}(z)$ is analytic everywhere that $e_\mathcal{P}(z)\neq 0$ and $e_\mathcal{P}(z)\neq z e_{\mathcal{P}-1}(z)$.
      \end{itemize}
    \end{lemma}
    \begin{proof}
      For any $z_0\in\mathbb{C}$, note $|e_\mathcal{P}(z+z_0)|\leq e_\mathcal{P}(|z+z_0|)\leq e_{\mathbb{Z}_{\geq0}}(|z+z_0|)=e^{|z+z_0|}$ converges for all $z\in\mathbb{C}$; in particular, $e_\mathcal{P}(z+z_0)$ is analytic at 0 and has a power series expansion $\sum_{n=0}^\infty c_n z^n$.
      Thus, $e_\mathcal{P}(z)=e_\mathcal{P}\big((z-z_0)+z_0\big)=\sum_{n=0}^\infty c_i (z-z_0)^i$ is analytic about $z_0$.

      Whenever $e_\mathcal{P}(z)\neq 0$, $e_\mathcal{P}(z)^{-1}$, and hence $\psi_\mathcal{P}(z)$, is analytic.

      Finally, Analytic Inversion Lemma \ref{analytic_inversion} gives a compositional inverse of $\psi_\mathcal{P}$ everywhere that $\psi_\mathcal{P}'(z)\neq 0$.
      Note $\psi_\mathcal{P}'(z)=\frac{e_\mathcal{P}(z)- z e_{\mathcal{P}-1}(z)}{e_\mathcal{P}(z)^2}$.
      Finally, by the uniqueness of compositional inverses, this analytic compositional inverse we just found must be $T^\mathcal{P}(z)$.
    \end{proof}

    \begin{lemma}
      \label{T_H_schema}
      Let $\mathcal{P}\subseteq\mathbb{Z}_{\geq0}$, with $0\in\mathcal{P}$ and $\mathcal{P}-2\neq\emptyset$, and work over the ring $R=\mathbb{C}$.
      Then $e_\mathcal{P}$ satisfies the H-schema.
    \end{lemma}
    \begin{proof}
      By Lemma \ref{e_analytic}, $e_\mathcal{P}$ is analytic everywhere.
      Since $0\in\mathcal{P}$, $e_\mathcal{P}(0)=1\neq 0$, and, since $\mathcal{P}-2\neq\emptyset$, $e_\mathcal{P}(u)=\sum_{n\in\mathcal{P}}\frac{u^n}{n!}$ cannot be written as $\epsilon_0+\epsilon_1 u$ for any $\epsilon_i\in\mathbb{C}$.
      Clearly, all of the coefficients are nonnegative.

      Let $f(t)=t e_\mathcal{P}'(t)-\left(e_\mathcal{P}(t)-1\right)$.
      Since $e_\mathcal{P}$ is analytic, so is $f$, and, in particular, $f|_{\mathbb{R}_{\geq 0}}$ is continuous and increasing.
      Since $0\in\mathcal{P}$, $f(0)=0$.
      Since $\mathcal{P}-2\neq\emptyset$, there is some $k\in\mathcal{P}\setminus\{0,1\}$.
      Since
      \begin{eqnarray*}
        f(t) &=& \sum_{n\in\mathcal{P}\setminus\{0\}}\frac{t^n}{(n-1)!}(1-\frac{1}{n})\\
        &\geq& \frac{t^k}{(k-1)!}(1-\frac{1}{k}),
      \end{eqnarray*}
      where this last expression goes to infinity as $t$ does, there is some $0<\tau\leq (\frac{(k-1)!}{1-\frac{1}{k}})^{\frac{1}{k}}$ such that $f(\tau)=1$.
      But $f$ is increasing, so this is the only such $\tau\geq 0$.
      By the definition of $f$, $f(\tau)=1$ iff
      $$e_\mathcal{P}(\tau)-\tau e_\mathcal{P}'(\tau)=0.$$
    \end{proof}

    Given $\mathcal{P}\subseteq\mathbb{Z}_{\geq 0}$ such that $0\in\mathcal{P}$ and $\mathcal{P}-2\neq\emptyset$, define $\tau_\mathcal{P}$ to be the unique $\tau\in\mathbb{R}_{> 0}$ such that
    \begin{eqnarray*}
      e_\mathcal{P}(\tau)-\tau e_\mathcal{P}'(\tau)&=&0.
    \end{eqnarray*}
    In this case, also define
    \begin{eqnarray}
      \label{eq:rho_def}
      \rho_\mathcal{P}&=&\frac{\tau_\mathcal{P}}{e_\mathcal{P}(\tau_\mathcal{P})}.
    \end{eqnarray}
    It is frequently useful to utilize the fact, immediate from the definitions, that
    \begin{eqnarray}
      \rho_\mathcal{P}&=&\label{eq:alternate_rho}\frac{1}{e_{\mathcal{P}-1}(\tau_\mathcal{P})}.
    \end{eqnarray}

    The expression $T_{\leq k-1}^\mathcal{P}(\rho_\mathcal{P})$ will make repeated appearances in exploring the expected value of $|f^k(\set{n})|$, so define
    \begin{eqnarray}
      \tau_k^\mathcal{P}&=&\label{tau_k_def}
      T_{\leq k-1}^\mathcal{P}(\rho_\mathcal{P}).
    \end{eqnarray}
    It turns out that $(\tau_k^\mathcal{P})$ is an increasing sequence that converges to $\tau_\mathcal{P}$.

    \begin{lemma}
      \label{tau_k_converge}
      Let $\mathcal{P}\subseteq\mathbb{Z}_{\geq0}$, with $0\in\mathcal{P}$ and $\mathcal{P}-2\neq\emptyset$, and work over the ring $R=\mathbb{C}$.
      Then, for all $k\geq 0$,
      $$\tau_k^\mathcal{P}<\tau_{k+1}^\mathcal{P}<\tau_\mathcal{P},$$
      and the sequence $\tau_k^\mathcal{P}$ converges to $\tau_\mathcal{P}$.
    \end{lemma}
    \begin{proof}
      For the first claim, we induct.
      Note that
      $$\tau_0^\mathcal{P}=0<\rho_\mathcal{P}=\tau_1^\mathcal{P}=\rho_\mathcal{P}=\frac{\tau_\mathcal{P}}{e_\mathcal{P}(\tau_\mathcal{P})}<\tau_\mathcal{P},$$
      giving the basis step.
      If $\tau_k^\mathcal{P}<\tau_{k+1}^\mathcal{P}<\tau_\mathcal{P}$, then
      $$\tau_{k+1}^\mathcal{P}=\rho_\mathcal{P}e_\mathcal{P}(\tau_{k}^\mathcal{P})<\rho_\mathcal{P}e_\mathcal{P}(\tau_{k+1}^\mathcal{P})=\tau_{k+2}^\mathcal{P}=\rho_\mathcal{P}e_\mathcal{P}(\tau_{k+1}^\mathcal{P})=\frac{\tau_\mathcal{P}}{e_\mathcal{P}(\tau_\mathcal{P})}e_\mathcal{P}(\tau_{k+1}^\mathcal{P})<\tau_\mathcal{P}.$$

      For the second claim, we have a monotone function bounded above by $\tau_\mathcal{P}$, so it must converge to some $\mu\leq \tau_\mathcal{P}$.
      Taking the limit as $k$ goes to infinity of the defining relation $\tau_{k+1}^\mathcal{P}=\rho_\mathcal{P}e_\mathcal{P}(\tau_{k}^\mathcal{P})$ says that
      \begin{eqnarray*}
        \mu&=&\rho_\mathcal{P}e_\mathcal{P}(\mu).
      \end{eqnarray*}
      Recalling $\psi_\mathcal{P}$ defined in Equation \ref{eq:psi}, we have
      $$\psi_\mathcal{P}(\mu)=\frac{\mu}{e_\mathcal{P}(\mu)}=\rho_\mathcal{P}=\frac{\tau_\mathcal{P}}{e_\mathcal{P}(\tau_\mathcal{P})}=\psi_\mathcal{P}(\tau_\mathcal{P}).$$
      Since $e_\mathcal{P}(z)- z e_{\mathcal{P}-1}(z)$ is a decreasing function that evaluates to 1 when $z=0$ and to 0 when $z=\tau_\mathcal{P}$, $\psi_\mathcal{P}'(z)=\frac{e_\mathcal{P}(z)- z e_{\mathcal{P}-1}(z)}{e_\mathcal{P}(z)^2}$ is strictly positive on the interval $(0,\tau_\mathcal{P})$ and the restriction of $\psi_\mathcal{P}$ to $(0,\tau_\mathcal{P})$ must be injective.
      It follows that $\mu=\tau_\mathcal{P}$.
    \end{proof}

    \begin{lemma}
      Let $\mathcal{P}_1\subseteq\mathcal{P}_2\subseteq \mathbb{Z}_{\geq0}$, with $0\in\mathcal{P}_1$ and $\mathcal{P}_1-2\neq\emptyset$, and work over the ring $R=\mathbb{C}$.
      Then
      $$\tau_{\mathcal{P}_1}\geq\tau_{\mathcal{P}_2}.$$
    \end{lemma}
    \begin{proof}
      It follows from $\mathcal{P}_1\subseteq\mathcal{P}_2$ that $0\in\mathcal{P}_2$ and $\mathcal{P}_2-2\neq\emptyset$, so $\tau_{\mathcal{P}_2}$ is actually defined.

      It is immediate from the definition of $\tau_{\mathcal{P}_i}$ that
      $$1=\sum_{n\in\mathcal{P}_i\setminus\{0\}}\frac{n-1}{n!}\tau_{\mathcal{P}_i}^n.$$
      In particular,
      $$\sum_{n\in\mathcal{P}_1}\frac{n-1}{n!}\tau_{\mathcal{P}_1}^n=\sum_{n\in\mathcal{P}_2}\frac{n-1}{n!}\tau_{\mathcal{P}_2}^n.$$

      Suppose, contrary to the claim, that $\tau_{\mathcal{P}_1}<\tau_{\mathcal{P}_2}$.
      Then, since $\emptyset\subset\mathcal{P}_1\subseteq\mathcal{P}_2$,
      \begin{eqnarray*}
        \sum_{n\in\mathcal{P}_1}\frac{n-1}{n!}\tau_{\mathcal{P}_2}^n &\leq& \sum_{n\in\mathcal{P}_2}\frac{n-1}{n!}\tau_{\mathcal{P}_2}^n\\
        &=& \sum_{n\in\mathcal{P}_1}\frac{n-1}{n!}\tau_{\mathcal{P}_1}^n\\
        &<& \sum_{n\in\mathcal{P}_1}\frac{n-1}{n!}\tau_{\mathcal{P}_2}^n,
      \end{eqnarray*}
      a contradiction.
    \end{proof}
  \end{subsection}

  \begin{subsection}{Results}

    We now begin the singularity analysis in earnest.
    For the next result, recall the definition of $\sim$ in Equation \ref{eq:sim_def}.
    The implicit limits are for $z\rightarrow \rho_\mathcal{P}$.

    The next result (as well as the first equation in Corollary \ref{F_Xi_coefficients} below) is basically an application of Theorem VII.2 on page 453 of \cite{FlajoletSedgewick}.

    \begin{corollary}
      \label{T_asymptotic}
      Let $\mathcal{P}\subseteq\mathbb{Z}_{\geq0}$, with $0\in\mathcal{P}$ and $\mathcal{P}-2\neq\emptyset$, and work over the ring $R=\mathbb{C}$.
      Then
      \begin{eqnarray*}
        T^\mathcal{P}(z) &\sim& \tau_\mathcal{P}-\sqrt{\frac{2e_\mathcal{P}(\tau_\mathcal{P})}{e_{\mathcal{P}-2}(\tau_\mathcal{P})}}\Big(1-\frac{z}{\rho_\mathcal{P}}\Big)^{\frac{1}{2}}\\
        \tau_\mathcal{P}-T^\mathcal{P}(z) &\sim& \sqrt{\frac{2e_\mathcal{P}(\tau_\mathcal{P})}{e_{\mathcal{P}-2}(\tau_\mathcal{P})}}\Big(1-\frac{z}{\rho_\mathcal{P}}\Big)^{\frac{1}{2}}
      \end{eqnarray*}
      If, in addition, $\gcd(\mathcal{P})=1$, then there is a $r>\rho_\mathcal{P}$ such that $T^\mathcal{P}(z)$ is analytic on the open disc of radius $r$, slit along $\mathbb{R}_{\geq\rho_\mathcal{P}}$; in other words, under this assumption, $T^\mathcal{P}(\rho_\mathcal{P} z)$ is $\Delta-analytic$.
    \end{corollary}
    \begin{proof}
      By Lemma \ref{T_H_schema}, we may apply Singular Inversion Theorem \ref{SIT}, and the result is immediate from the observation that $\gcd(\mathcal{P})=1$ iff $e_\mathcal{P}(z)$ is aperiodic.
    \end{proof}

    The following computation is not necessarily interesting in its own right, but it arises a few times.

    \begin{corollary}
      \label{e_of_T}
      Let $\mathcal{P},\mathcal{Q}\subseteq\mathbb{Z}_{\geq0}$ with $0\in\mathcal{P}$ and $\mathcal{P}-2\neq\emptyset$, and work over the ring $R=\mathbb{C}$. Then there is a slit open disc around $\rho_\mathcal{P}$ such that
      $$e_\mathcal{Q}\big(T^\mathcal{P}(z)\big) \sim e_\mathcal{Q}(\tau_\mathcal{P})-\sqrt{\frac{2e_\mathcal{P}(\tau_\mathcal{P})}{e_{\mathcal{P}-2}(\tau_\mathcal{P})}}\Big(1-\frac{z}{\rho_\mathcal{P}}\Big)^{\frac{1}{2}}e_{\mathcal{Q}-1}(\tau_\mathcal{P})$$
      and
      $$e_\mathcal{Q}\big(T^\mathcal{P}(z)\big) - e_\mathcal{Q}(\tau_\mathcal{P}) \sim -\sqrt{\frac{2e_\mathcal{P}(\tau_\mathcal{P})}{e_{\mathcal{P}-2}(\tau_\mathcal{P})}}\Big(1-\frac{z}{\rho_\mathcal{P}}\Big)^{\frac{1}{2}}e_{\mathcal{Q}-1}(\tau_\mathcal{P}).$$
      If, in addition, $\gcd(\mathcal{P})=1$, then there is a $r>\rho_\mathcal{P}$ such that $e_\mathcal{Q}(T^\mathcal{P})(z)$ is analytic on the open disc of radius $r$, slit along $\mathbb{R}_{\geq\rho_\mathcal{P}}$; in other words, under this assumption, $e_\mathcal{Q}\big(T^\mathcal{P}(\rho_\mathcal{P} z)\big)$ is $\Delta-analytic$.
    \end{corollary}
    \begin{proof}
      The first claim is immediate from the definition of $\sim$ and the fact that $T^\mathcal{P}(\rho_\mathcal{P})=\tau_\mathcal{P}$; alternately, it follows from Lemma \ref{sim_calculus}, once the second claim has been shown.

      An application of l'Hopital and the fact that
      $$T^\mathcal{P}(z)'=\frac{e_\mathcal{P}\big(T^\mathcal{P}(z)\big)}{1-z e_{\mathcal{P}-1}\big(T^\mathcal{P}(z)\big)},$$
      give that
      \begin{eqnarray*}
        && \lim_{z\rightarrow \rho_\mathcal{P}^-} \frac{e_\mathcal{Q}\left(T^\mathcal{P}\left(z\right)\right)-e_\mathcal{Q}\left(\tau_\mathcal{P}\right)}{1-z e_{\mathcal{P}-1}\big(T^\mathcal{P}\left(z\right)\big)}\\
        &=& \lim_{z\rightarrow \rho_\mathcal{P}^-} \frac{e_{\mathcal{Q}-1}\left(T^\mathcal{P}\left(z\right)\right)T^\mathcal{P}(z)'}{-e_{\mathcal{P}-1}\big(T^\mathcal{P}(z)\big)-z e_{\mathcal{P}-2}\big(T^\mathcal{P}(z)\big)T^\mathcal{P}(z)'}\\
        &=& \lim_{z\rightarrow \rho_\mathcal{P}^-} \frac{e_{\mathcal{Q}-1}\left(T^\mathcal{P}\left(z\right)\right)e_\mathcal{P}\left(T^\mathcal{P}\left(z\right)\right)}{-e_{\mathcal{P}-1}\big(T^\mathcal{P}(z)\big)\Big(1-z e_{\mathcal{P}-1}\big(T^\mathcal{P}\left(z\right)\big)\Big)-z e_{\mathcal{P}-2}\big(T^\mathcal{P}(z)\big)e_\mathcal{P}\big(T^\mathcal{P}(z)\big)}.
      \end{eqnarray*}
      Since $T^\mathcal{P}(\rho_\mathcal{P})=\tau_\mathcal{P}$, Equation \ref{eq:alternate_rho} simplifies this to
      \begin{eqnarray}
        \lim_{z\rightarrow \rho_\mathcal{P}^-} \frac{e_\mathcal{Q}\left(T^\mathcal{P}\left(z\right)\right)-e_\mathcal{Q}\left(\tau_\mathcal{P}\right)}{1-z e_{\mathcal{P}-1}\big(T^\mathcal{P}\left(z\right)\big)} &=& \nonumber \frac{e_{\mathcal{Q}-1}(\tau_\mathcal{P})e_\mathcal{P}\left(\tau_\mathcal{P}\right)}{-e_{\mathcal{P}-1}\left(\tau_\mathcal{P}\right)\big(1-\rho_\mathcal{P} e_{\mathcal{P}-1}\left(\tau_\mathcal{P}\right)\big)-\rho_\mathcal{P} e_{\mathcal{P}-2}(\tau_\mathcal{P}) e_\mathcal{P}\left(\tau_\mathcal{P}\right)}\\
        &=& \nonumber \frac{e_{\mathcal{Q}-1}\left(\tau_\mathcal{P}\right)e_\mathcal{P}\left(\tau_\mathcal{P}\right)}{0-\rho_\mathcal{P} e_{\mathcal{P}-2}(\tau_\mathcal{P}) e_\mathcal{P}\left(\tau_\mathcal{P}\right)}\\
        &=& \label{eq:temp_lim} - \frac{e_{\mathcal{Q}-1}\left(\tau_\mathcal{P}\right)}{\rho_\mathcal{P} e_{\mathcal{P}-2}(\tau_\mathcal{P}) }.
      \end{eqnarray}

      Next, apply l'Hopital again to see that
      \begin{eqnarray*}
        && \lim_{z\rightarrow \rho_\mathcal{P}^-} \frac{\left(e_\mathcal{Q}\left(T^\mathcal{P}\left(z\right)\right)-e_\mathcal{Q}\left(\tau_\mathcal{P}\right)\right)^2}{\frac{2e_\mathcal{P}(\tau_\mathcal{P})}{e_{\mathcal{P}-2}(\tau_\mathcal{P})}\left(1-\frac{z}{\rho_\mathcal{P}}\right)e_{\mathcal{Q}-1}\left(\tau_\mathcal{P}\right)^2}\\
        &=& \lim_{z\rightarrow \rho_\mathcal{P}^-} \frac{2\left(e_\mathcal{Q}\left(T^\mathcal{P}\left(z\right)\right)-e_\mathcal{Q}\left(\tau_\mathcal{P}\right)\right)e_{\mathcal{Q}-1}\left(T^\mathcal{P}(z)\right)T(z)'}{\frac{2e_\mathcal{P}(\tau_\mathcal{P})}{e_{\mathcal{P}-2}(\tau_\mathcal{P})}\left(-\frac{1}{\rho_\mathcal{P}}\right)e_{\mathcal{Q}-1}\left(\tau_\mathcal{P}\right)^2}\\
        &=& -\frac{\rho_\mathcal{P}e_{\mathcal{P}-2}(\tau_\mathcal{P})}{e_{\mathcal{P}}(\tau_\mathcal{P})}\lim_{z\rightarrow \rho_\mathcal{P}^-} \frac{\left(e_\mathcal{Q}\left(T^\mathcal{P}\left(z\right)\right)-e_\mathcal{Q}\left(\tau_\mathcal{P}\right)\right)}{1-z e_{\mathcal{P}-1}\big(T^\mathcal{P}(z)\big)}\frac{e_{\mathcal{Q}-1}\left(T^\mathcal{P}(z)\right)e_\mathcal{P}\big(T^\mathcal{P}(z)\big)}{e_{\mathcal{Q}-1}\left(\tau_\mathcal{P}\right)^2}.
      \end{eqnarray*}
      By Equation \ref{eq:temp_lim}, this is
      \begin{eqnarray*}
        \lim_{z\rightarrow \rho_\mathcal{P}^-} \frac{\left(e_\mathcal{Q}\left(T^\mathcal{P}\left(z\right)\right)-e_\mathcal{Q}\left(\tau_\mathcal{P}\right)\right)^2}{\frac{2e_\mathcal{P}(\tau_\mathcal{P})}{e_{\mathcal{P}-2}(\tau_\mathcal{P})}\left(1-\frac{z}{\rho_\mathcal{P}}\right)e_{\mathcal{Q}-1}\left(\tau_\mathcal{P}\right)^2} &=& -\frac{\rho_\mathcal{P}e_{\mathcal{P}-2}(\tau_\mathcal{P})}{e_{\mathcal{P}}(\tau_\mathcal{P})} \left(- \frac{e_{\mathcal{Q}-1}\left(\tau_\mathcal{P}\right)}{\rho_\mathcal{P} e_{\mathcal{P}-2}(\tau_\mathcal{P}) } \right)\frac{e_{\mathcal{Q}-1}\left(\tau_\mathcal{P}\right)e_\mathcal{P}\big(\tau_\mathcal{P}\big)}{e_{\mathcal{Q}-1}\left(\tau_\mathcal{P}\right)^2}\\
        &=& 1.
      \end{eqnarray*}

      In other words, we have shown from the definition of $\sim$ that, possibly up to sign,
      $$e_\mathcal{Q}\big(T^\mathcal{P}(z)\big) - e_\mathcal{Q}(\tau_\mathcal{P}) \sim \sqrt{\frac{2e_\mathcal{P}(\tau_\mathcal{P})}{e_{\mathcal{P}-2}(\tau_\mathcal{P})}}\Big(1-\frac{z}{\rho_\mathcal{P}}\Big)^{\frac{1}{2}}e_{\mathcal{Q}-1}(\tau_\mathcal{P}).$$
      It is clear from Corollary \ref{T_asymptotic} that $e_\mathcal{Q}\left(T^\mathcal{P}(z)\right)$ is an increasing function as $z$ approaches $\rho_\mathcal{P}$ from below, and so the negative root is the correct one, giving the second claim.

      The final claim is immediate from the last part of Corollary \ref{T_asymptotic} and the fact that $e_\mathcal{Q}(z)$ is analytic everywhere.
    \end{proof}

    For the next result, recall that $F^\mathcal{P}(z)$ denotes the generating function for functions subject to the preimage constraint $\mathcal{P}$; see Section \ref{functions}, especially the first result in Corollary \ref{when_u_equals_1_function}.
    Also, recall that $\Xi^{\CYC, \mathcal{P}}(z)$ denotes the generating function for all cyclic points across all functions that satisfy the constraint $\mathcal{P}$; see Section \ref{cycle_points}, especially Theorem \ref{cyclic_Xi_formula}.
    Finally, recall that $\Xi^{k,\mathcal{P}}(z)$ denotes the generating function of the number of points not in the $k$th image of $f$ as one sums across all $f$ subject to the constraint $\mathcal{P}$; see Section \ref{image_combinatorics}, especially Theorem \ref{Xi_formula}.

    \begin{corollary}
      \label{F_Xi_asymptotic}
      Let $k\geq 0$ and $\mathcal{P}\subseteq\mathbb{Z}_{\geq0}$, with $0\in\mathcal{P}$ and $\mathcal{P}-2\neq\emptyset$, and work over the ring $R=\mathbb{C}$.
      Then
      \begin{eqnarray*}
        F^\mathcal{P}(z) &\sim& \frac{1}{\sqrt{2 \tau_\mathcal{P} \rho_\mathcal{P} e_{\mathcal{P}-2}(\tau_\mathcal{P})}}\Big(1-\frac{z}{\rho_\mathcal{P}}\Big)^{-\frac{1}{2}}\\
        \Xi^{\CYC,\mathcal{P}}(z) &\sim& \frac{1}{2\tau_\mathcal{P}\rho_\mathcal{P} e_{\mathcal{P}-2}(\tau_\mathcal{P})} \Big(1-\frac{z}{\rho_\mathcal{P}}\Big)^{-1}\\
        \Xi^{k,\mathcal{P}}(z) &\sim& \frac{\tau_k^{\mathcal{P}}}{\tau_\mathcal{P}^2}\sqrt{\frac{e_\mathcal{P}(\tau_\mathcal{P})}{8 e_{\mathcal{P}-2}(\tau_\mathcal{P})}} \Big(1-\frac{z}{\rho_\mathcal{P}}\Big)^{-\frac{3}{2}}\\
        P^\mathcal{P}(z) &\sim& e^{\tau_\mathcal{P}}F^\mathcal{P}(z)\\
        \Xi^{{\rm partial},k,\mathcal{P}}(z) &\sim& e^{\tau_\mathcal{P}}\Xi^{k,\mathcal{P}}(z).
      \end{eqnarray*}
      If, in addition, $\gcd(\mathcal{P})=1$, then there is a $r>\rho_\mathcal{P}$ such that $F^\mathcal{P}(z)$, $\Xi^{\CYC, \mathcal{P}}(z)$, $\Xi^{k,\mathcal{P}}(z)$, and $\Xi^{{\rm partial},k,\mathcal{P}}(z)$ are analytic on the open disc of radius $r$, slit along $\mathbb{R}_{\geq\rho_\mathcal{P}}$; in other words, under this assumption, $F^\mathcal{P}(\rho_\mathcal{P}z)$, $\Xi^{\CYC, \mathcal{P}}(\rho_\mathcal{P}z)$, $\Xi^{k,\mathcal{P}}(\rho_\mathcal{P}z)$, and $\Xi^{{\rm partial},k,\mathcal{P}}(\rho_\mathcal{P} z)$ are $\Delta-analytic$.
    \end{corollary}
    \begin{proof}
      Corollary \ref{when_u_equals_1_function} says that
      \begin{eqnarray}
        \label{eq:FP}
        F^{\mathcal{P}}(z) &=& \Big(1-z e_{\mathcal{P}-1}\big(T^{\mathcal{P}}(z)\big)\Big)^{-1}.
      \end{eqnarray}
      Corollary \ref{T_asymptotic} says, in part, that $T^\mathcal{P}$ is analytic on the open disc $\{z \mid |z|<\rho_\mathcal{P}\}$.
      In particular, for $z$ in this region and as in Equation \ref{eq:der_T},
      $$T^\mathcal{P}(z)'=\frac{e_\mathcal{P}\big(T^\mathcal{P}(z)\big)}{1-z e_{\mathcal{P}-1}\big(T^\mathcal{P}(z)\big)}$$
      is analytic, so $1-z e_{\mathcal{P}-1}\left(T^\mathcal{P}(z)\right)$ is nonzero on the region, and $F^\mathcal{P}(z)$ is analytic on the same open disc.
      Then Theorems \ref{cyclic_Xi_formula} and \ref{Xi_formula} ensure that
      \begin{eqnarray}
        \Xi^{\CYC,\mathcal{P}}(z) &=& \label{eq:cyclicXiP} z e_{\mathcal{P}-1}\big(T^{\mathcal{P}}(z)\big) F^\mathcal{P}(z)^2\\
        \Xi^{k,\mathcal{P}}(z) &=& \label{eq:XiP} z T_{\leq k-1}^{\mathcal{P}}(z) e_{\mathcal{P}-2}\big(T^{\mathcal{P}}(z)\big) F^\mathcal{P}(z)^3\\
        \Xi^{{\rm partial}, k,\mathcal{P}}\left(z\right) &=& \nonumber T_{\leq k-1}^{\mathcal{P}}\left(z\right) F^{\mathcal{P}}\left(z\right)^2 e^{T^{\mathcal{P}}\left(z\right)} \Big(z e_{\mathcal{P}-2}\left(T^{\mathcal{P}}\left(z\right)\right)F^{\mathcal{P}}\left(z\right)+1\Big).
      \end{eqnarray}
      are all analytic on the open disc.
      The final observation that $\gcd(\mathcal{P})=1$ implies a $\Delta$-analytic condition is now immediate from Corollary \ref{T_asymptotic}.

      It is easy to check that
      $$ \lim_{z\rightarrow \rho_\mathcal{P}^-} \frac{1-z e_{\mathcal{P}-1}\left(T(z)\right)}{1-\rho_\mathcal{P} e_{\mathcal{P}-1}\left(T(z)\right)} = 1,$$
      so Equation \ref{eq:FP} yields that
      \begin{eqnarray*}
        F^\mathcal{P}(z) &\sim& \Big(1-\rho_\mathcal{P} e_{\mathcal{P}-1}\big(T^{\mathcal{P}}(z)\big)\Big)^{-1}\\
        &=& \rho_\mathcal{P}^{-1}\Big(e_{\mathcal{P}-1}\big(\tau_\mathcal{P}\big)- e_{\mathcal{P}-1}\big(T^{\mathcal{P}}(z)\big)\Big)^{-1}.
      \end{eqnarray*}
      Taking Corollary \ref{e_of_T} with $\mathcal{Q}=\mathcal{P}-1$ expresses this as
      \begin{eqnarray*}
        F^{\mathcal{P}}(z) &\sim& \rho_\mathcal{P}^{-1}\sqrt{\frac{e_{\mathcal{P}-2}(\tau_\mathcal{P})}{2e_\mathcal{P}(\tau_\mathcal{P})}}\Big(1-\frac{z}{\rho_\mathcal{P}}\Big)^{-\frac{1}{2}}e_{\mathcal{P}-2}(\tau_\mathcal{P})^{-1}
      \end{eqnarray*}
      Recalling that the definition of $\tau_\mathcal{P}$ gives $e_\mathcal{P}(\tau_\mathcal{P})=\tau_\mathcal{P}e_{\mathcal{P}-1}(\tau_\mathcal{P})$ and that $\rho_\mathcal{P}=\frac{1}{e_{\mathcal{P}-1}(\tau_\mathcal{P})}$ gives the $F^\mathcal{P}$ result.

      Now, plug the $F^\mathcal{P}$ result into Equation \ref{eq:cyclicXiP} and apply Corollary \ref{e_of_T} with $\mathcal{Q}=\mathcal{P}-1$ to see
      \begin{eqnarray*}
        \Xi^{\CYC,\mathcal{P}}(z) &\sim& \rho_\mathcal{P} \Big(e_{\mathcal{P}-1}(\tau_\mathcal{P})-\sqrt{\frac{2e_\mathcal{P}(\tau_\mathcal{P})}{e_{\mathcal{P}-2}(\tau_\mathcal{P})}}\Big(1-\frac{z}{\rho_\mathcal{P}}\Big)^{\frac{1}{2}}e_{\mathcal{P}-2}(\tau_\mathcal{P})\Big) \frac{\Big(1-\frac{z}{\rho_\mathcal{P}}\Big)^{-1}}{\big(2 \tau_\mathcal{P} \rho_\mathcal{P} e_{\mathcal{P}-2}(\tau_\mathcal{P})\big)} \\
        &\sim& \rho_\mathcal{P} e_{\mathcal{P}-1}(\tau_\mathcal{P}) \frac{\Big(1-\frac{z}{\rho_\mathcal{P}}\Big)^{-1}}{2 \tau_\mathcal{P} \rho_\mathcal{P} e_{\mathcal{P}-2}(\tau_\mathcal{P})} \\
        &=& \frac{1}{2 \tau_\mathcal{P} \rho_\mathcal{P} e_{\mathcal{P}-2}(\tau_\mathcal{P})} \Big(1-\frac{z}{\rho_\mathcal{P}}\Big)^{-1}.
      \end{eqnarray*}

      Next, plug the $F^\mathcal{P}$ estimate into Equation \ref{eq:XiP} and apply Corollary \ref{e_of_T} with $\mathcal{Q}=\mathcal{P}-2$ to see
      \begin{eqnarray*}
        \Xi^{k,\mathcal{P}}(z) &\sim& \rho_\mathcal{P} T_{\leq k-1}^{\mathcal{P}}(\rho_\mathcal{P}) \frac{\Big(1-\frac{z}{\rho_\mathcal{P}}\Big)^{-\frac{3}{2}}}{\big(2 \tau_\mathcal{P} \rho_\mathcal{P} e_{\mathcal{P}-2}(\tau_\mathcal{P})\big)^\frac{3}{2}}\\
        && \cdot \Big(e_{\mathcal{P}-2}(\tau_\mathcal{P})-\sqrt{\frac{2e_\mathcal{P}(\tau_\mathcal{P})}{e_{\mathcal{P}-2}(\tau_\mathcal{P})}}\Big(1-\frac{z}{\rho_\mathcal{P}}\Big)^{\frac{1}{2}}e_{\mathcal{P}-3}(\tau_\mathcal{P})\Big)\\
        &\sim& \rho_\mathcal{P} T_{\leq k-1}^{\mathcal{P}}(\rho_\mathcal{P}) \frac{\Big(1-\frac{z}{\rho_\mathcal{P}}\Big)^{-\frac{3}{2}}}{\big(2 \tau_\mathcal{P} \rho_\mathcal{P} e_{\mathcal{P}-2}(\tau_\mathcal{P})\big)^\frac{3}{2}} e_{\mathcal{P}-2}(\tau_\mathcal{P})\\
        &=& \frac{T_{\leq k-1}^{\mathcal{P}}(\rho_\mathcal{P})}{\sqrt{8 \tau_\mathcal{P}^3 \rho_\mathcal{P} e_{\mathcal{P}-2}(\tau_\mathcal{P})}} \Big(1-\frac{z}{\rho_\mathcal{P}}\Big)^{-\frac{3}{2}}.
      \end{eqnarray*}
      Plugging in the definitions of $\rho_\mathcal{P}$ and $\tau_k^\mathcal{P}$ yields
      \begin{eqnarray*}
        \Xi^{k,\mathcal{P}}(z) &\sim& \frac{\tau_k^{\mathcal{P}}}{\sqrt{8 \tau_\mathcal{P}^3 \frac{\tau_\mathcal{P}}{e_\mathcal{P}(\tau_\mathcal{P})} e_{\mathcal{P}-2}(\tau_\mathcal{P})}} \Big(1-\frac{z}{\rho_\mathcal{P}}\Big)^{-\frac{3}{2}}\\
        &=& \frac{\tau_k^{\mathcal{P}}}{\tau_\mathcal{P}^2}\sqrt{\frac{e_\mathcal{P}(\tau_\mathcal{P})}{8 e_{\mathcal{P}-2}(\tau_\mathcal{P})}} \Big(1-\frac{z}{\rho_\mathcal{P}}\Big)^{-\frac{3}{2}}.
      \end{eqnarray*}

      Corollary \ref{e_of_T} gives that $e^{T^\mathcal{P}(z)}\sim e^{\tau_\mathcal{P}}$, so Theorem \ref{partial_function_recursion_formula} and Lemma \ref{sim_calculus} give that $P^\mathcal{P}(z) \sim e^{\tau_\mathcal{P}}F^\mathcal{P}(z)$.
      
      Finally, Corollary \ref{e_of_T} and the $F^\mathcal{P}$ result give
      \begin{eqnarray*}
        z e_{\mathcal{P}-2}\left(T^\mathcal{P}(z)\right)F^\mathcal{P}(z) &\sim& \rho_\mathcal{P} e_{\mathcal{P}-2}(\tau_\mathcal{P}) \frac{1}{\sqrt{2 \tau_\mathcal{P} \rho_\mathcal{P} e_{\mathcal{P}-2}(\tau_\mathcal{P})}}\Big(1-\frac{z}{\rho_\mathcal{P}}\Big)^{-\frac{1}{2}}\\
        &=& \sqrt{\frac{\rho_\mathcal{P} e_{\mathcal{P}-2}(\tau_\mathcal{P})}{2 \tau_\mathcal{P}}}(1-\frac{z}{\rho_\mathcal{P}})^{-\frac{1}{2}}.
      \end{eqnarray*}
      Since
      $$\lim_{z\rightarrow \rho_\mathcal{P}} \frac{1}{\sqrt{\frac{\rho_\mathcal{P} e_{\mathcal{P}-2}(\tau_\mathcal{P})}{2 \tau_\mathcal{P}}}(1-\frac{z}{\rho_\mathcal{P}})^{-\frac{1}{2}}}=0,$$
      Theorem \ref{Xi_formula}, Lemma \ref{sim_calculus}, and the definition of $\tau_k^\mathcal{P}$ give that
      \begin{eqnarray*}
        \Xi^{{\rm partial}, k,\mathcal{P}}\left(z\right) &\sim& T_{\leq k-1}^{\mathcal{P}}\left(\rho_\mathcal{P}\right) \frac{(1-\frac{z}{\rho_\mathcal{P}})^{-1}}{2\tau_\mathcal{P}\rho_\mathcal{P}e_{\mathcal{P}-2}(\tau_\mathcal{P})} e^{\tau_{\mathcal{P}}} \Big(\sqrt{\frac{\rho_\mathcal{P} e_{\mathcal{P}-2}(\tau_\mathcal{P})}{2 \tau_\mathcal{P}}}(1-\frac{z}{\rho_\mathcal{P}})^{-\frac{1}{2}}+1\Big)\\
        &\sim& \frac{\tau_k^{\mathcal{P}} e^{\tau_{\mathcal{P}}}}{\sqrt{8\tau_\mathcal{P}^3\rho_\mathcal{P}e_{\mathcal{P}-2}(\tau_\mathcal{P})}}(1-\frac{z}{\rho_\mathcal{P}})^{-\frac{3}{2}},
      \end{eqnarray*}
      so the definition of $\rho_\mathcal{P}$ gives
      \begin{eqnarray*}
        \frac{\Xi^{k,\mathcal{P}}(z)}{\Xi^{{\rm partial},k,\mathcal{P}}(z)} &\sim& \frac{\frac{\tau_k^{\mathcal{P}}}{\tau_\mathcal{P}^2}\sqrt{\frac{e_\mathcal{P}(\tau_\mathcal{P})}{8 e_{\mathcal{P}-2}(\tau_\mathcal{P})}} \Big(1-\frac{z}{\rho_\mathcal{P}}\Big)^{-\frac{3}{2}}}{\frac{\tau_k^{\mathcal{P}} e^{\tau_{\mathcal{P}}}}{\sqrt{8\tau_\mathcal{P}^3\rho_\mathcal{P}e_{\mathcal{P}-2}(\tau_\mathcal{P})}}(1-\frac{z}{\rho_\mathcal{P}})^{-\frac{3}{2}}}\\
        &=& \frac{1}{e^{\tau_{\mathcal{P}}}}
      \end{eqnarray*}
      and the final claim.
    \end{proof}

    The first asymptotic in the next result is effectively an application of Theorem VII.2 on page 453 of \cite{FlajoletSedgewick}.
    The second asymptotic appears in Equation (7.8) in \cite{ArneyBender1982} and in Proposition VII.4 on page 464 of \cite{FlajoletSedgewick}.

    \begin{corollary}
      \label{F_Xi_coefficients}
      Let $k\geq 0$ and $\mathcal{P}\subseteq\mathbb{Z}_{\geq0}$, with $0\in\mathcal{P}$, $\mathcal{P}-2\neq\emptyset$, and $\gcd(\mathcal{P})=1$, and work over the ring $R=\mathbb{C}$.
      Then
      \begin{eqnarray*}
        \coeff{z^n}T^{\mathcal{P}}(z)&\sim& \sqrt{\frac{e_{\mathcal{P}}(\tau_\mathcal{P})}{2\pi e_{\mathcal{P}-2}(\tau_\mathcal{P})}}\frac{n^{-\frac{3}{2}}}{\rho_\mathcal{P}^{n}}\\
        \ \coeff{z^n}F^{\mathcal{P}}(z) &\sim& \frac{1}{\sqrt{2\pi \tau_\mathcal{P} \rho_\mathcal{P} e_{\mathcal{P}-2}(\tau_\mathcal{P})}}\frac{n^{-\frac{1}{2}}}{\rho_\mathcal{P}^{n}}\\
        \ \coeff{z^n}\Xi^{\CYC,\mathcal{P}}(z) &\sim& \frac{1}{2\tau_\mathcal{P} \rho_\mathcal{P}e_{\mathcal{P}-2}(\tau_\mathcal{P})}\frac{1}{\rho_\mathcal{P}^{n}}\\
        \ \coeff{z^n}\Xi^{k,\mathcal{P}}(z) &\sim& \frac{\tau_k^\mathcal{P}}{\tau_\mathcal{P}^2} \sqrt{\frac{e_{\mathcal{P}}(\tau_\mathcal{P})}{2\pi e_{\mathcal{P}-2}(\tau_\mathcal{P})}} \frac{n^{\frac{1}{2}}}{\rho_\mathcal{P}^{n}}\\
        \ \coeff{z^n}P^{\mathcal{P}}(z) &\sim&\ e^{\tau_\mathcal{P}}\coeff{z^n}F^{\mathcal{P}}(z)\\
        \ \coeff{z^n}\Xi^{{\rm partial}, k,\mathcal{P}}(z) &\sim&\ e^{\tau_\mathcal{P}}\coeff{z^n}\Xi^{k,\mathcal{P}}(z).
      \end{eqnarray*}
    \end{corollary}
    \begin{proof}
      We wish to apply sim-transfer Corollary \ref{sim_transfer} to each of
      \begin{eqnarray*}
        f_1(z) &=& \frac{\tau_\mathcal{P}-T^{\mathcal{P}}(\rho_\mathcal{P}z)}{\sqrt{\frac{2e_\mathcal{P}(\tau_\mathcal{P})}{e_{\mathcal{P}-2}(\tau_\mathcal{P})}}}\\
        f_2(z) &=& \sqrt{2\tau_\mathcal{P}\rho_\mathcal{P} e_{\mathcal{P}-2}(\tau_\mathcal{P})}F^{\mathcal{P}}(\rho_\mathcal{P}z)\\
        f_3(z) &=& 2\tau_\mathcal{P}\rho_\mathcal{P} e_{\mathcal{P}-2}(\tau_\mathcal{P})\Xi^{\CYC,\mathcal{P}}(\rho_\mathcal{P}z)\\
        f_4(z) &=& \frac{\tau_\mathcal{P}^2}{\tau_k^\mathcal{P}}\sqrt{\frac{8 e_{\mathcal{P}-2}(\tau_\mathcal{P})}{e_{\mathcal{P}}(\tau_\mathcal{P})}}\Xi^{k,\mathcal{P}}(\rho_\mathcal{P}z).
      \end{eqnarray*}
      The first hypothesis needed, that
      $$f_k(z)\sim (1-z)^{-\alpha_k},$$
      where
      \begin{eqnarray*}
        \alpha_1&=&-\frac{1}{2}\\
        \alpha_2&=&\frac{1}{2}\\
        \alpha_3 &=& 1\\
        \alpha_4&=&\frac{3}{2}
      \end{eqnarray*}
      is immediate from Corollaries \ref{T_asymptotic} and \ref{F_Xi_asymptotic}.
      The second and final hypothesis needed is that $f_k$ is $\Delta$-analytic, but this is immediate from the assumption that $\gcd(\mathcal{P})=1$ and the final statements in the same corollaries.

      Thus, the sim-transfer Corollary \ref{sim_transfer} says that
      \begin{eqnarray}
        \label{eq:sim_transfer}
        \coeff{z^n} f_k(z)&\sim& \frac{n^{\alpha_k-1}}{\Gamma(\alpha_k)}.
      \end{eqnarray}

      For $k=1$, Equation \ref{eq:sim_transfer} says that
      \begin{eqnarray*}
        \rho_\mathcal{P}^n\coeff{(\rho_\mathcal{P}z)^n}f_1(z)&\sim& \frac{n^{-\frac{1}{2}-1}}{\Gamma(-\frac{1}{2})}\\
        \ \coeff{(\rho_\mathcal{P}z)^n}\frac{\tau_\mathcal{P}-T^{\mathcal{P}}(\rho_\mathcal{P}z)}{\sqrt{\frac{2e_\mathcal{P}(\tau_\mathcal{P})}{e_{\mathcal{P}-2}(\tau_\mathcal{P})}}}&\sim& \rho_\mathcal{P}^{-n}\frac{n^{-\frac{3}{2}}}{-2\sqrt{\pi}}\\
        \ \coeff{z^n}\frac{T^{\mathcal{P}}(z)}{\sqrt{\frac{2e_\mathcal{P}(\tau_\mathcal{P})}{e_{\mathcal{P}-2}(\tau_\mathcal{P})}}}&\sim& \rho_\mathcal{P}^{-n}\frac{n^{-\frac{3}{2}}}{2\sqrt{\pi}},
      \end{eqnarray*}
      and the first result follows.

      For $k=2$, Equation \ref{eq:sim_transfer} says that
      \begin{eqnarray*}
        \rho_\mathcal{P}^n\coeff{(\rho_\mathcal{P}z)^n}f_2(z)&\sim& \frac{n^{\frac{1}{2}-1}}{\Gamma(\frac{1}{2})}\\
        \ \coeff{(\rho_\mathcal{P}z)^n}\sqrt{2\tau_\mathcal{P}\rho_\mathcal{P} e_{\mathcal{P}-2}(\tau_\mathcal{P})}F^{\mathcal{P}}(\rho_\mathcal{P}z) &\sim& \rho_\mathcal{P}^{-n}\frac{n^{-\frac{1}{2}}}{\sqrt{\pi}},
      \end{eqnarray*}
      and the second result follows.

      For $k=3$, Equation \ref{eq:sim_transfer} says that
      \begin{eqnarray*}
        \rho_\mathcal{P}^n\coeff{(\rho_\mathcal{P}z)^n}f_3(z)&\sim& \frac{n^{1-1}}{\Gamma(1)}\\
        \ \coeff{(\rho_\mathcal{P}z)^n}2\tau_\mathcal{P}\rho_\mathcal{P} e_{\mathcal{P}-2}(\tau_\mathcal{P})\Xi^{\CYC,\mathcal{P}}(\rho_\mathcal{P}z) &\sim& \rho_\mathcal{P}^{-n}\frac{1}{1},
      \end{eqnarray*}
      and the third result follows.

      For $k=4$, Equation \ref{eq:sim_transfer} says that
      \begin{eqnarray*}
        \rho_\mathcal{P}^n\coeff{(\rho_\mathcal{P}z)^n}f_4(z)&\sim& \frac{n^{\frac{3}{2}-1}}{\Gamma(\frac{3}{2})}\\
        \ \coeff{(\rho_\mathcal{P}z)^n}\frac{\tau_\mathcal{P}^2}{\tau_k^\mathcal{P}}\sqrt{\frac{8 e_{\mathcal{P}-2}(\tau_\mathcal{P})}{e_{\mathcal{P}}(\tau_\mathcal{P})}}\Xi^{k,\mathcal{P}}(\rho_\mathcal{P}z)&\sim& \rho_\mathcal{P}^{-n}\frac{n^{\frac{1}{2}}}{\frac{\sqrt{\pi}}{2}}\\
        \ \coeff{z^n}\Xi^{k,\mathcal{P}}(z)&\sim& \frac{\tau_k^\mathcal{P}}{\tau_\mathcal{P}^2} \sqrt{\frac{e_{\mathcal{P}}(\tau_\mathcal{P})}{8 e_{\mathcal{P}-2}(\tau_\mathcal{P})}} \frac{2n^{\frac{1}{2}}}{\sqrt{\pi}\rho_\mathcal{P}^{n}}\\
        &=& \frac{\tau_k^\mathcal{P}}{\tau_\mathcal{P}^2} \sqrt{\frac{e_{\mathcal{P}}(\tau_\mathcal{P})}{2\pi e_{\mathcal{P}-2}(\tau_\mathcal{P})}} \frac{n^{\frac{1}{2}}}{\rho_\mathcal{P}^{n}},
      \end{eqnarray*}
      and the fourth result follows.

      The final two results are immediate from the arguments above and Corollary \ref{F_Xi_asymptotic}.
    \end{proof}

    These estimates become accurate very quickly.
    Figure \ref{tree_function_034} compares the approximate counts of trees and functions subject to the constraint $\mathcal{P}=\{0,3,4\}$ to the corresponding exact counts taken from Figures \ref{exact_tree_counts} and \ref{exact_function_counts}.
    Note that the approximate curves are only visually distinguishable from the exact curves for small $n$, where the constraint is very restrictive.
    Constructing a similar graph for, say, $\mathcal{P}=\{0,1,2,3,4\}$, which is not so restrictive for small $n$, results in approximate curves that are completely obscured by the exact counts.

    \begin{figure}
      \begin{center}
        \includegraphics[width=7in]{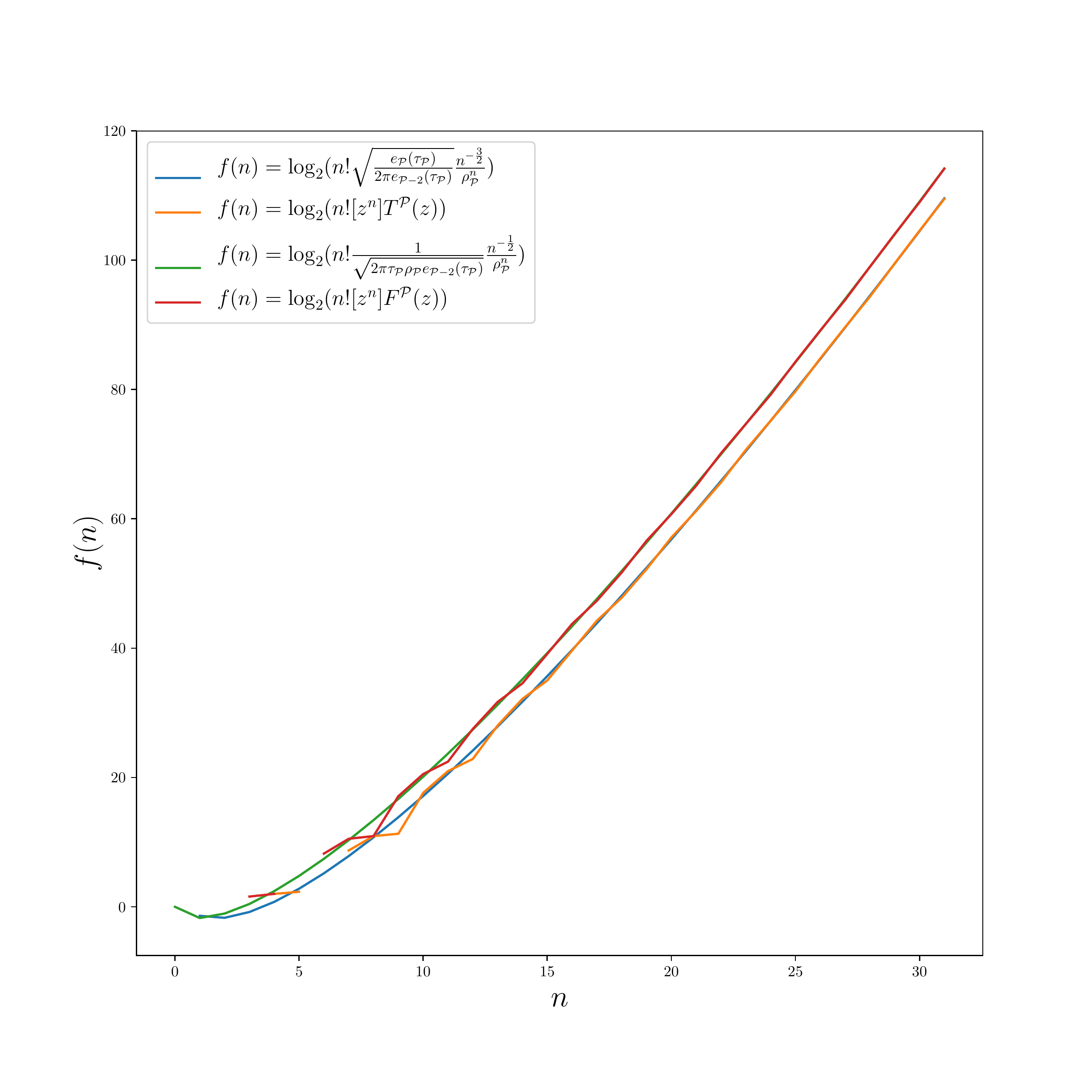}
      \end{center}
      \caption{Trees and Functions with Preimage Constraint $\mathcal{P}=\{0,3,4\}$}
      \label{tree_function_034}
    \end{figure}

    Even though Corollary \ref{F_Xi_coefficients} does not rigorously apply to $\mathcal{P}=\{0,4\}$, we can still compare the conclusions of that result to the exact counts.
    See Figure \ref{tree_function_04}.

    \begin{figure}
      \begin{center}
        \includegraphics[width=7in]{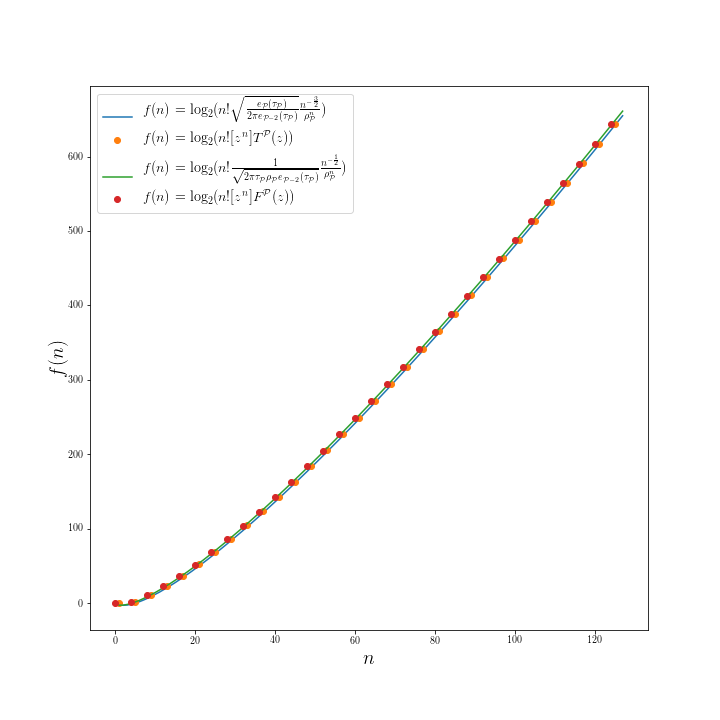}
      \end{center}
      \caption{Trees and Functions with Preimage Constraint $\mathcal{P}=\{0,4\}$}
      \label{tree_function_04}
    \end{figure}

    The next result appears in the entry for $c(\phi)$ in Table II on page 273 of \cite{ArneyBender1982}.

    \begin{corollary}
      \label{ave_cycle}
      Let $\mathcal{P}\subseteq\mathbb{Z}_{\geq0}$, with $0\in\mathcal{P}$, $\mathcal{P}-2\neq\emptyset$, and $\gcd(\mathcal{P})=1$, and work over the ring $R=\mathbb{C}$.
      Then the average number of cyclic points under $f$, where $f$ is subject to the preimage constraint $\mathcal{P}$, is
      $$\frac{\sum_{f\in\mathfrak{F}^\mathcal{P}, f:\set{n}\rightarrow\set{n}}|\{x\in \set{n}| f^k(x)=x \text{ for some }k\})|}{|\{f\in\mathfrak{F}^\mathcal{P}\mid f:\set{n}\rightarrow\set{n}\}|} \sim \sqrt{\frac{\pi}{2\tau_\mathcal{P}\rho_\mathcal{P} e_{\mathcal{P}-2}(\tau_\mathcal{P})}}n^{\frac{1}{2}}.$$
    \end{corollary}
    \begin{proof}
      Since the average number of cyclic points of a function on $n$ points that satisfies the preimage constraint $\mathcal{P}$ is $\frac{n!\coeff{z^n}\Xi^{\CYC,\mathcal{P}}(z)}{n!\coeff{z^n}F^{\mathcal{P}}(z)}$, the second and third results in Corollary \ref{F_Xi_coefficients} yield that this is asymptotically
      \begin{eqnarray*}
        \frac{n!\coeff{z^n}\Xi^{\CYC,\mathcal{P}}(z)}{n!\coeff{z^n}F^{\mathcal{P}}(z)}&\sim& \frac{\frac{1}{2\tau_\mathcal{P} \rho_\mathcal{P} e_{\mathcal{P}-2}(\tau_\mathcal{P})}\frac{1}{\rho_\mathcal{P}^{n}}}{\frac{1}{\sqrt{2\pi \tau_\mathcal{P} \rho_\mathcal{P} e_{\mathcal{P}-2}(\tau_\mathcal{P})}}\frac{n^{-\frac{1}{2}}}{\rho_\mathcal{P}^{n}}}\\
        &=& \frac{\sqrt{2\pi \tau_\mathcal{P} \rho_\mathcal{P} e_{\mathcal{P}-2}(\tau_\mathcal{P})}}{2\tau_\mathcal{P} \rho_\mathcal{P}e_{\mathcal{P}-2}(\tau_\mathcal{P})}n^{\frac{1}{2}}\\
        &=& \sqrt{\frac{\pi}{2\tau_\mathcal{P}\rho_\mathcal{P} e_{\mathcal{P}-2}(\tau_\mathcal{P})}}n^{\frac{1}{2}},
      \end{eqnarray*}
    \end{proof}

    The next result is the main goal of this paper.
    For the impatient reader who jumped straight here from the beginning of the document, we recall a few definitions, without context or justification.
    Equations \ref{eq:eP_def} and \ref{eq:Pminus} defined
    \begin{eqnarray*}
      e_\mathcal{P}(z)&=&\sum_{n\in\mathcal{P}}\frac{z^n}{n!}\\
      \mathcal{P}-i&=&\{n-i\mid n\in\mathcal{P}\}\setminus\mathbb{Z}_{<0}.
    \end{eqnarray*}
    The paragraph containing Equation \ref{eq:rho_def} defines $\tau_\mathcal{P}$ to be the unique $\tau\in\mathbb{R}_{> 0}$ such that $$e_\mathcal{P}(\tau)-\tau e_\mathcal{P}'(\tau)=0$$ and defines
    $$\rho_\mathcal{P}=\frac{\tau_\mathcal{P}}{e_\mathcal{P}(\tau_\mathcal{P})}.$$
    Lemma \ref{tree_recursion_formula}, Theorem \ref{when_u_equals_1_tree}, and Equation \ref{tau_k_def} give that $\tau_{k}^\mathcal{P}$ can be computed recursively by
    \begin{eqnarray*}
      \tau_0^\mathcal{P} &=&0\\
      \tau_{k+1}^\mathcal{P} &=&\rho_\mathcal{P} e_\mathcal{P}(\tau_k^\mathcal{P})
    \end{eqnarray*}
    for $k\geq 0$.

    \begin{theorem}
      \label{kth_image}
      Let $k\geq 0$ and $\mathcal{P}\subseteq\mathbb{Z}_{\geq0}$, with $0\in\mathcal{P}$, $\mathcal{P}-2\neq\emptyset$, and $\gcd(\mathcal{P})=1$, and work over the ring $R=\mathbb{C}$.
      Then the average size of a $k$th image of a function, where the function is subject to the preimage constraint $\mathcal{P}$, is
      $$\frac{\sum_{f\in\mathfrak{F}^\mathcal{P}, f:\set{n}\rightarrow\set{n}}|f^k(\set{n})|}{|\{f\in\mathfrak{F}^\mathcal{P}\mid f:\set{n}\rightarrow\set{n}\}|} \sim n\Big(1-\frac{\tau_k^\mathcal{P}}{\tau_\mathcal{P}}\Big).$$
      If one instead averages over partial functions subject to $\mathcal{P}$, one gets the same result, namely,
      $$\frac{\sum_{f\in\mathfrak{P}^\mathcal{P}, f:\set{n}\rightarrow\set{n}}|f^k(\set{n})|}{|\{f\in\mathfrak{P}^\mathcal{P}\mid f:\set{n}\rightarrow\set{n}\}|} \sim n\Big(1-\frac{\tau_k^\mathcal{P}}{\tau_\mathcal{P}}\Big).$$
    \end{theorem}
    \begin{proof}
      Recall that $\Xi^{k,\mathcal{P}}$ was defined in such a way that
      $$n!\coeff{z^n}\Xi^{k,\mathcal{P}}(z)=\sum_{f\in\mathfrak{F}^\mathcal{P}}|\set{n}\setminus f^k(\set{n})|,$$
      so
      \begin{eqnarray}
        \frac{\sum_{f\in\mathfrak{F}^\mathcal{P}, f:\set{n}\rightarrow\set{n}}|\set{n}\setminus f^k(\set{n})|}{|\{f\in\mathfrak{F}^\mathcal{P}\mid f:\set{n}\rightarrow\set{n}\}|} &=& \nonumber \frac{n!\coeff{z^n}\Xi^{k,\mathcal{P}}(z)}{n!\coeff{z^n}F^{\mathcal{P}}(z)}\\
        n-\frac{\sum_{f\in\mathfrak{F}^\mathcal{P}, f:\set{n}\rightarrow\set{n}}|f^k(\set{n})|}{|\{f\in\mathfrak{F}^\mathcal{P}\mid f:\set{n}\rightarrow\set{n}\}|}  &=& \label{eq:Xi_F_ratio}\frac{\coeff{z^n}\Xi^{k,\mathcal{P}}(z)}{\coeff{z^n}F^{\mathcal{P}}(z)}\\
      \end{eqnarray}
      By Corollary \ref{F_Xi_coefficients}, this is
      \begin{eqnarray*}
        n-\frac{\sum_{f\in\mathfrak{F}^\mathcal{P}, f:\set{n}\rightarrow\set{n}}|f^k(\set{n})|}{|\{f\in\mathfrak{F}^\mathcal{P}\mid f:\set{n}\rightarrow\set{n}\}|} &\sim& \frac{\frac{\tau_k^\mathcal{P}}{\tau_\mathcal{P}^2} \sqrt{\frac{e_{\mathcal{P}}(\tau_\mathcal{P})}{2\pi e_{\mathcal{P}-2}(\tau_\mathcal{P})}} \frac{n^{\frac{1}{2}}}{\rho_\mathcal{P}^{n}}}{\frac{1}{\sqrt{2\pi \tau_\mathcal{P} \rho_\mathcal{P} e_{\mathcal{P}-2}(\tau_\mathcal{P})}}\frac{n^{-\frac{1}{2}}}{\rho_\mathcal{P}^{n}}}\\
        &=& \frac{\tau_k^\mathcal{P}}{\tau_\mathcal{P}^2}\sqrt{\frac{e_\mathcal{P}(\tau_\mathcal{P})2\pi \tau_\mathcal{P} \rho_\mathcal{P} e_{\mathcal{P}-2}(\tau_\mathcal{P})}{2\pi e_{\mathcal{P}-2}(\tau_\mathcal{P})}}n\\
        &=& \frac{\tau_k^\mathcal{P}}{\tau_\mathcal{P}^2}\sqrt{ \tau_\mathcal{P} \rho_\mathcal{P} e_\mathcal{P}(\tau_\mathcal{P})}n\\
        &=& \frac{\tau_k^\mathcal{P}}{\tau_\mathcal{P}}\sqrt{\frac{\rho_\mathcal{P} e_\mathcal{P}(\tau_\mathcal{P})}{\tau_\mathcal{P}}}n.
      \end{eqnarray*}
      Recalling the definition of $\rho_\mathcal{P}$ gives
      \begin{eqnarray*}
        n-\frac{\sum_{f\in\mathfrak{F}^\mathcal{P}, f:\set{n}\rightarrow\set{n}}|f^k(\set{n})|}{|\{f\in\mathfrak{F}^\mathcal{P}\mid f:\set{n}\rightarrow\set{n}\}|} &\sim& \frac{\tau_k^\mathcal{P}}{\tau_\mathcal{P}}n,
      \end{eqnarray*}
      from which the first result follows.

      The partial function version of Equation \ref{eq:Xi_F_ratio} is
      \begin{eqnarray*}
        n-\frac{\sum_{f\in\mathfrak{P}^\mathcal{P}, f:\set{n}\rightarrow\set{n}}|f^k(\set{n})|}{|\{f\in\mathfrak{P}^\mathcal{P}\mid f:\set{n}\rightarrow\set{n}\}|} &=& \frac{\coeff{z^n}\Xi^{{\rm partial}, k,\mathcal{P}}(z)}{\coeff{z^n}P^{\mathcal{P}}(z)}.
      \end{eqnarray*}
      By Corollary \ref{F_Xi_coefficients}, the expression on the right is exactly $\frac{\coeff{z^n}\Xi^{k,\mathcal{P}}(z)}{\coeff{z^n}F^{\mathcal{P}}(z)}$ and so the calculations above apply verbatim.
    \end{proof}

    Note that taking $k=1$ in Theorem \ref{kth_image} and subtracting the result from $n$ gives that, on average (and when the preimage constraints satisfy the hypotheses), there are asymptotically $\frac{\rho_\mathcal{P}}{\tau_\mathcal{P}}n$ points that have no preimages; this observation is Equation (7.9) in \cite{ArneyBender1982}.

    We close this section by exploring the asymptotics of the calculations in the paragraph containing Equation \ref{eq:F_gen_seq}.
    Take $\mathcal{P}=\mathbb{Z}_{\geq 0}$, so $\tau_{\mathbb{Z}_{\geq 0}}=1$ is the unique positive solution to $e^\tau-\tau e^\tau=0$ and $\rho_{\mathbb{Z}_{\geq 0}}=\frac{1}{e^{\tau_{\mathbb{Z}_{\geq 0}}}}=e^{-1}$.
    Then $\sqrt{\frac{2e_\mathcal{P}(\tau_\mathcal{P})}{e_{\mathcal{P}-2}(\tau_\mathcal{P})}}=\sqrt{\frac{2e^1}{e^1}}=\sqrt{2}$, and Corollary \ref{T_asymptotic} says
    $$T^{\mathbb{Z}_{\geq 0}}(z) \sim 1-\sqrt{2}\big(1-\frac{z}{e^{-1}}\big)^{\frac{1}{2}}=1-\sqrt{2}\left(1-e z\right)^{\frac{1}{2}},$$
    which is consistent with Proposition 1 of \cite{FlajoletOdlyzko_nolink}.
    Since there are $n^n$ functions from $\set{n}$ to itself, $F^{\mathbb{Z}_{\geq 0}}(z)=\sum_{n=0}^\infty \frac{n^n}{n!} z^n$.
    But Corollary \ref{F_Xi_coefficients} ensures that
    \begin{eqnarray*}
      \frac{n^n}{n!} &=& \coeff{z^n}F^{\mathbb{Z}_{\geq 0}}(z)\\
      &\sim& \frac{1}{\sqrt{2\pi e^{-1} e^{1}}}\frac{n^{-\frac{1}{2}}}{e^{-n}}\\
      &=& \frac{n^n}{\sqrt{2\pi n}\big(\frac{n}{e}\big)^n}\\
      \sqrt{2\pi n}\Big(\frac{n}{e}\Big)^n &=& n!;
    \end{eqnarray*}
    in particular, we have repeated a derivation of Stirling's approximation given in Section 2 of \cite{FlajoletOdlyzko_nolink}.
    Finally, Theorem \ref{kth_image} yields
    \begin{eqnarray}
      \frac{1}{n^n}\sum_{f:\set{n}\rightarrow\set{n}}|f^k(\set{n})|  &\sim& \label{eq:FO_iterates}
      n\Big(1-T_{\leq k-1}^{\mathbb{Z}_{\geq0}}(e^{-1})\Big);
    \end{eqnarray}
    appealing to the recursion in Lemma \ref{tree_recursion_formula}, we have reproduced result (v) in \cite{FlajoletOdlyzko_nolink}'s Direct Parameters Theorem 2.
  \end{subsection}
\end{section}

\begin{section}{Coalescence}
  \label{coalescence}

  The motivation for the preceding work is to understand coalescence of random functions subject to preimage constraints.
  Theorem \ref{kth_image} gives an answer that addresses how fast the coalescence occurs, and Corollary \ref{ave_cycle} gives an answer that addresses the ultimate size after the coalescence has stabilized.

  It is tempting, but fallacious, to try to recover the second answer from the first.
  For a concrete function, the iterated image will eventually stabilize, so one might take the limit as $k$ approaches infinity in Theorem \ref{kth_image}.
  By Lemma \ref{tau_k_converge}, the result is 0, and one erroneously concludes that the number of cyclic points is asymptotically 0.
  An explanation as to why this does not match Corollary \ref{ave_cycle} is that the corollary, in effect, first takes the limit as $k$ goes to infinity, and then looks at the asymptotics of $n$; in other words, one cannot swap the order one takes the limits.

  Recall that Theorem \ref{kth_image} can be interpreted as an estimate for the size of the $k$th iterate of a function $f:\set{n}\rightarrow \set{n}$ that is subject to the preimage constraint $\mathcal{P}\subseteq\mathbb{Z}_{\geq 0}$; more explicitly,
  $$|f^k(n)|\approx n\Big(1-\frac{\tau_k^\mathcal{P}}{\tau_\mathcal{P}}\Big).$$
  To demonstrate this concretely, Figure \ref{tau_k} plots the sequence $(1-\frac{\tau_k^\mathcal{P}}{\tau_\mathcal{P}})$, on a $\log_2$ scale, for all $\mathcal{P}\subseteq\{0,1,2,3,4\}$ for which Theorem \ref{kth_image} applies, as well as for the unconstrained $\mathcal{P}=\mathbb{Z}_{\geq 0}$.

  \begin{figure}
    \begin{center}
      \includegraphics[width=7in]{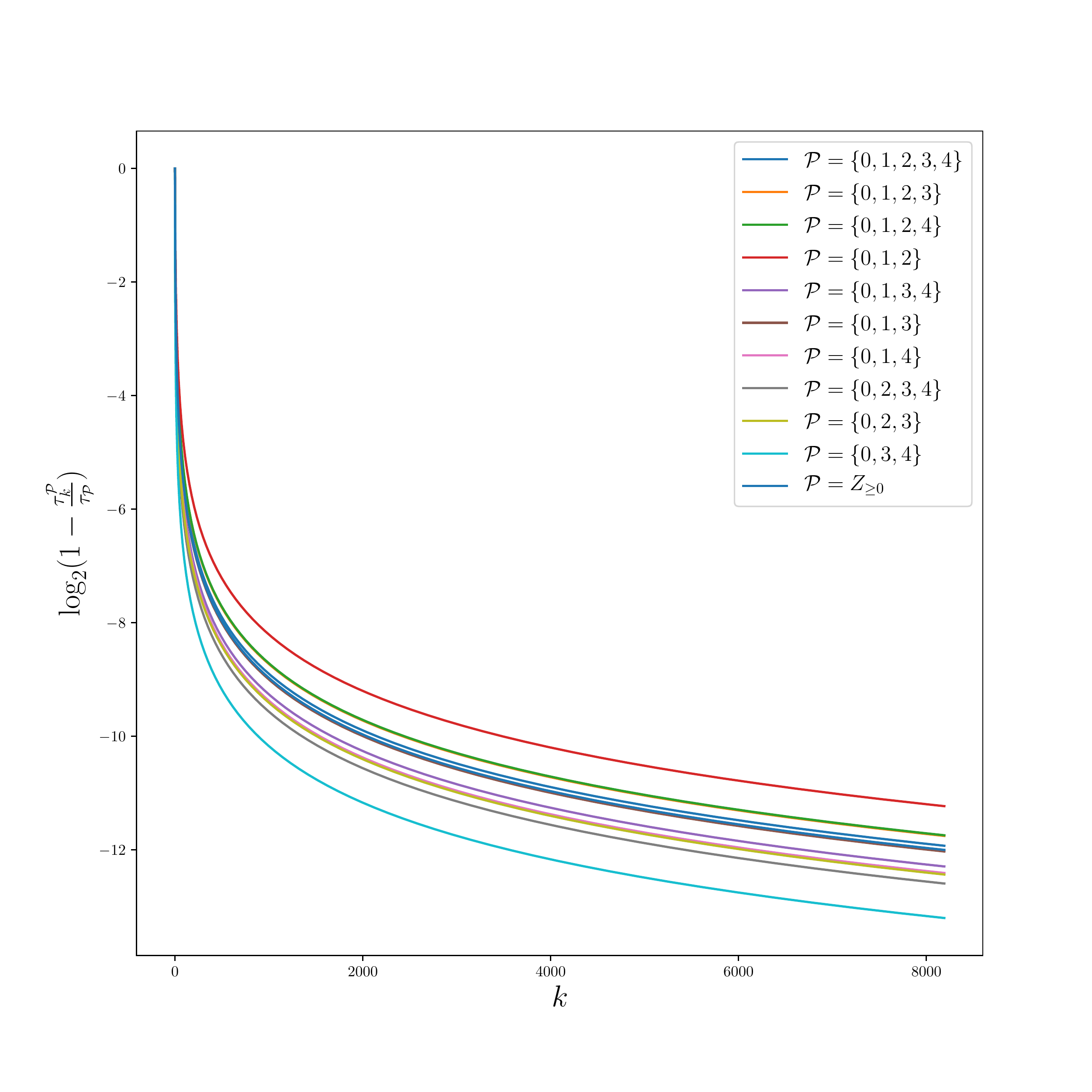}
    \end{center}
    \caption{$\log_2$ of Ratio of Points in $k$th Image to All Points}
    \label{tau_k}
  \end{figure}

  Figure \ref{from_unconstrained} subtracts $\log_2(1-\frac{\tau_k^{\mathbb{Z}_\geq 0}}{\tau_{\mathbb{Z}_{\geq 0}}})$ from each $\log_2(1-\frac{\tau_k^\mathcal{P}}{\tau_\mathcal{P}})$ appearing in Figure \ref{tau_k}; in other words, it views the unconstrained case as the baseline from which all other cases' coalescence is judged.
  Since the resulting functions become constant very quickly, one may use $\log_2(1-\frac{\tau_k^\mathcal{P}}{\tau_\mathcal{P}})-\log_2(1-\frac{\tau_k^{\mathbb{Z}_\geq 0}}{\tau_{\mathbb{Z}_{\geq 0}}})$ for some fixed $k$ as a proxy metric for how fast a random function with preimage constraint $\mathcal{P}$ coalesces.

  \begin{figure}
    \begin{center}
      \includegraphics[width=7in]{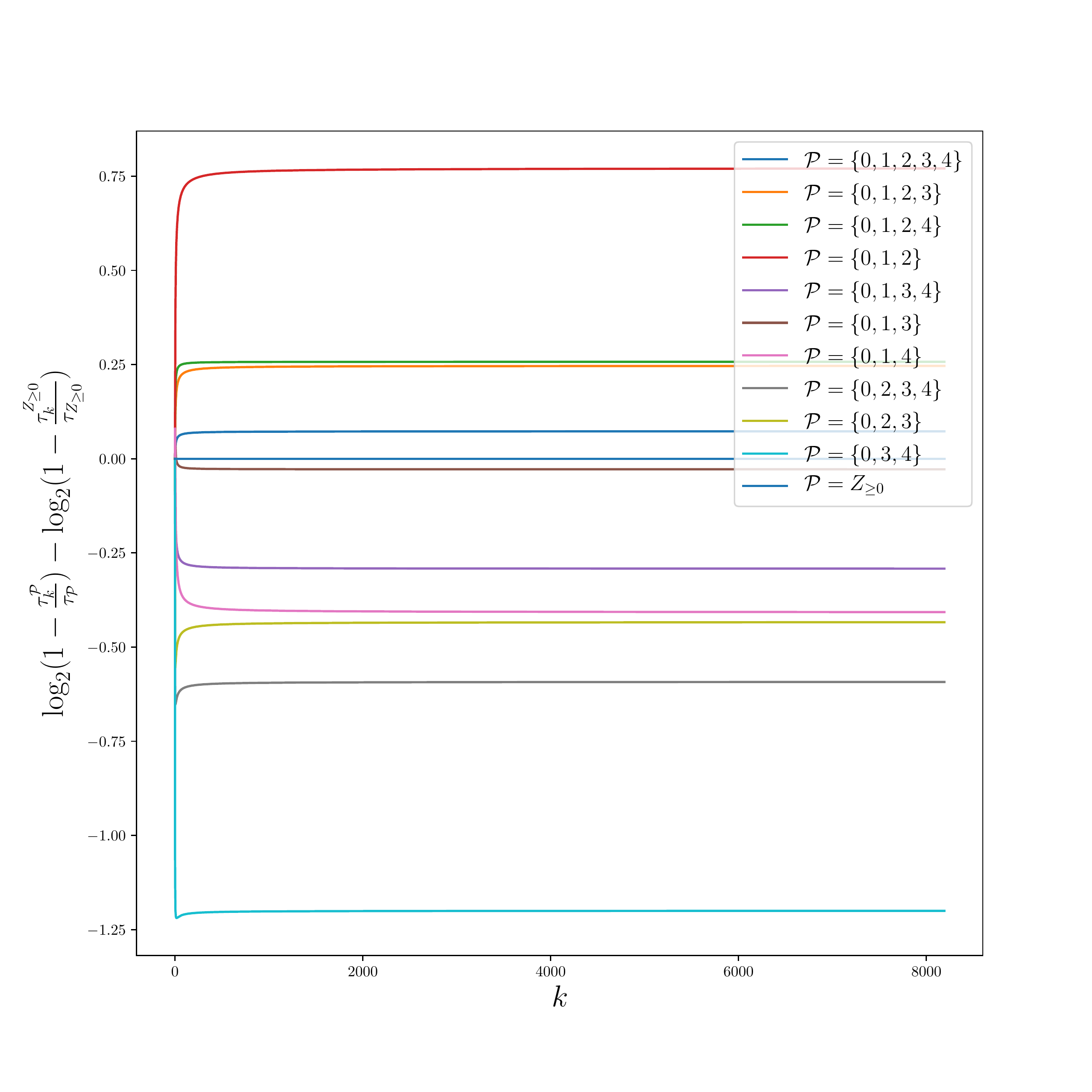}
    \end{center}
    \caption{Coalescence}
    \label{from_unconstrained}
  \end{figure}

  Recall that Corollary \ref{ave_cycle} says that, under certain conditions, random functions subject to preimage constraints asymptotically have $c \sqrt{n}$ cyclic points, where $c=\sqrt{\frac{\pi}{2\tau_\mathcal{P}\rho_\mathcal{P} e_{\mathcal{P}-2}(\tau_\mathcal{P})}}$ is a constant that depends on the constraint.
  For the sake of completeness, Figure \ref{coalesce_scatter} compares the coalescence measure of Figure \ref{from_unconstrained} to the $c$s.
  One immediately notices that the $c$s are correlated to the rate of coalescence.
  In other words, the a priori different questions of how fast a function coalesces and where that coalescence stabilizes actually have related answers.

  \begin{figure}
    \begin{center}
      \includegraphics[width=7in]{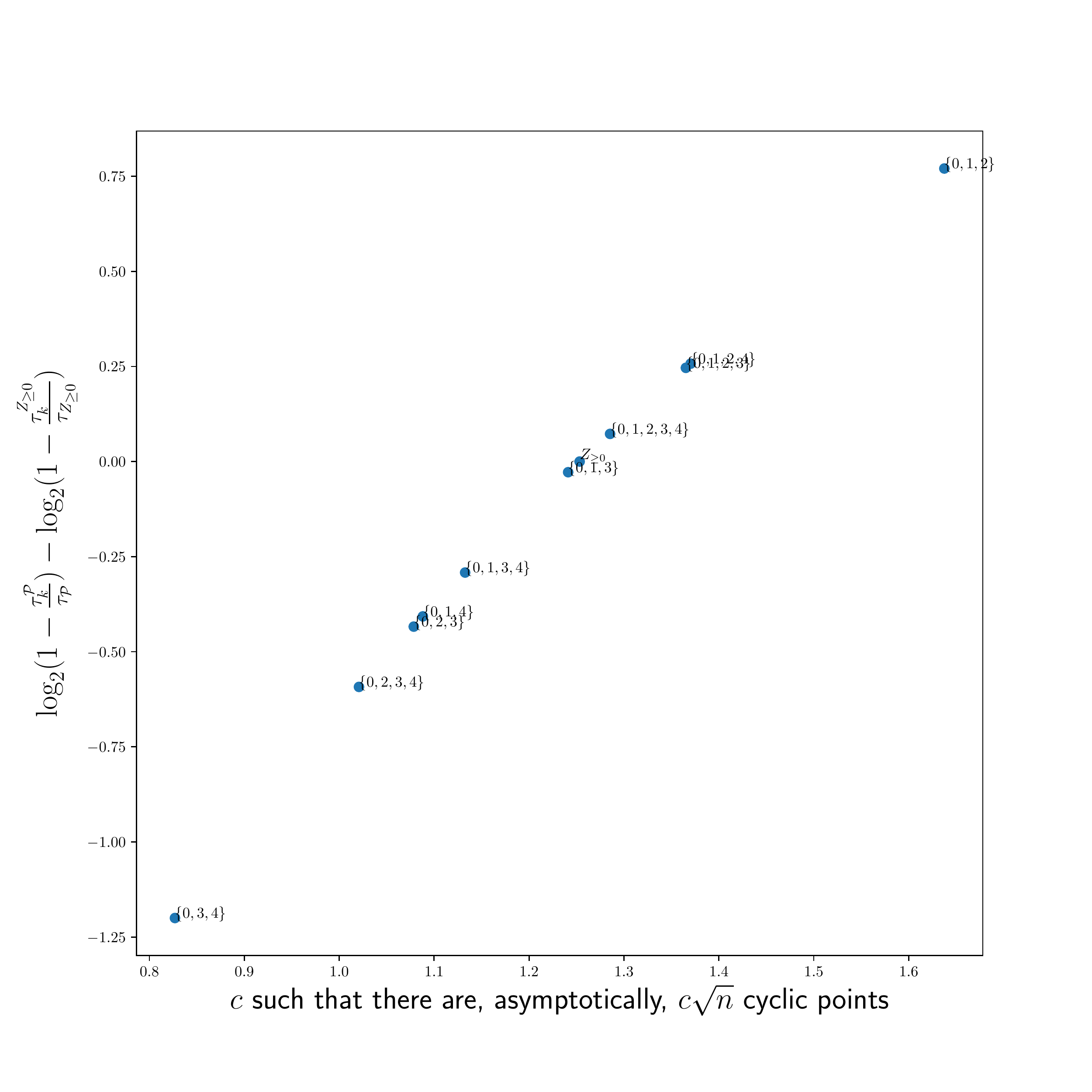}
    \end{center}
    \caption{Cyclic Points Versus Coalescence}
    \label{coalesce_scatter}
  \end{figure}
\end{section}

\begin{section}{Final Thoughts}
  Cheyne Homberger has pointed out that Corollaries \ref{tree_coefficient} and \ref{function_coefficients} provide useful tools for explicit enumeration.
  A superficial search of the Online Encyclopedia of Integer Sequences for the number of functions that satisfy specific constraints $\mathcal{P}$ yielded some hits and some misses; regardless, it does not appear that these sequences are being approached in a unified manner.

  The biggest deficiency in Theorem \ref{kth_image} is the requirement that $\gcd(\mathcal{P})=1$.
  \cite{FlajoletSedgewick} addresses how to drop this aperiodic condition to extend results for $T_{\leq k}^\mathcal{P}$, and one might be able to push those through to results about $\Xi^{k,\mathcal{P}}$ and $F^\mathcal{P}$.
  \cite{ArneyBender1982} includes a more cavalier claim that one will get the same results as with the aperiodic case, since moving to the periodic situation leads to the same constant appearing in the coefficients of all relevant generating functions, which in turn cancel out when averaging over $\mathfrak{F}^\mathcal{P}$ (that is, the same constant appears in both the numerator and denominator).
  While a strict interpretation of these claims is false, since the periodic condition results in many 0 coefficients in the resulting generating functions, Figure \ref{tree_function_04} certainly supports it.
\end{section}

%%%%%%%%%%%%%%%%%%%%%%%%%%%%%%%%%%%%%%%%%%%%%%%%%%%%%%%%%%%%%%%%%%%%%%%%%%%%%%%%%%%%%%%%%%%%%%%%%%%
% Acknowledgements
%%%%%%%%%%%%%%%%%%%%%%%%%%%%%%%%%%%%%%%%%%%%%%%%%%%%%%%%%%%%%%%%%%%%%%%%%%%%%%%%%%%%%%%%%%%%%%%%%%%

\begin{section}*{Acknowledgments}

  The author is grateful to a number of people, most especially to Cheyne Homberger for several illuminating discussions and to Art Drisko and Art Pittenger for careful comments on an earlier draft.
\end{section}

\bibliographystyle{apalike}
\bibliography{preimage}

\end{document}